\documentclass[]{interact}

\usepackage{epstopdf}% To incorporate .eps illustrations using PDFLaTeX, etc.
\usepackage[caption=false]{subfig}% Support for small, `sub' figures and tables

\usepackage[numbers,sort&compress]{natbib}% Citation support using natbib.sty
\bibpunct[, ]{[}{]}{,}{n}{,}{,}% Citation support using natbib.sty
\theoremstyle{plain} % Theorem-like structures provided by amsthm.sty
\newtheorem{theorem}{Theorem}[section]
\newtheorem{lemma}[theorem]{Lemma}
\newtheorem{corollary}[theorem]{Corollary}
\newtheorem{proposition}[theorem]{Proposition}

\theoremstyle{definition}

\theoremstyle{remark}
\newtheorem{remark}{Remark}

\usepackage{amsmath,amssymb}
\usepackage{amsfonts}
\usepackage{amsthm}
\usepackage{newlfont}
\usepackage{epsfig}
\usepackage{graphics}
\usepackage{graphicx}
\usepackage{comment}
\usepackage{color}
\usepackage{hyperref}
\usepackage[pagewise]{lineno}
\newcommand{\R}{\mathbb{R}}

\begin{document}
\thispagestyle{empty} \setcounter{page}{1}

%line numbers%%%%%%%%
%\setpagewiselinenumbers
%\modulolinenumbers[5]
%\linenumbers
%%%%%%%%%%%%%%%%%%%%%

%\noindent

%\begin{center}
%{\Large\bf
%Chaotic dynamics in a
%simple predator-prey model with discrete delay
%}

%\setcounter{page}{\pageref{0}}
%\setcounter{page}{1}

\title{Chaotic dynamics in a simple predator-prey model with discrete delay}

% {\raggedleft\today\par}

%\setcounter{footnote}{0}

%\vskip.20in

\author{
\name{Guihong Fan\textsuperscript{a} and Gail S. K.
	Wolkowicz\textsuperscript{b}\thanks{CONTACT G.~S.~K. Wolkowicz.
	Email: wolkowic@mcmaster.ca}}
\affil{\textsuperscript{a}Department of Mathematics, Columbus State
	University, Columbus, Georgia 31907;
	\textsuperscript{b}Department of Mathematics and
	Statistics, McMaster University, Hamilton, Ontario,
	Canada L8S 4K1}
}

\maketitle

\begin{abstract}
	A discrete delay is included to model the time between the capture of the prey and its conversion to viable biomass in  the simplest classical Gause type predator-prey model that has  equilibrium dynamics without delay. As the delay increases from zero, the coexistence equilibrium undergoes a supercritical Hopf bifurcation, two saddle-node bifurcations of limit cycles, and a cascade of period doublings, eventually leading to chaos. The resulting periodic orbits and the  strange attractor resemble their counterparts for the Mackey-Glass equation.  Due to the global stability of the system without delay, this complicated dynamics can be solely attributed to the introduction of the delay.  Since many   models include predator-prey like interactions as submodels, this study emphasizes the importance of understanding the implications of overlooking delay in such models  on the reliability of the model-based predictions, especially since temperature  is known to have an effect on the length of certain delays.
\end{abstract}



\section{Introduction}\label{intro}
\noindent
A Gause type predator-prey model with response
function $f(x)$ is given by
\begin{equation}\label{ppnodelay}
\left\{
\begin{aligned}
&\dot x(t)=rx(t)\bigg(1-\frac{x(t)}{K}\bigg)-y(t) f(x(t)),\\
&\dot y(t)=-sy(t)+Y y(t) f(x(t)),
\end{aligned}\right.
\end{equation}
where $x(t)$ denotes the density of the prey population and $y(t)$
the density of predators. Parameters $r,\ K,\ s$, and $Y$ are
positive constants denoting the intrinsic growth rate and the
carrying capacity of the prey, the death rate of the predator in the
absence of prey, and the growth yield constant for the conversion of prey
to viable predator density, respectively.

If  $f(x)$ is
of Holling type I form in model (\ref{ppnodelay}) (i.e. $f(x)=mx$ where
$m$  is a positive constant denoting the maximal growth rate of
the predator), it is well-known (see e.g.
 \cite{Bazykin,Freedman1980}) that  either
the predator population approaches extinction and the prey population
approaches its carrying capacity, or the predator population and the
prey population coexist and their density approaches a positive
equilibrium. Hence, for
all choices of the parameters, all solutions of this system approach a
globally asymptotically stable equilibrium, and so any nontrivial
oscillatory behaviour that arises due to the introduction of delay in
the model
can be attributed solely to the delay. For this reason, we choose
Holling type I response functions instead
of the more realistic Holling type II form, since the Holling type II
form  results in a model that gives rise to  nontrivial period solutions  without delay (see Rosenzweig \cite{Rosenzweig1971}).
One would also expect that any  exotic dynamics that the model with  Holling
type I form admits due to the introduction of delay would  be shared by
the model with Holling type II form.   Li et al. \cite{MR3180725} studied this model with the Holling type II
response function of Monod form and showed that stability switches caused by
varying the time delay are accompanied by bounded global Hopf branches, and they proved that when multiple Hopf branches exist, they are nested and the overlap produces coexistence of two or possibly more stable limit cycles.
However, they did not go on to discover the even richer  dynamics that
our analysis suggests exists in that case.

Incorporating a time delay in (\ref{ppnodelay}) to model the time between the
capture of the prey by the predator and its conversion to viable
predator biomass, in the case of Holling type I functional response,
$f(x)=mx, \ m>0$, we  obtain the following system:
\begin{equation}
\left\{\begin{aligned}
&\dot x(t)=rx(t)\left(1-\frac{x(t)}{K}\right)-m y(t)x(t),\\
&\dot y(t)=-sy(t)+Ye^{-s\tau}my(t-\tau)x(t-\tau). \label{ppreydelay}
%&\mbox{with initial data for} \ x(t) \ \mbox{and} \ y(t) \ \mbox{in} \ C([-\tau,0],\mbox{int}\R_+^2)
\end{aligned}\right.
\end{equation}
The term $e^{-s\tau}y(t-\tau)$ represents those predators that
 survive the $\tau \geq 0$ units of time required to process the
 prey captured at time  $t-\tau$  in the past. Thus, we have incorporated the delay in the  growth term of the predator equation in a manner that is
consistent with  its  decline rate given by the model, as described in Arino,
Wang, and Wolkowicz \cite{Gail2006Alter}.

We define
$\R_+\equiv\{x\in\R : \, x\geqslant 0\}$,
int$\R_+\equiv\{x\in\R : \, x>0 \}$.   We denote by
 $C([-\tau,0],  \R_+)$,  the Banach space of continuous functions   from
 the interval $[-\tau,0]$ into  $\R_+$, equipped with the uniform norm.
 We assume initial data for model (\ref{ppreydelay}) is taken from
 \begin{equation}\label{eq:X}
 X= C([-\tau,0],  \R_+) \times  C([-\tau,0],  \R_+).
\end{equation}

%%%%%%%%%%%%%%%%%%%%%%%%%%%%%%%%%%%%%%%%%%%%%%%%%
Model (\ref{ppreydelay}) also has other interpretations.
Gourley and Kuang \cite{Kuang2004} studied
a stage-structured predator-prey model
in which they included an equation for the juvenile
predators and assumed a constant maturation time
delay, i.e.,  they assumed that the juvenile predators
take a fixed  time to mature.
Using the approach developed in Beretta and Kuang
\cite{Kuang2002},
the authors considered the
possibility of stability switches, and concluded that
there is a range of the parameter modeling the time delay for which
there are periodic solutions.
 If the juveniles in their model suffer the
same mortality rate as adult predators, their model decouples and
yields model (\ref{ppreydelay}).
Forde \cite{forde2005delay} also considered this model and conjectured
that there are periodic orbits whenever the interior equilibrium exists
and is unstable. He also noted that if the interior equilibrium
exists and is asymptotically stable without delay, then for small
delays it remains globally asymptotically stable. We show that whenever the
interior equilibrium exists and is unstable, the system is uniformly
persistent. Gourley and Kuang
\cite{Kuang2004} had already showed that a Hopf bifurcation eventually occurs
if the delay is increased, destabilizing this equilibrium and  giving
birth to a nontrivial periodic solution. Forde  \cite{forde2005delay}
 left as an open question whether more than one periodic orbit is
 possible and provided a numerical example suggesting chaos is possible
 but does not consider the route to chaos.
 We give numerical evidence that there is a range of parameters for which
 two stable periodic orbits and an unstable periodic orbit all exist and
we show a period-doubling route to chaos followed by a
 period-halfing route back to stability of the interior equilibrium.

 Cooke, Elderkin, and Huang \cite{Cooke2006} considered a model
similar to the one in Gourley and Kuang \cite{Kuang2004}, and
obtained results concerning Hopf bifurcation of a scaled version.
The scaling they used  eliminated the parameter modelling the time
delay, the parameter that we focus on and use as a bifurcation
parameter.  This simplified their analysis,
since then, unlike in our case, the components of the coexistence
equilibrium are independent of the time delay.

In this manuscript we show that the introduction of
time delay cannot only destabilize the  globally asymptotically stable
coexistence equilibrium of model (\ref{ppreydelay}),  it can also be responsible for exotic dynamics for
intermediate values of the delay as well as the  eventual   disappearance of
the coexistence equilibrium with the extinction of the predator
for large enough delays.
Although chaotic dynamics has been observed in other models of
predator-prey interactions, the other models either require at least
three trophic levels, or the  response functions  are not as simple and
so the models admit oscillatory behavior even in the absence of delay,
or the other models incorporate the delay  in such a way  that   the
predators   still contribute to
 population growth even if the time required
 to process the prey is longer than the life-span of the predator
 (i.e., the factor $e^{-s\tau}$ is missing in the $\dot{y}$ equation),
 or  the delay is used to model different mechanisms (see e.g.,
\cite{Ginoux2005,Hastings1991,Morozov2004,JWang2012}). %Nakaoka2006,
The observation that the resulting strange attractor   resembles the
strange attractor for the Mackey-Glass equation  \cite{Mackey:2009} is also new.

This paper is organized as follows.  In section~\ref{basic}, we scale the
model and show that it is well-posed. In section~\ref{analysis}, we
consider the existence and stability of equilibria.
If parameters are set so that it is possible for the predator to survive
when there is no delay, it is well-known that the
equilibrium at which both the prey and the predator survive  is globally asymptotically stable with respect to positive
initial conditions   (i.e. solutions
approach this equilibrium for any choice of positive initial data).
In the case of delay,
  the components of this coexistence equilibrium depend on the
delay.
We prove that for positive delay,  when this equilibrium exists, both the predator and the prey
populations persist uniformly. However, a sufficiently long delay
results in   the disappearance of this equilibrium, resulting in  the extinction
of the predator and convergence to a globally asymptotically stable
equilibrium with the prey at carrying capacity.
We give criteria which when satisfied imply that there are at least two Hopf
bifurcations that occur before the extinction of the predator, resulting
in sustained oscillatory behaviour for intermediate values of the delay. Finally, in
section~\ref{sec:example},  by means of
time series, time delay embeddings, and orbit (bifurcation) diagrams we show
that there are saddle-node bifurcations of limit cycles resulting in
bistability as well as sequences of period doubling bifurcations leading to
chaos,
with a strange attractor resembling the strange attractor for the
Mackey-Glass equation \cite{Mackey:2009}.
We conclude with a brief discussion.

\section{Scaling and basic properties of solutions}\label{basic}
In order to simplify the analysis, we introduce the following change of
variables:
\begin{equation}
\begin{aligned}
&\breve{t}=rt,\qquad \breve{x}(\breve{t})=x(t)/K,\qquad
\breve{y}(\breve{t})=m y(t)/r, \\
&\breve{\tau}=r\tau,\qquad \breve{s}=\frac{s}{r},\qquad\qquad\qquad
\breve{Y}=YKm/r.\label{rescale}
\end{aligned}
\end{equation}
We drop the $\breve{} \ $'s for convenience and study the
equivalent scaled version of model (\ref{ppreydelay}):
\begin{equation}
\left\{
\begin{aligned}
&\dot x(t)=x(t)(1-x(t))-y(t)x(t),\\
&\dot y(t)=-sy(t)+Ye^{-s\tau}y(t-\tau)x(t-\tau),\\\label{dimenless}
& (x(t),y(t))=(\phi(t),\psi(t))\in X, \, \mbox{for} \,\, t\in[-\tau,0],
\end{aligned}\right.
\end{equation}
where $X$ was defined in  (\ref{eq:X}).

First we address well-posedness of system (\ref{dimenless}).  For positive delay $\tau$, the existence and uniqueness of solutions of
system (\ref{dimenless})
was shown in Gourley and Kuang \cite{Kuang2004}. The following
proposition, proved in \ref{app:bounded}, indicates that for positive
delay the solutions remain
nonnegative and provides an  upper bound for each component.

\begin{proposition}\label{bounded}%lemma 2
 Consider model (\ref{dimenless}) with initial data in $X$.
\begin{enumerate}
\item The solutions  exist, are unique, and  remain
  nonnegative for all $t\geqslant 0$.
\item  $ \, \limsup_{t\rightarrow \infty}x(t)\leqslant 1 \, $ and
$ \, \limsup_{t\rightarrow \infty}y(t)\leqslant
\frac{1}{4s}Ye^{-s\tau}(s+1)^2$.
  \item Consider model (\ref{dimenless}) with initial data in $X^0$ where
  \begin{equation}\label{eq:X0}
    X^0=\{ (\phi(t),\psi(t)) \in X : \, \phi(0)>0 \,\,
    \& \,\, \exists \,  \theta \in[-\tau,0] \, \mbox{s.t.} \, \,
  \phi(\theta)\psi(\theta)>0.\}
\end{equation}
Then,  $x(t)>0$ for all $t>0$ and there exists $T\geqslant 0$ such that
$y(t)>0$ for all $t>T$.
\end{enumerate}
\end{proposition}

\section{Existence and stability of equilibria and uniform persistence} \label{analysis}

 Model
(\ref{dimenless}) can have up to three distinct
equilibria:
\begin{eqnarray}\label{notation1}
  E_0=(0,0), \  E_1=(1,0), \     E_+=(x_+(\tau),y_+(\tau))= \left(\frac{s}{Y}e^{s\tau},
1-\frac{s}{Y}e^{s\tau}\right).
\end{eqnarray}

The components of $E_+$ are nonnegative and $E_+$ is distinct from
$E_1$, if,
and only if, $0\leqslant \tau < \tau_c$, where
\begin{equation}\label{notation2}
\tau_{c}=\frac{1}{s}\ln\bigg( \frac{Y}{s} \bigg).
\end{equation}
Thus, when $\tau_c>0$, i.e., when $Y>s$, the components of $E_+$ are both positive, and   $E_+$ is referred to as the coexistence equilibrium.

When there is no delay, i.e. $\tau=0$ in (\ref{dimenless}), $E_0$ is
always a saddle attracting solutions with $x(0)=0$.
If $Y<s$, one of the components of $E_+$ is negative and so it is not
relevant, and $E_1$ is globally asymptotically stable with
respect to initial conditions satisfying $x(0)>0$ and $y(0)\geqslant 0$. When $Y=s$,  $E_1$ and $E_+$ coalesce  and are globally
attracting provided $x(0)>0$. If $Y>s$,  then $E_1$ is  a saddle
attracting solutions with $x(0)>0$ and $y(0)=0$  and
$E_+$ sits in int$\mathbb{R}^2_+$ and is global asymptotically stable
with respect to initial conditions in  int$\mathbb{R}^2_+$.

When $\tau>0$, to determine the local stability of each equilibrium
solution, we use the linearization technique for differential equations with
discrete delays (see Hale and Lunel  \cite{Hale1993}). After
linearizing (\ref{dimenless}) about any one of these equilibria,
 $(x^{\star},y^{\star})$, the characteristic equation,
 $P(\lambda)|_{(x^{\star},y^{\star})}=0$,
 is given by,
\begin{equation}
  (\lambda+s)(\lambda+y^{\star}-(1-2x^{\star}))
+Ye^{-(s+\lambda)\tau}x^{\star}(1-2x^{\star})-\lambda
Ye^{-(s+\lambda)\tau}x^{\star}=0 .
\label{chareqn}
\end{equation}

We summarize the results on local and global stability of
the equilibrium points and uniform persistence of the populations in the following theorem.  The proof can be found
in \ref{app:gasEi}.
\begin{theorem}%theorem 2
  \label{th:gasEi}
Consider (\ref{dimenless}).
\begin{enumerate}
  \item  Equilibrium $E_0$ is always unstable.
  \item  Equilibrium $E_1$ is
    \begin{enumerate}
      \item unstable if
$0\leqslant \tau <\tau_{c}$, and
\item globally asymptotically stable (with respect to $X^0$) if $\tau>
  \tau_{c}$, (i.e. if $\frac{se^{s\tau}}{Y}>1$).
\end{enumerate}
  \item Both components of $E_+$ are positive (i.e., $E_+$ exists), if,
    and only if,  $0\leqslant \tau < \tau_c$,  (i.e. $\frac{se^{s\tau}}{Y}<1$).
    \begin{enumerate}
  \item When $E_+$ exists and $\tau=0$, $E_+$ is globally asymptotically
    stable with respect to int$R^2_+$.
  \item When $E_+$ exists, and $\tau\geqslant 0$, model (\ref{dimenless}) is uniformly
    persistent with respect to initial data in $X^0$, i.e., there exists
    $\epsilon>0$ independent of $(\phi(t),\psi(t))\in X^0$ such that \\ $\liminf_{t\rightarrow\infty}
    x(t)>\epsilon$ and $\liminf_{t\rightarrow\infty}
    y(t)>\epsilon$.
\end{enumerate}
\end{enumerate}
\end{theorem}

%\medskip

Thus, for any fixed time delay $\tau$, if $\frac{s}{Ye^{-s\tau}}>1$,  only the prey
population survives and it converges to a steady state. On the other
hand, if the inequality is reversed,  for appropriate initial data both
the prey and the predator populations are uniformly persistent, i.e.
survive indefinitely.
However, we have not yet addressed what form the dynamics takes in the latter
case.

\subsection{Local stability  of $E_+$} \label{HB}

%%%%%%%%%%%%%%%%%%%%%%%%%%%%%%%%%%%%%%%%%%%%%%%%%%%%%%%%%%%
When $E_+$ exists, by Theorem~\ref{th:gasEi} both populations survive
indefinitely.  To address the possible forms the dynamics can take, we begin
by investigating the local stability of $E_+$ when it exists, i.e., when
$0\leqslant\tau<\tau_c$, and hence both components are positive.
Evaluating the characteristic equation (\ref{chareqn}) at $E_+$ gives
\begin{equation*}
P(\lambda)|_{E_+}=\lambda^2 +\lambda
s\left(1+\frac{e^{s\tau}}{Y}\right)+\frac{s^2}{Y}e^{s\tau}+
e^{-\lambda\tau}s\left(-\lambda +\left(1-\frac{2 s
e^{s\tau}}{Y}\right)\right)=0.
\end{equation*}
Therefore, $P(\lambda)|_{E_+}=0$ is of the form
\begin{equation}
P(\lambda)|_{E_+}=\lambda^2+p(\tau)\lambda+(q\lambda+c(\tau))e^{-\lambda\tau}+\alpha(\tau)=0,\label{CharE+}
\end{equation}
where
\begin{equation} \label{eqn:coeffs}
p(\tau)=s\left(1+\frac{e^{s\tau}}{Y}\right),\ \ q=-s,\ \ c(\tau)=
s\left(1-2\frac{s e^{s\tau}}{Y}\right),\ \ \mbox{and} \ \
\alpha(\tau)=\frac{s^2 e^{s\tau}}{Y},
\end{equation}
First assume that $\tau=0$. Then (\ref{CharE+}) reduces to
\begin{equation*}
\lambda^2 + ( p(0) + q )\lambda + (\alpha(0) +c(0) )=0.
\end{equation*}
Since $\alpha(0)+c(0)=s\left(1-\frac{s}{Y}\right)=s y_+(0)>0$
and $p(0)+q=\frac{s e^{s\tau}}{Y}>0$, by the Routh-Hurwitz
criterion  \cite{Gantmacher1959b},
 all roots of (\ref{CharE+}) have negative real part.  Therefore, $E_+$ is locally
asymptotically stable when $\tau=0$ and hence also for $\tau>0$ sufficiently
small.

We  consider the stability of $E_+$ as $\tau$ varies in the interval
$0<\tau<\tau_c$. Here,  $P(0)|_{E_+}=\alpha(\tau)+c(\tau)=s\
y_+(\tau)>0$ and so   $\lambda=0$ is not a root of (\ref{CharE+}).
Therefore, the only ways that $E_+$ can lose stability is:  (i)
  when one of the characteristic roots equals zero.  This only occurs
  when $\tau=\tau_c$.  This gives rise to a
  transcritical bifurcation
where $E_+$ coalesces with $E_1$   and then
disappears  as  $\tau$ increases
through $\tau_c$; (ii)
if characteristic roots bifurcate in from infinity; or (iii)  if a pair of
complex  roots with negative real parts and non-zero imaginary parts cross the imaginary axis as
$\tau$ increases from $0$, potentially resulting in  Hopf bifurcation.
 In \ref{app:imaginarycross} we prove that (ii) is impossible to obtain
 the following lemma.
%First we show that  roots bifurcating in from infinity  is impossible.

\begin{lemma}\label{lem:imaginarycross} %lemma 3
As $\tau$ increases from zero,  the number of
roots  of (\ref{CharE+}) with positive
real part can change only if a root appears on or crosses the
imaginary axis as $\tau$ varies.
\end{lemma}

%\medskip
In order to determine when Hopf bifurcations occur, we first
 determine for what values of $\tau$
pure imaginary roots of (\ref{CharE+})  exist so that (iii) can occur.
We will also be interested in secondary Hopf bifurcations.
%resulting in  more than one Hopf bifurcation as $\tau$ varies.

Suppose that $\lambda=i\omega$ $(\omega > 0)$ is a root of
$P(\lambda)|_{E_+}=0$, where $i=\sqrt{-1}$. Then
%\begin{equation}
$$
P(i\omega)|_{E_+}=-\omega^2+ i p(\tau) \omega +(i q \omega +
c(\tau)) e^{-i\tau\omega}+\alpha(\tau)=0. %\label{chpt1omegai}
$$
%\end{equation}
Using Euler's identity, $e^{i\theta} = \cos\theta + i\sin\theta$, and
%\begin{equation*}
%-\omega^2 +\alpha(\tau)+q
%\omega\sin(\tau\omega)+c(\tau)\cos(\tau\omega)+i\left(p(\tau)\omega+q
%\omega\cos(\tau\omega)-c(\tau)\sin(\tau\omega)\right)=0.
%\end{equation*}
equating the real and imaginary parts,   this is equivalent to
%\begin{subequations}\label{sincosleft}
%\left\{\begin{array}{ll}
\begin{eqnarray*}
%\begin{align}
       c(\tau)\cos(\tau\omega) + q \omega\sin(\tau\omega)=&\omega^2
       -\alpha(\tau), %\label{sincosleft:sin}
       \\
       c(\tau)\sin(\tau\omega) - q \omega\cos(\tau\omega)=&p(\tau)\omega.
        %\label{sincosleft:cos}
%\end{align}
\end{eqnarray*}
%     \end{array}\right.
%\end{subequations}
Solving for $\cos(\tau\omega)$ and $\sin(\tau\omega)$ gives
\begin{subequations}\label{chpt1sincos}
%\left\{
%\begin{array}{ll}
\begin{align}
&\sin(\tau\omega)=\frac{c(\tau)(p(\tau)\omega)+q\omega(\omega^2-\alpha(\tau))}{c(\tau)^2+q^2\omega^2}%:=h_1(\omega,\tau)
, \label{chpt1sincos:sin}\\
%&\\
&\cos(\tau\omega)=\frac{c(\tau)(\omega^2-\alpha(\tau))+q\omega(-p(\tau)\omega)}{c(\tau)^2+q^2\omega^2}.
\label{chpt1sincos:cos}%:=h_2(\omega,\tau)
%\end{array}\right.\label{chpt1sincos}
\end{align}
\end{subequations}
%Since $\sin(\tau\omega)=h_1(\omega,\tau)$,
%$\cos(\tau\omega)=h_2(\omega,\tau)$ and $\omega_+(\tau)$ solves the
%Therefore,
%\begin{equation}\label{h1h2square1}
%h^2_1(\omega,\tau)+h_2^2(\omega,\tau)=1,
%\end{equation}
%and so  $-1\leqslant h_i(\omega,\tau)\leqslant 1$ for $i=1,2$.
%Note that $p(\tau)=s+x_+(\tau),\
%q=-s,\ c(\tau)=s(1-2 x_+(\tau))$ and $\alpha(\tau)=s x_+(\tau)$.
%In the following theorem, we give an explicit form for $\omega$.
%Recalling that $\sin^2(\tau\omega)+\cos^2(\tau\omega)=1$,
Squaring
both sides of the equations in (\ref{chpt1sincos}), adding, and
rearranging gives
 \begin{equation}
\omega^4+(p(\tau)^2-q^2-2\alpha(\tau))\omega^2+\alpha(\tau)^2-c(\tau)^2=0.\label{omegapower4}
\end{equation}
Noting that \eqref{omegapower4} is a quadratic function of
$\omega^2$, we use  the quadratic formula to obtain
\begin{equation*}
	\omega_{\pm}^2(\tau)=\frac{1}{2}\bigg(q^2-p^2(\tau)+2\alpha(\tau) \pm
\sqrt{(q^2-p^2(\tau)+2\alpha(\tau))^2
-4\left(\alpha^2(\tau)-c^2(\tau)\right)}\bigg).\label{omega+-2}
\end{equation*}
Substituting using \eqref{eqn:coeffs}, it follows that 
%$u=\omega^2$ in \eqref{omegapower4} and solve for $u$ we obtain 
  \begin{equation}\label{chpt1omegapm_squared}
\omega_{\pm}^2(\tau)=\frac{1}{2}\bigg(-\bigg(\frac{se^{s\tau}}{Y}\bigg)^2
              \pm\sqrt{\bigg(\frac{se^{s\tau}}{Y}\bigg)^4+s^2\bigg(12\frac{s^2 e^{2s\tau}}{Y^2}
              -16\frac{se^{s\tau}}{Y}+4\bigg)}\bigg)  \ .	
\end{equation}
%From our analysis (please add where we did that), we know $u_+$ is
%positive and $u_-$ is negative. Therefore, $\omega_{+}(\tau)=\sqrt{u_+}$
%give us the only positive real root.  $\omega_{-}(\tau)=+\sqrt{u_}$ is
%a negative real root, and  $\pm\sqrt{u_{-}}$ are two
%complex roots, which is irrelevant to our Hopf bifurcation. }

In order to determine for what values of $\tau$ there are  positive real
roots of (\ref{omegapower4}), and hence candidates for  pure imaginary
roots, and possibly Hopf bifurcations, we define
\begin{equation}\label{tstardefined}
\tau^*=\frac{1}{s}\ln\left(\frac{Y}{3s}\right).
\end{equation}
We will prove that for $\tau \geqslant \tau^*$, there are no postive real
roots and  for $0 \leqslant \tau < \tau^*$ 
  \begin{equation}\label{chpt1omegapm}
\omega_{+}(\tau)=\sqrt{\frac{1}{2}\bigg(-\bigg(\frac{se^{s\tau}}{Y}\bigg)^2
              +\sqrt{\bigg(\frac{se^{s\tau}}{Y}\bigg)^4+s^2\bigg(12\frac{s^2 e^{2s\tau}}{Y^2}
              -16\frac{se^{s\tau}}{Y}+4\bigg)}\bigg)}
\end{equation}
is the only positive real root.

\begin{remark}
  Note that $\tau^*>0$, if, and only if, $\frac{s}{Y}<\frac{1}{3}$, and
  then $x_+(\tau^*)=\frac{s e^{s\tau^*}}{Y}=\frac{1}{3}$.
\end{remark}

In the following theorem, proved in \ref{app:omega+posicond}, we give
necessary conditions on $\tau$ for  Hopf bifurcations to occur.

\begin{theorem}\label{omega+posicond}
  Consider (\ref{dimenless}). Assume that both components of $E_+$ are
  positive.
  %, and $\frac{3}{Y}<\frac{1}{3}$.
  \begin{enumerate}
\item   $\omega_+(\tau^*)=0$.
If $\tau\geqslant\tau^*$, then (\ref{omegapower4}) has no positive real root.
Therefore, there can be no pure imaginary roots of (\ref{CharE+})  in this case. In
particular, if $\frac{s}{Y}\geqslant \frac{1}{3}$, then $\tau^*\leqslant 0$, and so
there can be no Hopf bifurcation of $E_+$ for any $\tau\geqslant 0$.
\item Assume that $\frac{s}{Y}<\frac{1}{3}$, and hence, $\tau^*>0$.
      If $\tau\in[0,\tau^*)$, then
$x_+(\tau)\in[\frac{s}{Y},\frac{1}{3})$, and  (\ref{omegapower4}) has
exactly one positive real
root, $\omega_+(\tau)$, given by (\ref{chpt1omegapm}). If
(\ref{CharE+}) has pure imaginary roots at $\tau$, and hence $\tau$ is a
candidate for Hopf bifurcation of $E_+$, then  $\tau\in(0,\tau^*)$ and $\omega_+(\tau)$ must
satisfy (\ref{chpt1omegapm}).
%\item $\omega_+(\tau)=0$ if, and only if, $\tau=\tau^*$.
  \end{enumerate}
%Moreover, $-1\leqslant h_i(\omega_+(\tau),\tau)\leqslant 1$ for $i=1,2$.
\end{theorem}

%\medskip

\begin{remark} \label{rem:candidates}
	\begin{enumerate}
		\item  Note  that even though $\omega_+(0)>0$ when
  $\frac{s}{Y}<\frac{1}{3}$, $\tau=0$ is not a candidate for a Hopf bifurcation,
  since  in this case all roots of (\ref{CharE+}) have negative real
  parts. Also, $\tau^*$ is not a candidate, since $\omega_+(\tau^*)=0$.
		\item   Haque \cite{Haque2012} finds an expression for
			$\omega_+(\tau^*)$ for a different scaling of
			the model.  However, he does not go on as we do
			in what follows, to determine when
			the equations given in  (\ref{chpt1sincos})
are both simultaneously  satisfied for 	$\omega_+(\tau^*)$, and to
			show that the Hopf bifurcations are nested and
			how the number of Hopf bifurcations increases as
			the death rate of the predator decreases.
	\end{enumerate}
\end{remark}

Substituting   the values of the coefficients given by
(\ref{eqn:coeffs}), in the right-hand side of  (\ref{chpt1sincos}), and recalling
that $x_+(\tau)=se^{s\tau}/Y$,  we define
%$h_1(\omega,\tau)$ and $h_2(\omega,\tau)$, we obtain
\begin{subequations} \label{h1_h2_defined}
  \begin{align}
h_1(\omega,\tau)&=\frac{\omega}{s}\left(\frac{s+x_+(\tau)-s
x_+(\tau)-2x_+^2(\tau)-\omega^2}{(1-2
x_+(\tau))^2+\omega^2}\right),\label{h1defined} \\
h_2(\omega,\tau)&=\frac{\omega^2(1+ s -x_+(\tau))-(1-2 x_+(\tau)) s
x_+(\tau)}{s\left((1-2x_+(\tau))^2+\omega^2\right)}.\label{h2defined}
\end{align}
\end{subequations}

Properties of the functions $h_1$ and $h_2$ are summarized
in \ref{app:h1_h2}.

\begin{remark}\label{rem:char_pure_imag}
By Theorem~\ref{omega+posicond} and Lemma~\ref{lem:h1_pos}, $\tau$ satisfies
(\ref{chpt1sincos}) for $\omega>0$, (and hence (\ref{CharE+}) has a pair of pure
imaginary eigenvalues) if and only if $\tau\in(0,\tau^*)$, and
\begin{subequations} \label{eqn:sin_h1_cos_h2}
%\left\{
%\begin{array}{ll}
   \begin{align}
&\sin(\tau\omega_+(\tau)) \, \, =h_1(\omega_+(\tau),\tau)
,  \label{eqn:sin_h1_cos_h2:sin}
\\
%&\\
&\cos(\tau\omega_+(\tau)) \, \,  =h_2(\omega_+(\tau),\tau).
\label{eqn:sin_h1_cos_h2:cos}
\end{align}
%\end{array}\right.
\end{subequations}
\end{remark}

 % Since  by part 1, %$\omega_+(\tau)$ satisfies (\ref{chpt1sincos}),
% $h_1^2(\omega_+(\tau),\tau)+h_2^2(\omega_+(\tau),\tau)=1$ and
%by part 2, $h_1(\omega_+(\tau)>0$ for $\tau\in[0,\tau^*)$, the result follows.
% \begin{comment}
% $h_2(\omega_+(\tau),\tau)\leqslant 1$, if, and only if,
%$$\omega^2_+(\tau)(1+ s -x_+(\tau))-(1-2 x_+(\tau)) s
%x_+(\tau)\leqslant s\left((1-2x_+(\tau))^2+\omega^2_+(\tau)\right) .$$
%After expansion and simplification, this is equivalent to
%$$\omega_+^2(\tau)\leqslant s(1-2 x_+(\tau)).$$
%Replacing $\omega_+^2(\tau)$ using its definition given in
%(\ref{chpt1omegapm}), after calculations similar to those in part 1.,
%this is equivalent to
%$$-s\leqslant(1-2x_+(\tau)),$$
%which clearly holds, and so $h_2(\omega_+(\tau),\tau)\leqslant
%1$ for all $\tau\in[0,\tau^*].$
%
%That $h_2(\omega_+(\tau),\tau)\geqslant
%-1$  for $\tau\in[0,\tau^*]$ is equivalent to
%$$\omega_+^2(\tau)(2s+1-x_+(\tau))+s(1-2x_+(\tau)(1-3x_+(\tau))\geqslant
%0.$$ Again, since this clearly holds, the result follows.
%\end{comment}

Define the function
\begin{equation}\label{eq:theta}
\theta:[0,\tau^*]\rightarrow [0,\pi]
\end{equation}
%\begin{equation}\label{notation3}
%\theta(\tau)=\arccos(h_2(\omega_+(\tau),\tau)),
%\end{equation}
 %\begin{equation}\label{eqn:def_theta}
$$\theta(\tau):=\arccos(h_2(\omega_+(\tau),\tau)).$$
%\end{equation}
 By part 3 of Lemma~\ref{lem:h1_pos}, stated and proved in \ref{app:h1_h2}, $\theta(\tau)$ is a well-defined, continuously
differentiable function.   Replacing  $\tau\omega_+(\tau)$ by
$\theta(\tau)+2 n \pi$ in the left hand side of
(\ref{eqn:sin_h1_cos_h2}),
we obtain
\begin{subequations} \label{eqn:sin_h1_cos_h2_theta}
    \begin{align}
&\sin(\theta(\tau)+2 n\pi) \, \, =h_1(\omega_+(\tau),\tau)
,  \label{eqn:sin_h1_cos_h2:sin_theta}
\\
&\cos(\theta(\tau)+2 n\pi) \, \,  =h_2(\omega_+(\tau),\tau).
\label{eqn:sin_h1_cos_h2:cos_theta}
\end{align}
\end{subequations}
% where $\theta(\tau)+2 n \pi$ is in the position of
% $\tau\omega_+(\tau)$ in  \eqref{eqn:sin_h1_cos_h2}.
%but represents a larger set of angels.
 Equation \eqref{eqn:sin_h1_cos_h2:cos_theta} is satisfied directly by the
 definition of $\theta(\tau)$.   Equation
 \eqref{eqn:sin_h1_cos_h2:sin_theta} is also satisfied,
 from  parts 1 and 2 of Lemma~\ref{lem:h1_pos}, since $0\leqslant
 \theta(\tau)\leqslant \pi .$
%and $0<\theta(\tau)\leqslant \pi$.
%In the following lemma,  we  define a continuous function $\theta(\tau)$ on
%$[0,\tau^*]$ that will  help determine the critical
%values of the delay that give a pair of pure imaginary eigenvalues of the
%characteristic equation.
%\begin{lemma}\label{omegataul} %theorem 4
%Assume $\tau\in[0,\tau^*]$.
%Define
%\begin{equation}\label{notation3}
%\theta(\tau)=\arccos(h_2(\omega_+(\tau),\tau)),
%\end{equation}
%(taking the branch where $\theta(\tau)\in[0,\pi]$).  Then, $\theta(\tau)$
%is continuous, and
%$0<\theta(\tau)\leqslant \pi$.

By comparing \eqref{eqn:sin_h1_cos_h2} and
\eqref{eqn:sin_h1_cos_h2_theta}, it follows that solutions of
\eqref{eqn:sin_h1_cos_h2} occur at precisely those points
where the curves $\tau\omega_+(\tau)$ and $\theta(\tau)+2n\pi$
intersect.
\begin{remark} \label{rem:characterization}
  By Remark~\ref{rem:char_pure_imag}, (\ref{CharE+}) has a pair of pure imaginary roots at precisely those points
where the curves  $\tau\omega_+(\tau)$ and $\theta(\tau)+2n\pi$
intersect  for $\tau\in(0,\tau^*)$,  where $n$ is a
nonnegative integer.
\end{remark}

%%%%%%%%%%%%%%%%%%%%%%%%%%%%%%
For each integer $n\geqslant 0$, denote the $j_n$ points of intersection of
the  curves $\tau\omega_+(\tau)$ and $\theta(\tau)+2n\pi$
 for $\tau\in(0,\tau^*)$, in increasing order,  by $\tau_n^j, \, \,
 j=1,2\dots,j_n$,
i.e., for each $n=0,1,2,\dots$,
$$\tau_n^j\omega_+(\tau_n^j)=\theta(\tau_n^j)+2n\pi, \quad
\ j=1,2,\dots,j_n.$$

\begin{theorem}\label{chpt3intersections}% theorem 3.8
Consider system (\ref{dimenless}).
%(where $\tau^*$ was defined by (\ref{tstardefined})).
The characteristic equation (\ref{CharE+}) has a  pair of pure
imaginary eigenvalues, if and only if, $\tau=\tau_n^j\in(0,\tau^*)$, a point of
intersection of the curves $\tau\omega_+(\tau)$ and
$\theta(\tau)+2n\pi$, for some integer $n\geqslant0.$
At  all  such intersections, the pair of pure imaginary eigenvalues is
 simple
and no other root of (\ref{CharE+})
is an integer multiple of $ i \omega_+(\tau_n^j).$
If in addition, $\frac{\mathrm{d}}{\mathrm{d}\tau}(\tau\omega_+(\tau))
\Big|_{\tau=\tau_n^j} \neq \frac{\mathrm{d}}{\mathrm{d}\tau}\theta(\tau)
\Big|_{\tau=\tau_n^j}$,
the transversality condition  for Hopf bifurcation,
$\frac{\mathrm{d}}{\mathrm{d}\tau}\mathrm{Re}(\lambda(\tau))
\Big|_{\tau=\tau_{n}^j}\neq0$,
holds.
\end{theorem}

\noindent
The proof is given in \ref{app:chpt3intersections}.

\begin{remark} \label{rem:test_transversality}
   If the slope of the curve  $\theta(\tau)+2n\pi$
  is less than the slope of   $\tau\omega_+(\tau)$ at an intersection
  point $\tau_n^j$, then   a pair of complex  roots of
  (\ref{CharE+}) crosses the imaginary axis from left to right as $\tau$
  increases through  $\tau_n^j$.
  On the other hand, if the slope of the curve  $\theta(\tau)+2n\pi$
  is greater than the slope of   $\tau\omega_+(\tau)$ at an intersection
  point $\tau_n^j$, then   a  pair of complex roots of
  (\ref{CharE+}) crosses the imaginary axis from right to left as $\tau$
  increases through  $\tau_n^j$.
\end{remark}

\begin{corollary}\label{cor:2k+1intersection}%theorem 3.9
Consider system (\ref{dimenless}). Assume that $\tau\in[0,\tau^*]$ and
that there exists $N\geqslant 0$ such that
$(2N+1)\pi\leqslant\max_{\tau\in[0,\tau^*]}\tau\omega_+(\tau)\leqslant
2(N+1)\pi$.
\begin{enumerate}
\item For $0\leqslant n\leqslant N$, $\theta(\tau)+2n\pi$ and
$\tau\omega_+(\tau)$ have at least two intersections in
$(0,\tau^*)$.
 \item For $n\geqslant N+1$,  $\theta(\tau)+2n\pi$ and
$\tau\omega_+(\tau)$  do not intersect in $(0,\tau^*)$.
\item
If $\theta(\tau)+2n\pi$ and
$\tau\omega_+(\tau)$ intersect for any $n\geqslant 0$, then $\tau_0^1$ is the
smallest  and $\tau_0^{j_0}$ the largest
value of $\tau$ for which
(\ref{CharE+}) has a pair of pure imaginary eigenvalues.
\item
  The
coexistence equilibrium $E_+$ is locally asymptotically stable
for $\tau\in[0,\tau_0^1)\cup(\tau_0^{j_0},\tau_c)$.
%If ???
%$E_+$ is unstable for $\tau\in(\tau_0^1,\tau_0^{j_0})$.
\end{enumerate}
\end{corollary}

\noindent
The proof is given in \ref{app:2k+1intersection}.

 \begin{figure}[bt!]
  \begin{center}
    \includegraphics[width=5.25cm]{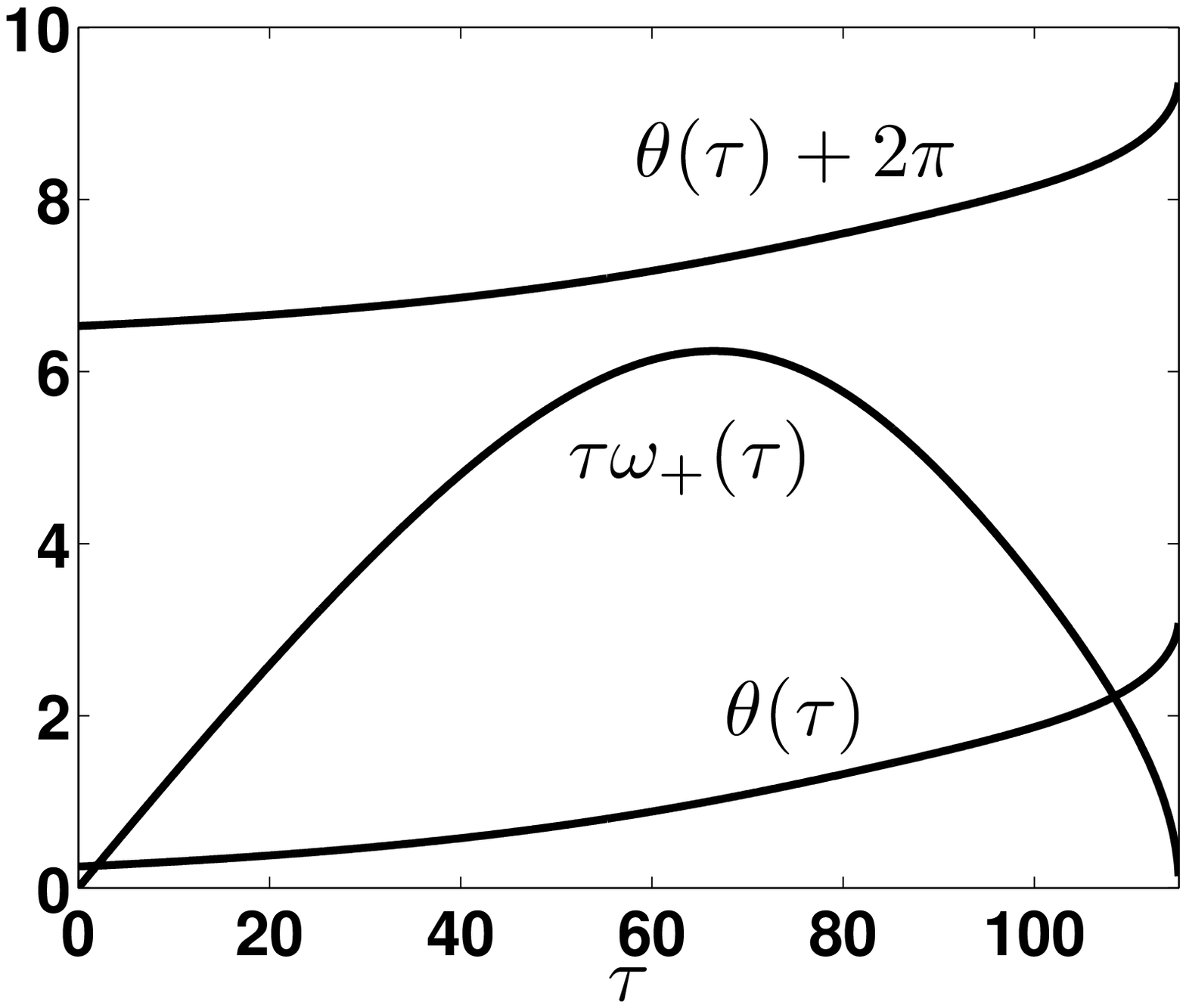}
 \includegraphics[width=5.25cm]{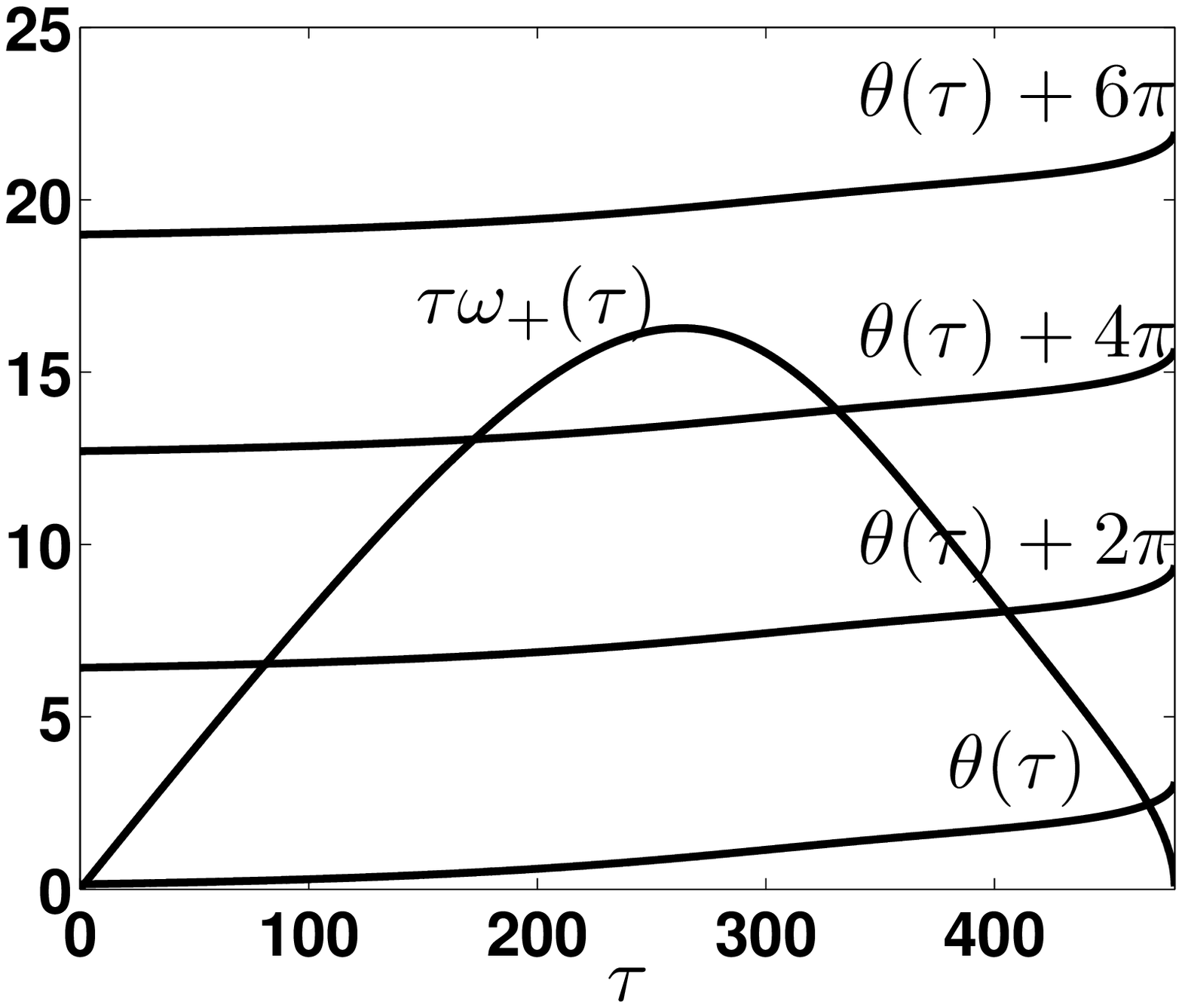}
\end{center}
\caption[  Intersections of $\theta(\tau)+2n\pi$ and $\tau
\omega(\tau)$] { Intersections of $\theta(\tau)+2n\pi$
and $\tau \omega_+(\tau), \ n=0,1,\dots, $.
Values of $\tau$ at which the characteristic equation has pure imaginary eigenvalues, and hence candidates for critical values
of $\tau$ at which there could be Hopf bifurcations.
In both graphs, at all such intersections, transversality holds, since
the slope of these curves
at these intersections are different.
Parameters: $ m=1,\ r=1,\ K=1, \ Y=0.6 . $
{\scriptsize\bf (LEFT)}  $s=0.02$.
 For $n=0$ there are two
intersections (i.e. $j_0=2$), at $\tau_0^1$ and $\tau_0^2$, but for $n=1$, and hence
$n\geqslant 1$, there are no
intersections. {\scriptsize\bf (RIGHT)}  $s=0.007$.
 There are two
intersections each (i.e. $j_n=2, \ n=0,1,2$), at $\tau_n^1$ and $\tau_n^2$, for $n=0,1$ and $2$,
but for $n=3$, and hence
$n\geqslant 3$, there are no
intersections. In both {\scriptsize\bf (LEFT)} and  {\scriptsize\bf (RIGHT)},  $E_+$ is
asymptotically stable for $\tau\in[0,\tau_0^1)\cup(\tau_0^2,\tau_c)$
and unstable for $\tau\in(\tau_0^1,\tau_0^2)$.}
\label{fig:thetatauomega2}
\end{figure}

 Our results  differ from those in Gourley and Kuang \cite{Kuang2004},
since  we give explicit formulas   for solutions of \eqref{omegapower4}
and define $\theta(\tau)$ explicitly. These explicit formulas play an
important role in analysis and make numerical simulations more
straightforward.  Although $\omega_{-}^2(\tau)$ in
equation
\eqref{chpt1omegapm_squared} is negative in
this model and hence its square root is not real, in other models
equation
\eqref{omegapower4} can have two positive real
solutions (see \cite{FanW}, where  both $\omega_+(\tau)$ and $\omega_-(\tau)$
are positive). In that case, double Hopf bifurcations are possible.

  In the next section we will demonstrate
  numerically that in the example shown in
  Figure~\ref{fig:thetatauomega2},
  $E_+$ first loses its stability through a supercritical Hopf
  bifurcation as $\tau$  increases through $\tau_0^1$ and then
  restabilizes  as a result of a second supercritical Hopf
  bifurcation as $\tau$   increases through $\tau_0^2$.  We will also
  show that between these two values of $\tau$ there is a
  sequence of bifurcations resulting in   interesting
  dynamics, including a strange attractor.

%%%%%%%%%new%%%%%%%%%%%%%%%%%%%%%%%%%%%%%%%%%%%%%%%%

\section{An example demonstrating complex dynamics}\label{sec:example}
In this section, unless specified otherwise,  we select the following
values for the parameters in model (\ref{ppreydelay}):
\begin{equation}
  m=1,\ r=1,\ K=1, Y=0.6,\ s=0.02, \label{example}
\end{equation}
  and consider $\tau$ as a bifurcation parameter. Since, with this
  selection of parameters, the model is already in the form of the scaled
 version of the model (\ref{dimenless}), it is not necessary to apply the scaling given by (\ref{rescale}).
We first use this example  to  illustrate the analytic
results given in section~\ref{analysis} where we provided necessary and sufficient
conditions for  a simple  pair of pure imaginary eigenvalues of the
characteristic equation  to occur as $\tau$ varies. We then provide
bifurcation diagrams with $\tau$ as the bifurcation parameter,
simulations  including time series and time delay embeddings
for various values of $\tau$, and a return map at a value of $\tau$
at which there is a chaotic attractor,
in order  to illustrate the wide variety of dynamics displayed by the
model,   even in the case when there are only two Hopf
bifurcations.
This includes two supercritical Hopf bifurcations, saddle-node  bifurcations of limit cycles and  sequences of period doublings that
appear to lead to chaotic dynamics with a strange attractor reminiscent
of the strange attractor for the well-known Mackey-Glass equation
\cite{Mackey:1977,Mackey:2009}
$$  \frac{dx}{dt}=\beta \frac{x(t-\tau)}{1+(x(t-\tau))^n}-\gamma x(t).$$

\subsection{Illustration of analytic results}

For the parameters given by (\ref{example}), the model has three
equilibria: $E_0=(0,0)$, $E_1=(1,0)$, that always exist, and
$E_+=(x_+(\tau),y_+(\tau))=( 0.0\dot{3} e^{0.02\tau},1-0.0\dot{3} e^{0.02\tau})$, given by
(\ref{notation1}).
The components of $E_+$ are both positive, if, and only if, $\tau\in[0,\tau_c)$,
where $\tau_c\approx 170$ (see (\ref{notation2})), and by part 3(b) of
Theorem~\ref{th:gasEi}  the model is uniformly persistent
for $\tau\in[0,170)$.

   $E_0$ is always a saddle, and hence unstable. $E_1$ is globally
asymptotically stable for $\tau>\tau_c$, and unstable for
$\tau\in[0,\tau_c)$.  For $\tau=0$, $E_+$ is asymptotically stable, and
by Lemma~\ref{lem:imaginarycross}, can only lose stability by means of a Hopf
bifurcation.  In order to determine th:e
critical values of $\tau$ at which
 there is a pair of  pure imaginary eigenvalues:
$\lambda(\tau)=\pm i \omega(\tau)$,  we consider the interval
$(0,115)$, since by Theorem~\ref{omega+posicond}, the interval on
which $\omega(\tau)$ is positive, is bounded above  by $\tau^*\approx
115$ (defined in (\ref{tstardefined})).

%1.234567901
With the parameters  given in (\ref{example}),  by \eqref{chpt1omegapm}  the  function
%\[\begin{array}{l}
	\[	\omega_+(\tau)= 
	 \sqrt{ -.\dot{5} \ 10^{-3}(e^{0.02\tau})^2 + 
	 \frac{
\sqrt{1.235\  10^{-6}(e^{0.02\tau})^4 +
	.5\dot{3} \ 10^{-5}e^{0.04\tau} - 0.5\dot{3}e^{0.02\tau} +
	4}}{2} } \]  
%\end{array}\]
 and recall,  by \eqref{eq:theta}
 $$\theta(\tau)=\arccos(\cos(\tau \omega_+(\tau))). $$
The intersections of the functions $\tau\omega_+(\tau)$  and
$\theta(\tau)+2n\pi,$  $n$ a nonnegative integer, in the interval
$(0,\tau^*)$,   give  the critical values of
$\tau$ for  which there is a pair of pure imaginary eigenvalues.
There are only two intersections as can be seen in
Figure~\ref{fig:thetatauomega2}~{\scriptsize\bf (LEFT)}.
%Since,  $0<\bar{\tau}=\frac{1}{s}\ln\left(\frac{Y}{6s}\right)
%<\tau^*$,
%$\max_{\tau\in[0,\tau^*]} \tau\omega(\tau)\geqslant\bar{\tau}(\omega_+(\bar{\tau}))
%\approx 5.73>\pi$,
Since,  $\pi<\max_{\tau\in[0,\tau^*]} \tau\omega(\tau)<2\pi$,
by Corollary~\ref{cor:2k+1intersection},
we are guaranteed at least
 two  values of $\tau>0$ at which the characteristic
 equation has a pair of pure imaginary roots.  In fact,
 (see Figure~\ref{fig:thetatauomega2} ({\scriptsize (LEFT})),  there are
 precisely two such values, $\tau_0^1$ and
 $\tau_0^2$.
% Here, $\omega_+(\tau)$ is defined by
%(\ref{chpt1omegapm}),  and $\theta(\tau)=\arccos(h_2(\omega_+(\tau),\tau)$
%where  $h_2(\omega_+(\tau),\tau)$ is given
% in (\ref{h2defined}) with  $x_+(\tau)=0.0\dot{3} e^{0.02\tau}.$
%The graphs of $\tau\omega_+(\tau), \ \theta(\tau)$, and
%$\theta(\tau)+2\pi$ are given in
%Figure~\ref{fig:thetatauomega2}~({\scriptsize TOP-LEFT}),  and show that
%there are exactly two intersections.
Using Maple \cite{maple},  we found that $\tau^1_0\approx 1.917$ and
$\tau^2_0 \approx 108.365$. By
Theorem~\ref{chpt3intersections},    the slopes of the curves
$\theta(\tau)$ and $\tau\omega_+(\tau)$ are different at these
intersection points, since these curves cross transversally (see
Figure~\ref{fig:thetatauomega2} {\scriptsize (LEFT)}), and hence the
transversality required for Hopf bifurcation holds at each
root.  By part 4 of Corollary~\ref{cor:2k+1intersection} and
Remark~\ref{rem:test_transversality}, the stability of $E_+$ changes
from asymptotically stable to unstable as $\tau$ increases through
$\tau_0^1$ and from unstable to asymptotically stable as $\tau$
increases through $\tau_0^2$, and is unstable for
$\tau\in(\tau_0^1,\tau_0^2)$.  But recall, even though $E_+$ is unstable
here, by  part 3(b) of Theorem~\ref{th:gasEi},  the model is still uniformly
persistent in this interval.
%since at these two values of $\tau$, all
%other roots of the characteristic equation have negative real parts.
%(see Remark~{rem: ??}.
% and hence there are two Hopf bifurcations (see \cite{Smith_delay},
%Theorem~6.1 page 90).
%  The graphs of $S_0(\tau)$ and $S_1(\tau)$ are
% given in  Figure~\ref{fig:thetatauomega2}~({\scriptsize TOP-RIGHT}), also
%demonstrating that $S^{\prime}_0(\tau_0^1)$
%and $S^{\prime}_0(\tau_0^2)$ are nonzero.

 Thus, we have shown that for the parameters chosen, there are exactly
  two candidates for  Hopf bifurcations (see  \cite{Smith_delay},
Chapter 6, Theorem~6.1, page 89-90).  That both Hopf bifurcations are
supercritical (involving first the birth, and then the disappearance of
orbitally asymptotically stable periodic solutions)
  will be demonstrated in the next section.  As $\tau$ increases through
  $\tau=1.917$, a family of orbitally asymptotically stable periodic
  orbits is born.  We will see that these periodic orbits  undergo additional bifurcations
 as $\tau$ increases, and that they disappear when $\tau$
  increases through the critical value $\tau=108.365$, at the second supercritical  Hopf bifurcation.

 By decreasing $s$,  the value of $n$ for which the curves
 %the value of
 % $N$ in Corollary~\ref{CoroIntersec} gets larger and larger
 %resulting
   $\theta(\tau)+2n\pi$ and $\tau\omega_+(\tau)$ intersect can increase.
   See
  Figure~\ref{fig:thetatauomega2} {\scriptsize (RIGHT)} for an example
  with $s=0.007$ for which the curves $\theta(\tau)+2n\pi$ and
  $\tau\omega_+(\tau)$ intersect when $n=0,1,2$. In fact, there are $6$
  points of intersection. We see by
  Remark~\ref{rem:test_transversality}, that for this example, $E_+$ is
  asymptotically stable until $\tau=\tau_0^1$, is
  unstable for $\tau\in(\tau_0^1,\tau_0^2)$, and finally becomes stable
  again for $\tau\in(\tau_0^2,\tau_c)$. Again, even though the model is
  unstable for $\tau\in (\tau_0^1,\tau_0^2)$, by part 3(b) of
  Theorem~\ref{th:gasEi}, it is uniformly persistent
  in this interval, since the model is uniformly persistent whenever
  both components of   $E_+$ are positive, i.e. for $\tau\in[0,\tau_c)$.

%%%%%%%%%%%%%%%%%%%%%%%%%%%%%%%%%%%%%%%%%%%%%%%%%%%%%%%
\subsection{Saddle-node of limit cycles, period doublings, and chaotic dynamics}
\label{numerics}

The computations and figures in this section
were done using Maple \cite{maple},
MATLAB  \cite{MATLAB:2018}, and
XPPAUT \cite{Ermentrout2002}

%We demonstrate the  occurrence of Hopf bifurcation and other possible bifurcations ????
%We use software MATLAB
%by calling XPPAUT in silence to plot the

\begin{figure}[tbhp!]
  \begin{center}
 \includegraphics[width=10cm]{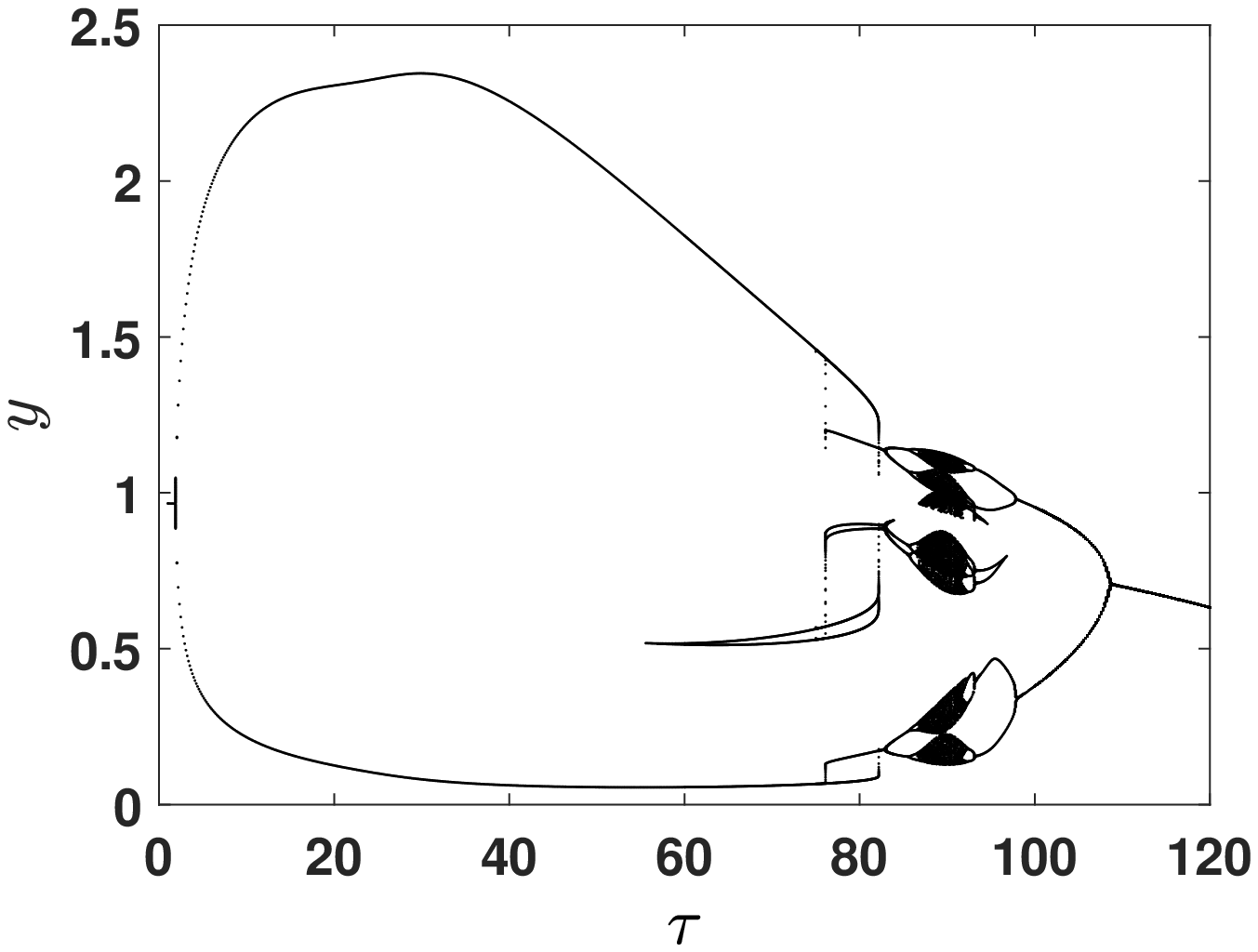}\\
 \includegraphics[width=10cm]{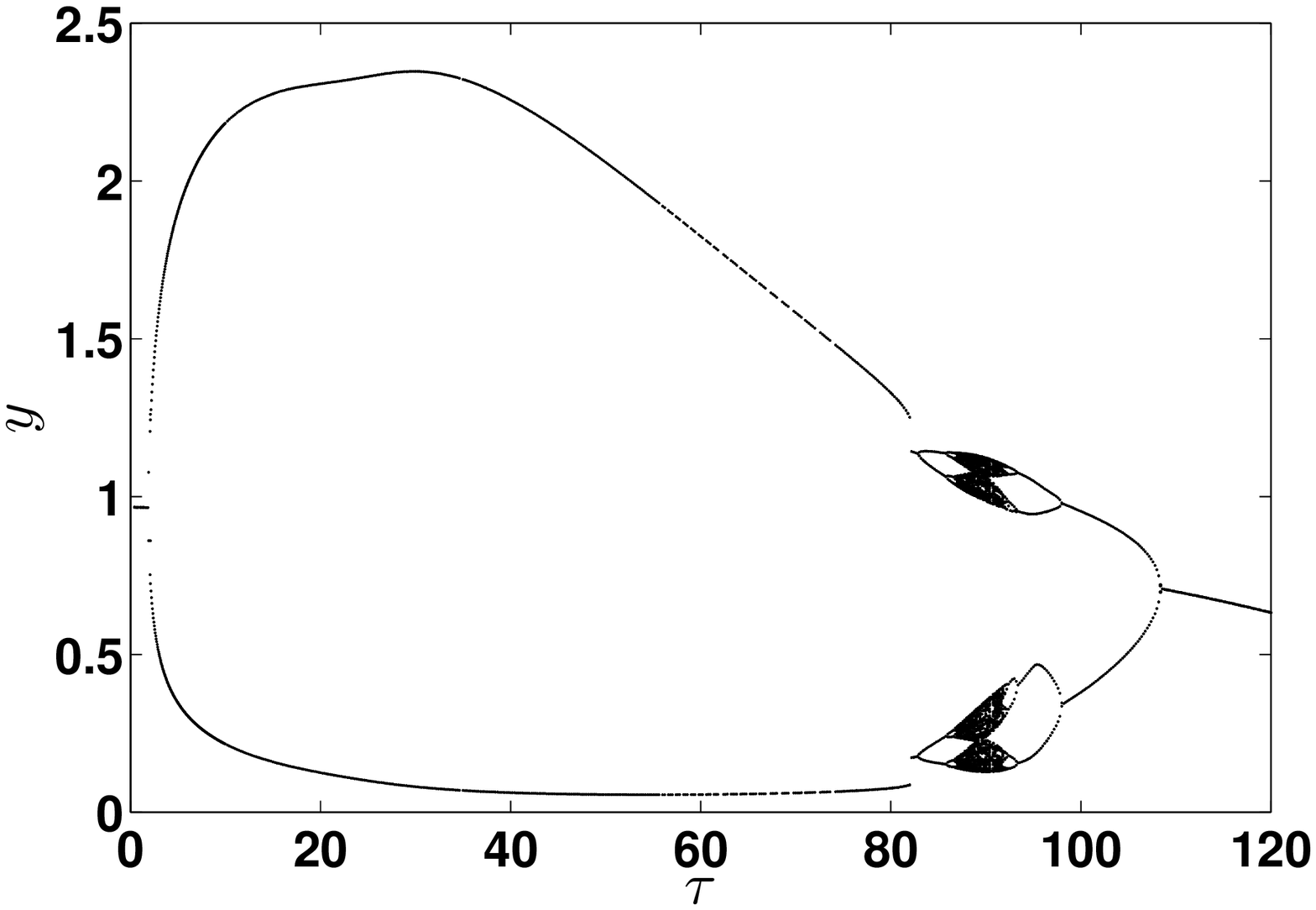}

  \end{center}
\caption[Parameter range for
	interesting dynamics]{Orbit diagrams.  Initial
	data was taken to be $x(t)=y(t)=0.1$ for
	$t\in[-\tau,0]$.  However, we found bistability in the portion
	of the diagram beween the vertical dots and varied the initial
	data as explained below.   
	%when increasing $\tau$ from the left    
	%and $x(t)=y(t)=0.1$ for
	%$t\in[-\tau,0)$ but $x(0) = 0.3, y(0) = 0.83.$ 
	Except for the
	portion between the vertical dots, the rest of the diagram was the same 
	for all of the initial conditions we tried (not shown).
{\scriptsize\bf (TOP)} All local maxima and minima
for the $y(t)$ coordinate of the attractor as $\tau$ varies,
including kinks.  {\scriptsize\bf (BOTTOM)}
Diagram including local maxes and mins
for the $y(t)$ coordinate  as $\tau$ varies, but with kinks eliminated.
 	There are two saddle-node of limit cycle
	bifurcations. They occur for
	$\tau$ approximately equal to  76 and 82, where  the curves  in the orbit diagrams stop abruptly and there
	appear to be vertical dots. For $\tau$ between these values, there
	is is bistability. Two orbitally asymptotically stable
	periodic orbits (with their maximum and minimum amplitudes
	shown) and an unstable periodic orbit with amplitudes between
	them (not shown).  The two stable periodic orbits were found
	by producing this part of the orbit diagram varying $\tau$
	forward and then varying  it backwards but startng at the last
	point of the attractor for the previous value of $\tau$.}
 \label{bif_diag_full}
 \end{figure}

\begin{figure}[tbhp!]
  \begin{center}
    \includegraphics[angle=270, width=7cm]{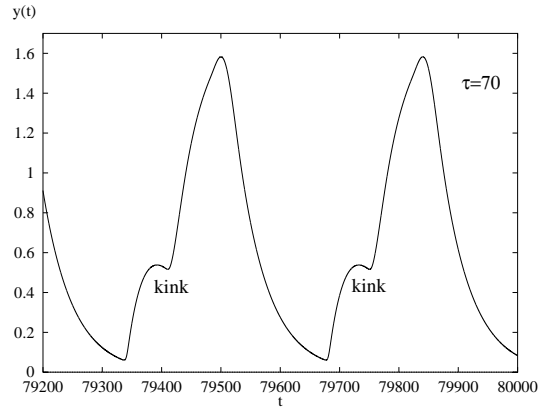}
  \end{center}
  \caption[kinks]{The time series for $y(t)$ when $\tau=70$, depicting
   kinks. There are two local maxima and two local minima over each period as shown  in  Figure~\ref{bif_diag_full} {\scriptsize (TOP)}, but only one local
  maxima and one local minima in Figure~\ref{bif_diag_full} {\scriptsize
  (BOTTOM)} in which kinks have been removed. }\label{kinks}
\end{figure}

In Figure~\ref{bif_diag_full}~{\scriptsize (TOP)}, for each value of
$\tau\in[0,120]$, starting with initial data $x(t)=y(t)=0.1$ for
$t\in[-\tau,0]$, we  integrate long enough for the solution to converge
to an attractor (e.g., equilibrium, periodic orbit, or strange
attractor), and then
 plot the local minima and maxima of the $y(t)$ coordinate
 on the attractor.
 Because we are interested in period doubling bifurcations, we then
 eliminate certain local maxima and minima that are due to kinks
 in the
 solutions  (see
 Figure~\ref{kinks}) rather than actual bifurcations, to obtain the graph in
 Figure~\ref{bif_diag_full} {\scriptsize (BOTTOM)}).
 Our solutions have
kinks for values of $\tau\in[55,98]$.  Kinks were also observed in
the Mackey-Glass equation  \cite{Mackey:2009}.

\begin{figure}[tbhp!]
  \begin{center}
	\includegraphics[width=10cm]{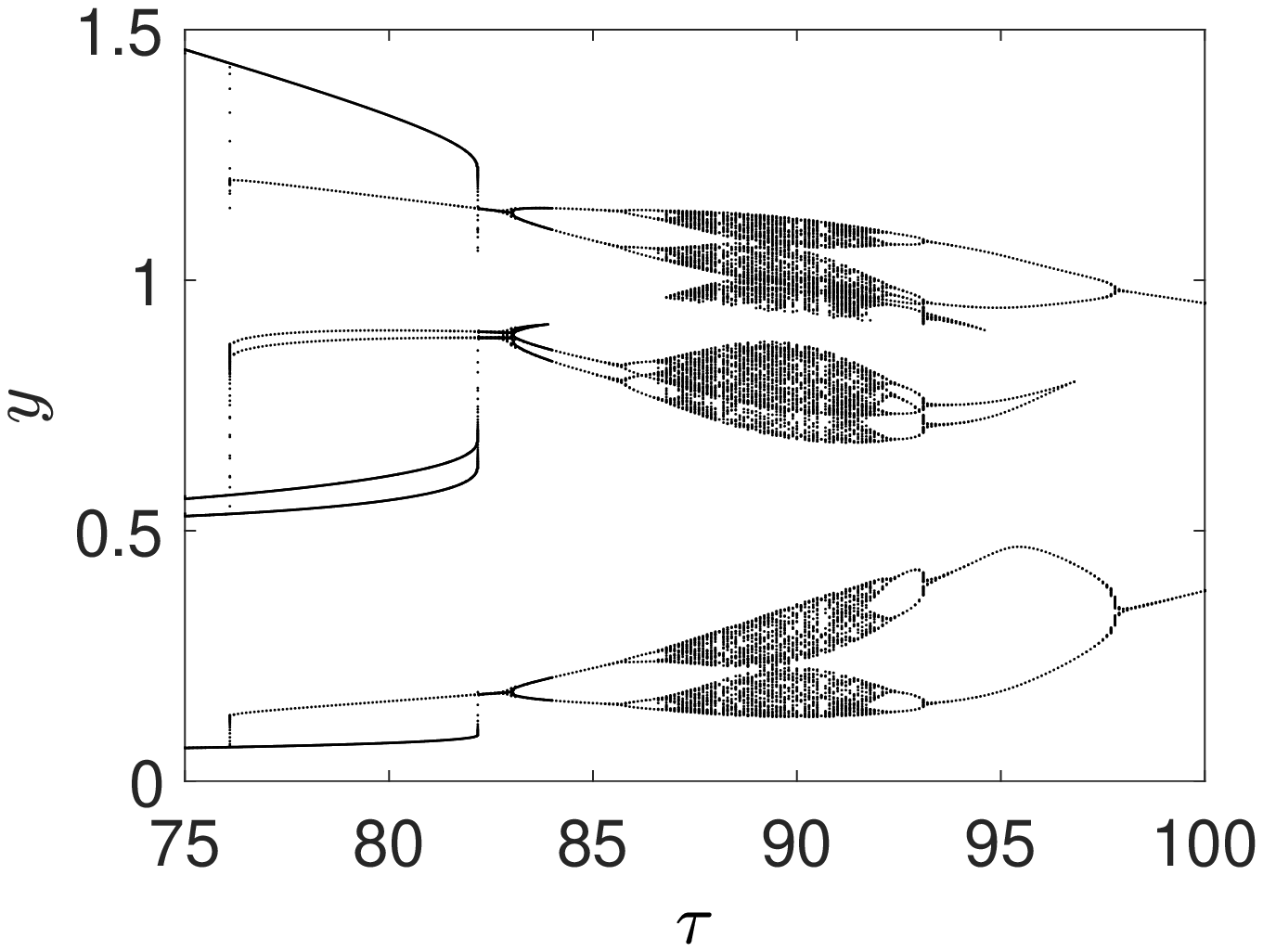}
  \end{center}
 \caption{ Zoom-in of orbit diagram shown in  Figure
\ref{bif_diag_full} for
	$\tau\in[75,100]$ % ({\scriptsize\bf TOP})
	including kinks.   The vertical dots indicate the boundary of the
	region of bistability, where the  two saddle-node of limit cycle
	bifurcations occur.
}
	%and ({\scriptsize\bf BOTTOM})
%without kinks.}
	\label{zoom-in_80_100}
\end{figure}

Figure \ref{bif_diag_full} confirms that there are two
Hopf bifurcations at $\tau\approx 1.917$ and $\tau\approx 108.365$, and
allows us to conclude that these Hopf bifurcations are both
supercritical, since they involve a family of orbitally asymptotically
stable periodic orbits.

\begin{figure}[htbp!]
	\begin{center}
	\includegraphics[angle=270, scale=.25]{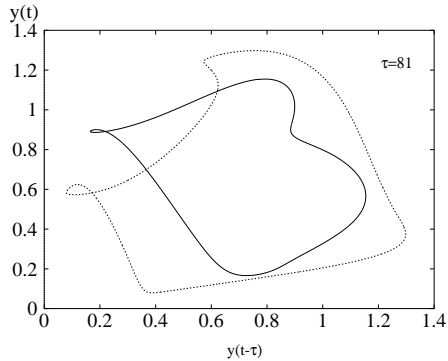}
\caption{  Time delay embedding  of two orbitally asymptotically stable periodic orbits
demonstrating bistability for $\tau=81$.  The one with larger amplitude (dashed) has initial data
$x(t)=y(t)=0.1$, for $t \in [-\tau,0]$, and period approximately
$345$.
%$353$.
The one with smaller amplitude (solid) has initial data
$x(t)=y(t)=0.1$, for $t \in [-\tau,0)$ and $x(0)=0.3, \ y(0)=0.83$ and has
period approximately
$273.7$.}\label{bistability}
\end{center}
\end{figure}

Next we focus on the  more interesting
dynamics observed for $\tau\in[80,100]$ (see  Figure~\ref{zoom-in_80_100}).
There appears to be a discontinuity in the
bifurcation diagram for $\tau\approx 82.225$.  Upon further
investigation we have determined that there is a saddle-node bifurcation
of limit cycles at this value of $\tau$ and another saddle-node
bifurcation of limit cycles for a value of $\tau$ smaller than
$\tau=81$.  For values of $\tau$ between these two saddle-node bifurcations,
there is bistability.  There are two orbitally asymptotically stable
period orbits.  An example of two such orbits is given in
Figure~\ref{bistability}, where $\tau=81$.

\begin{figure}[tbhp!]
   \begin{center}
\includegraphics[angle=270, width=6cm]{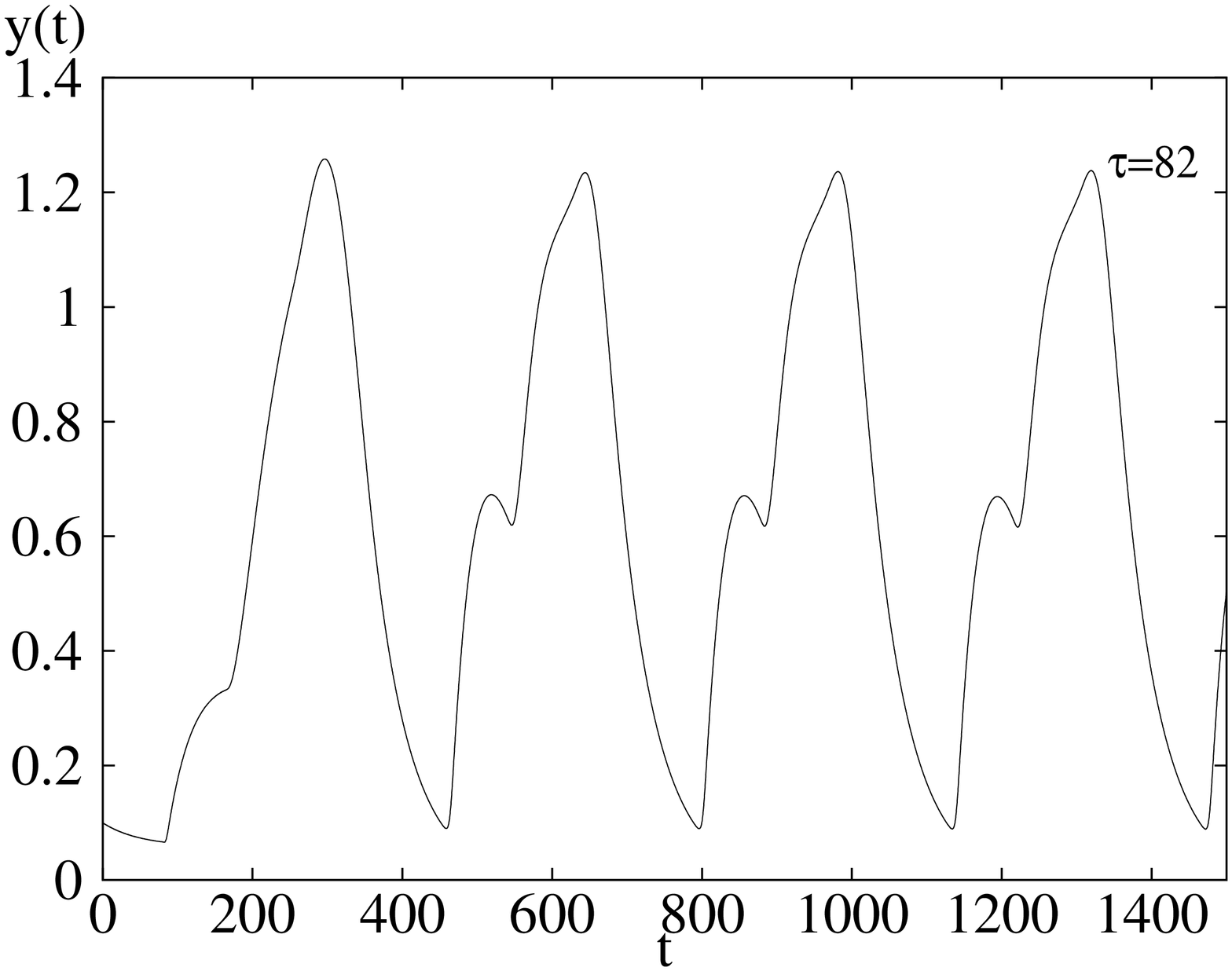}
\includegraphics[angle=270, width=6cm]{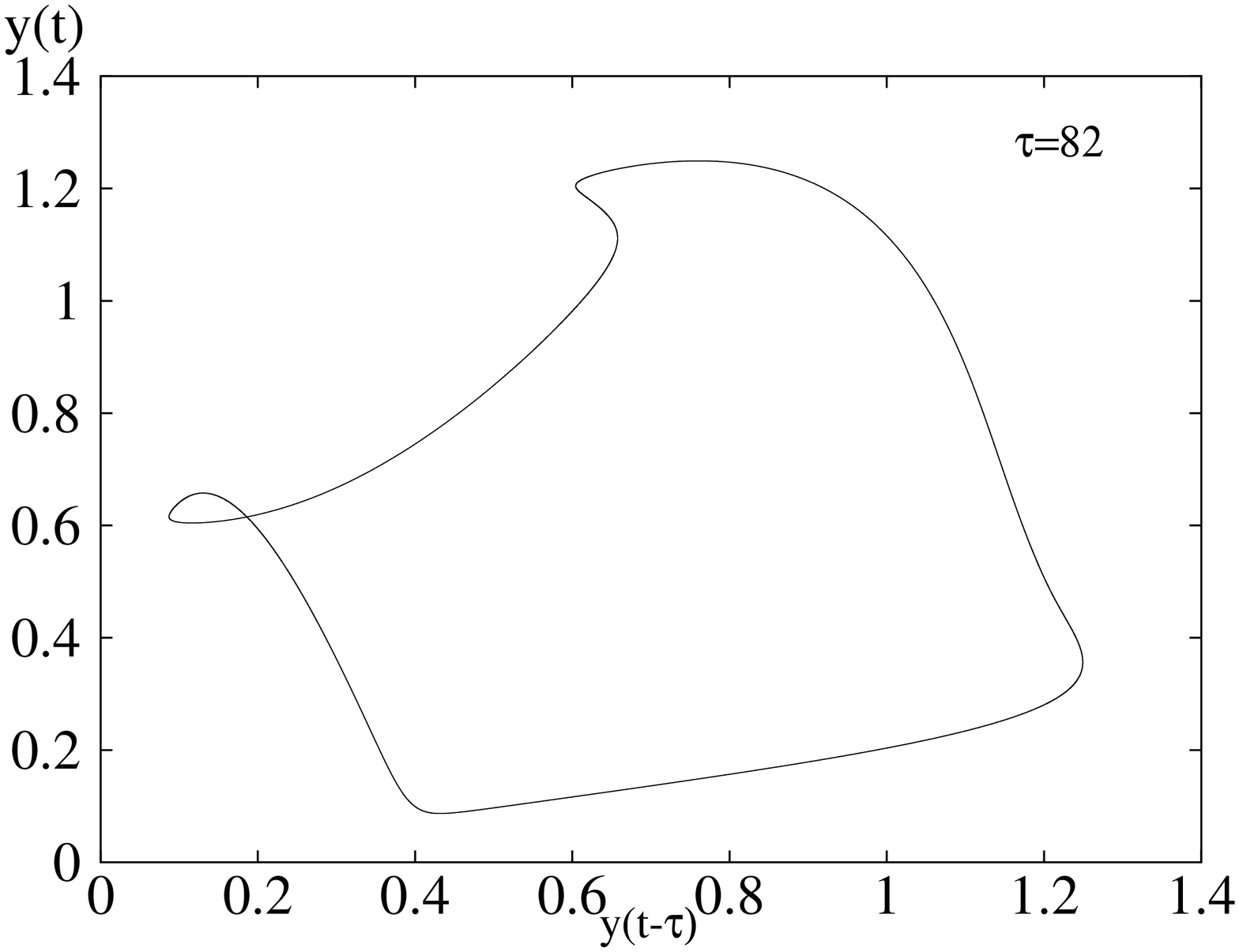}
\includegraphics[angle=270, width=6cm]{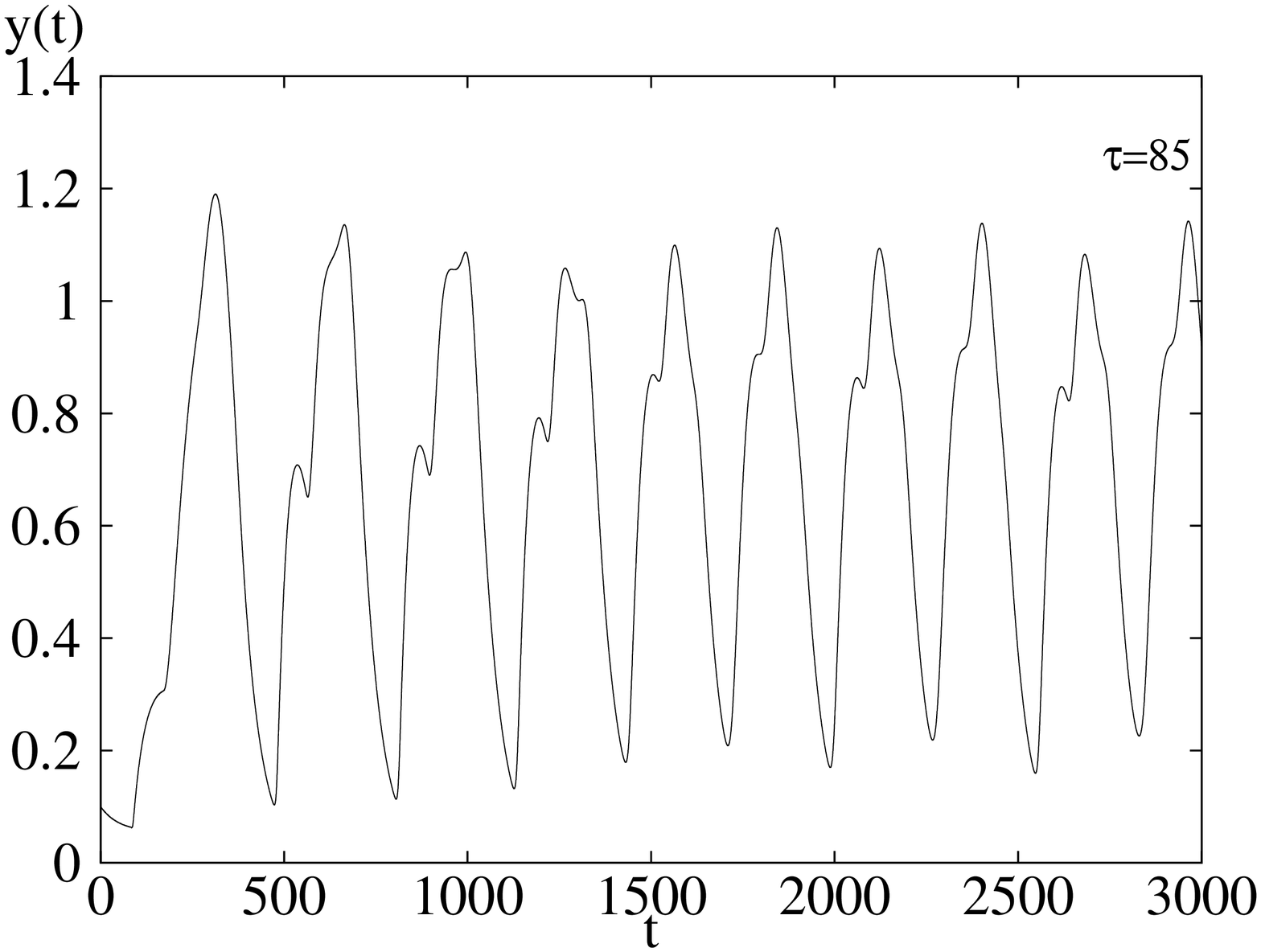}
\includegraphics[angle=270, width=6cm]{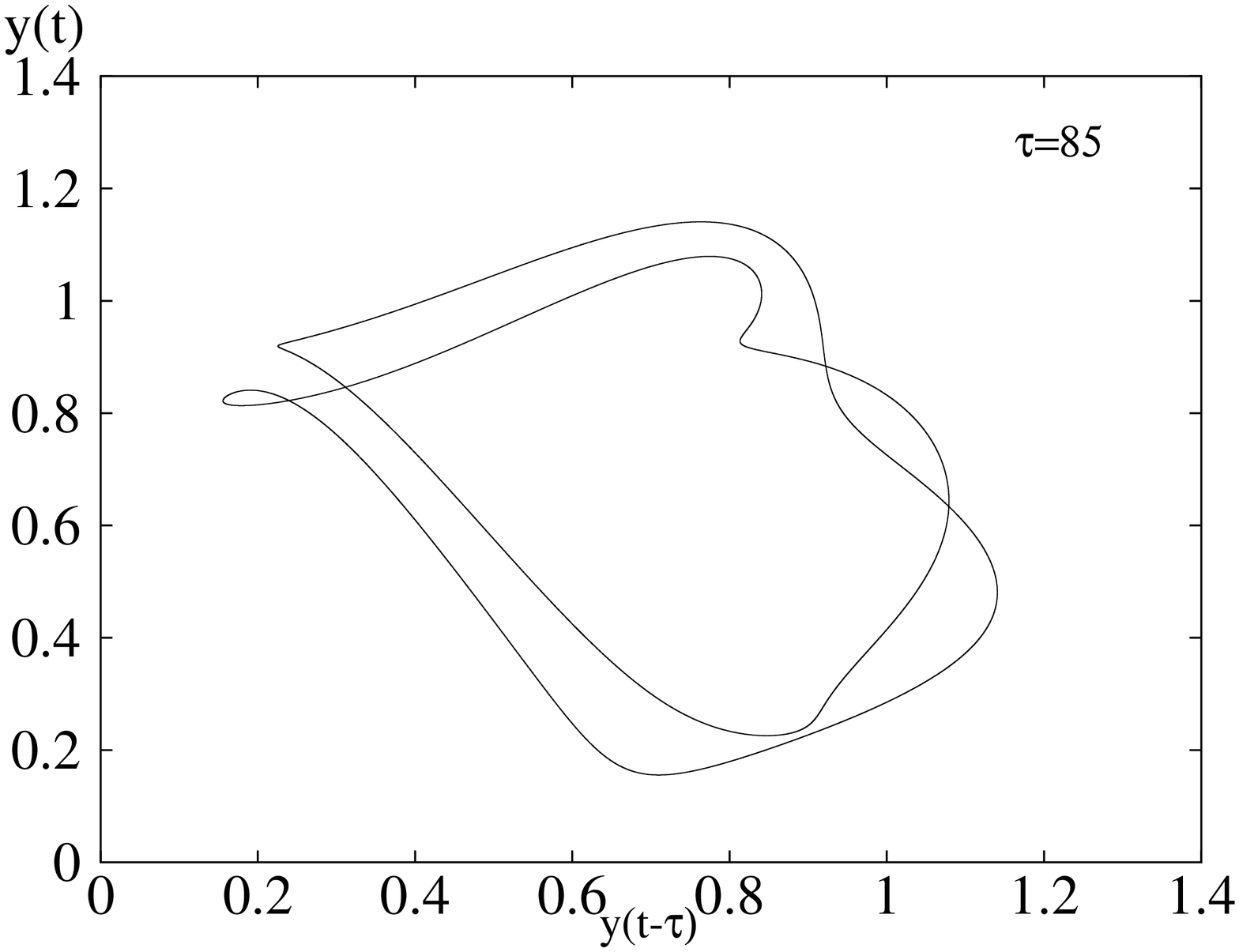}
\includegraphics[angle=270, width=6cm]{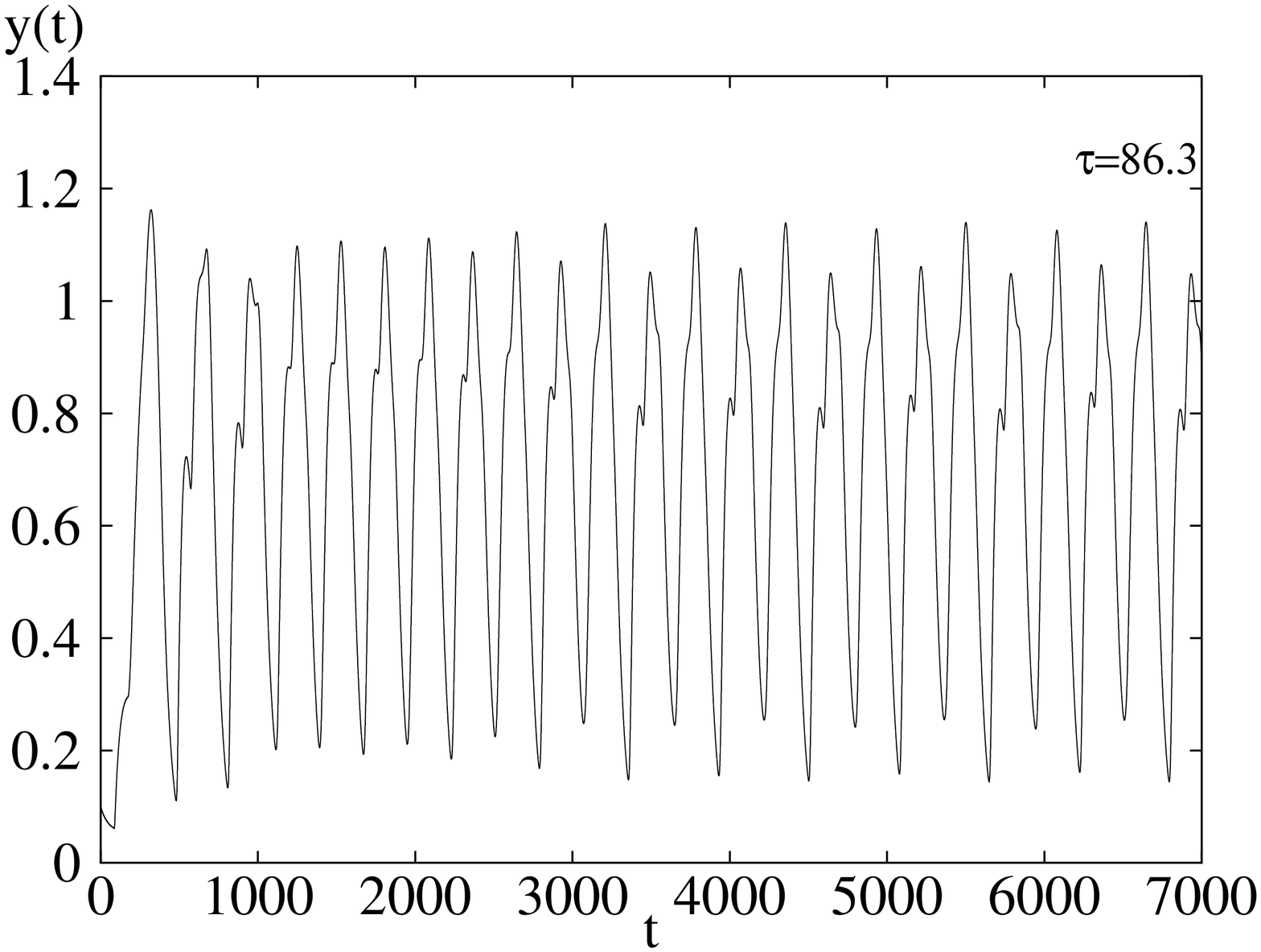}
\includegraphics[angle=270, width=6cm]{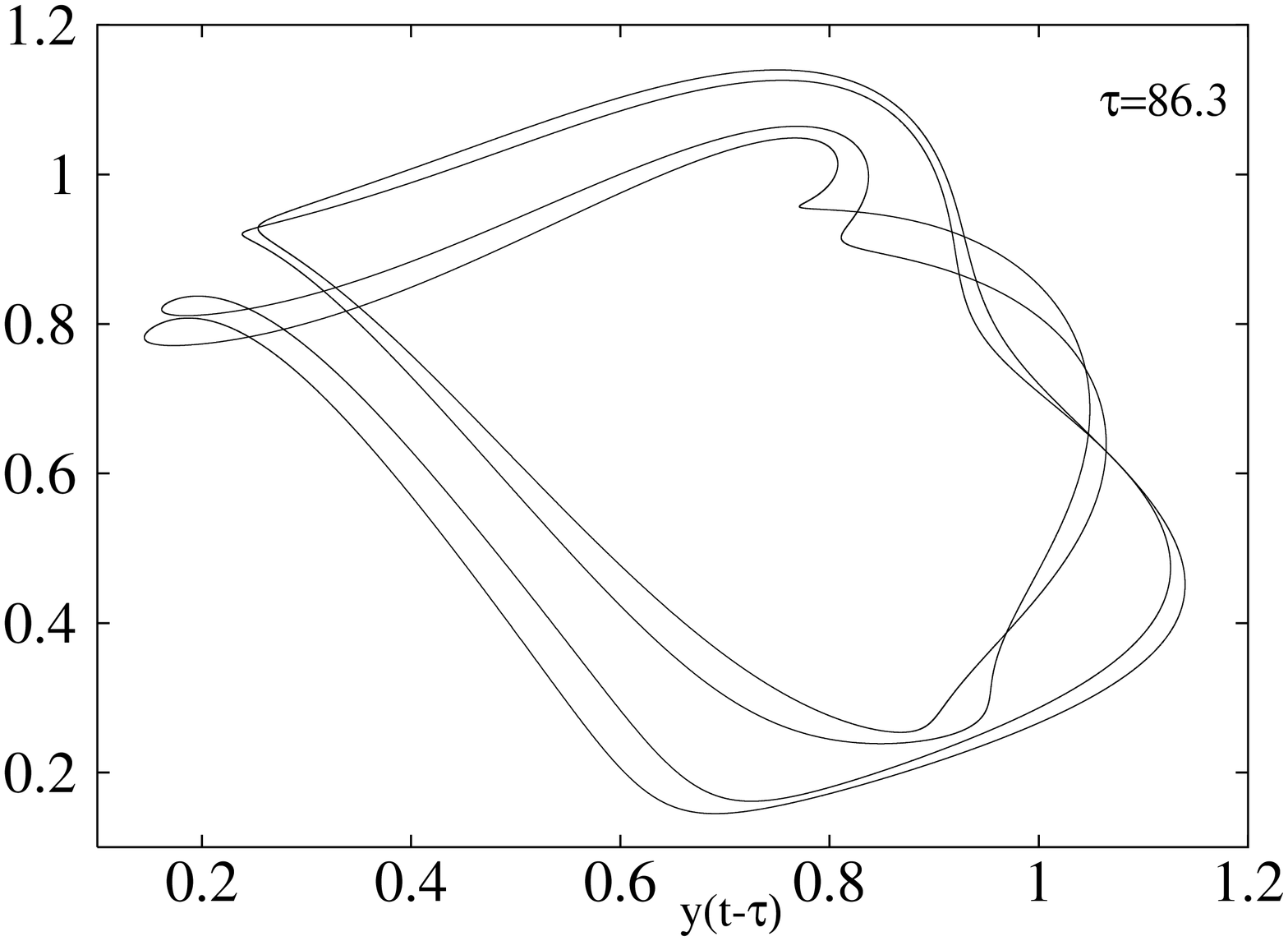}
\includegraphics[angle=270, width=6cm]{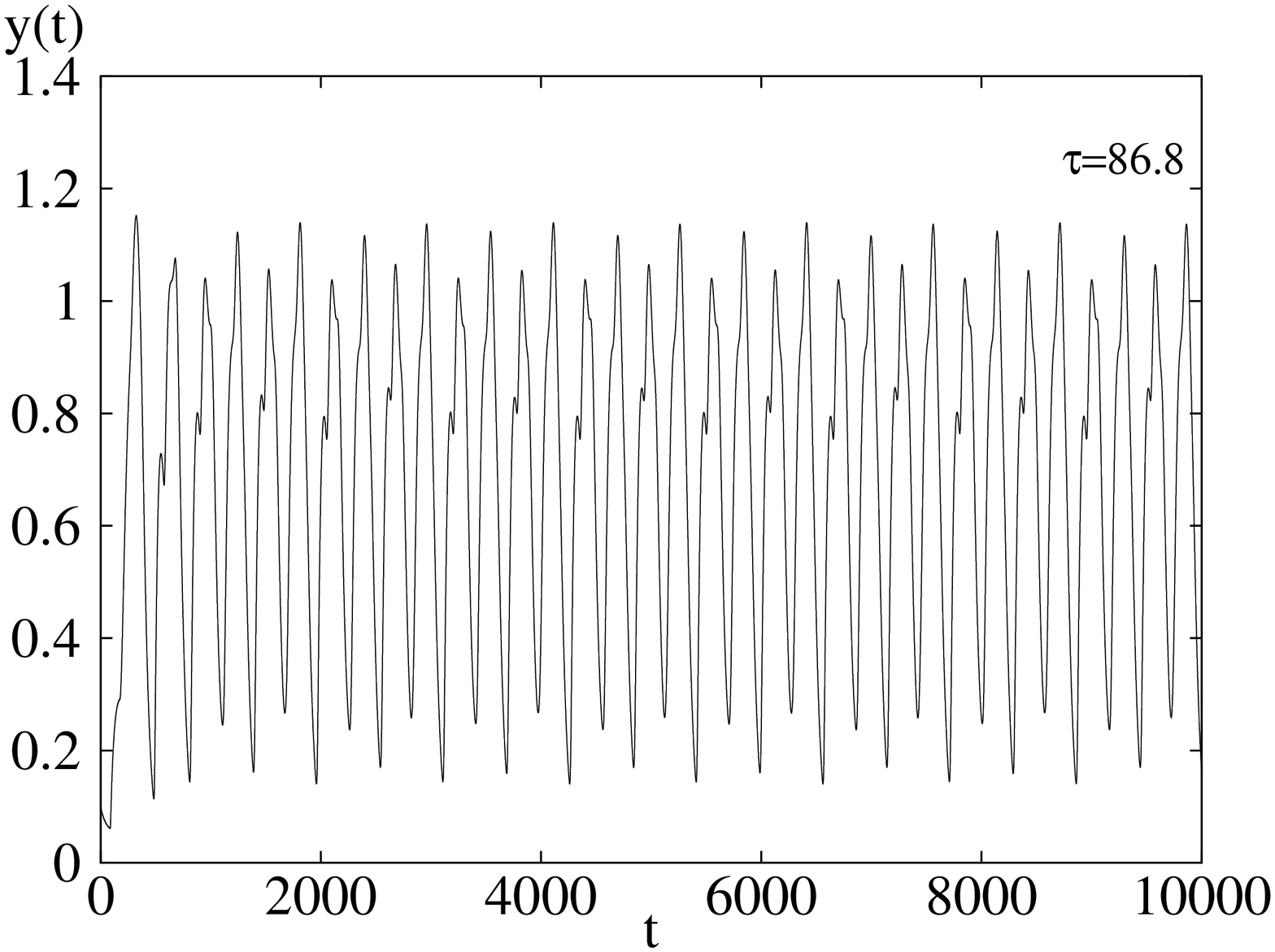}
\includegraphics[angle=270, width=6cm]{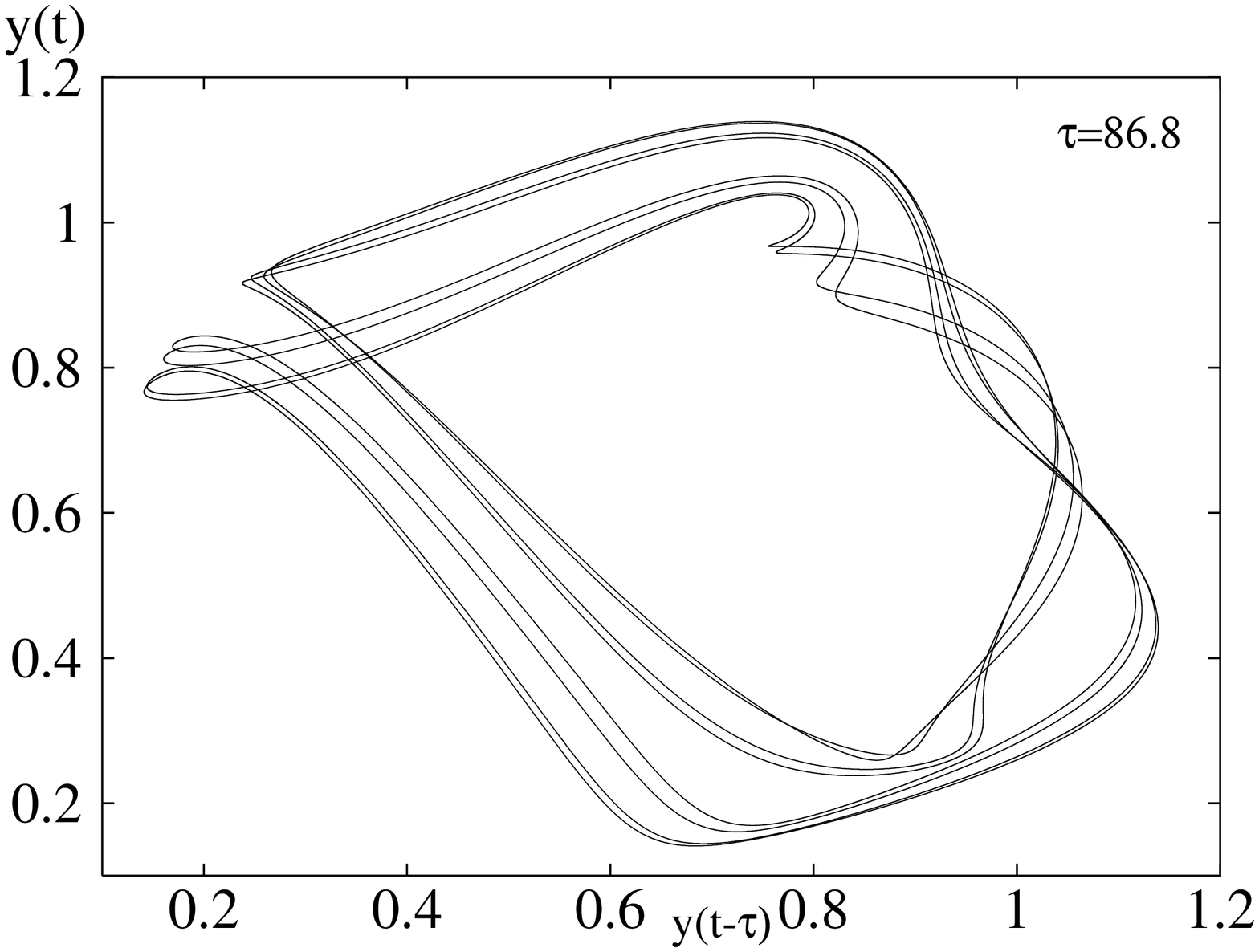}
   \end{center}
   \caption[ ]{ {\scriptsize\bf (LEFT)} Time series starting from the initial data $x(t)=y(t)=0.1, t\in[-\tau,0]$
   indicating how quickly the orbit gets close to the periodic attractor and
   {\scriptsize\bf (RIGHT)} time delay embeddings of
   the periodic attractors, demonstrating the sequence of  period doubling bifurcations
   initiating from the left at
   $\tau\approx 83,\ 86$, and $86.6$. Values of $\tau$ selected between these bifurcations:
    $\tau=82,\ 85,\
   86.3$, and $86.8$, with periods of the periodic attractor approximately equal to:
   $340.2,\ 564.6,\ 1144.5$, and $2298.3$, respectively, are
   shown. } \label{PD_left}
\end{figure}

\begin{figure}[tbhp!]
   \begin{center}
\includegraphics[angle=270, width=5.cm]{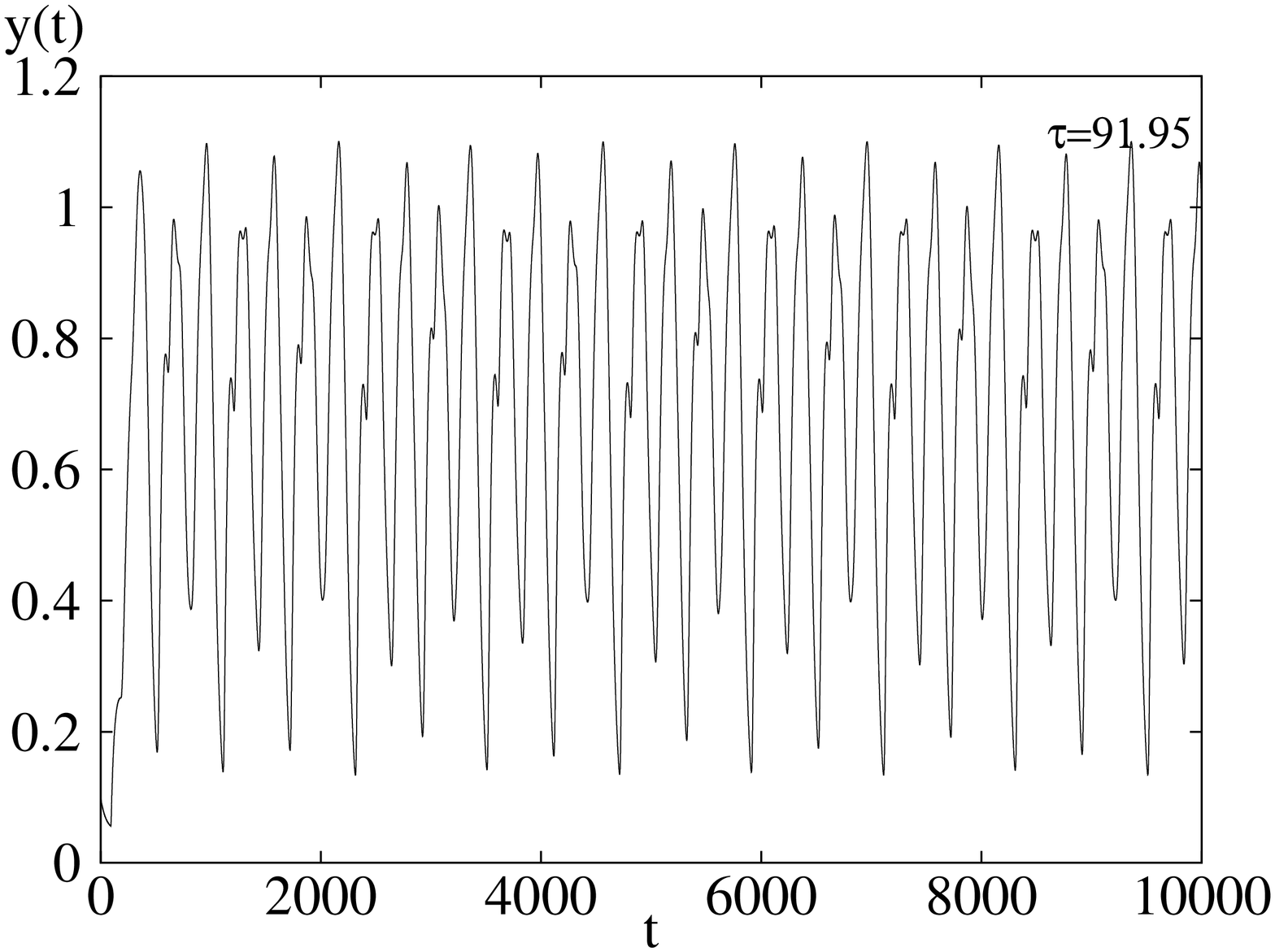}
\includegraphics[angle=270, width=5.25cm]{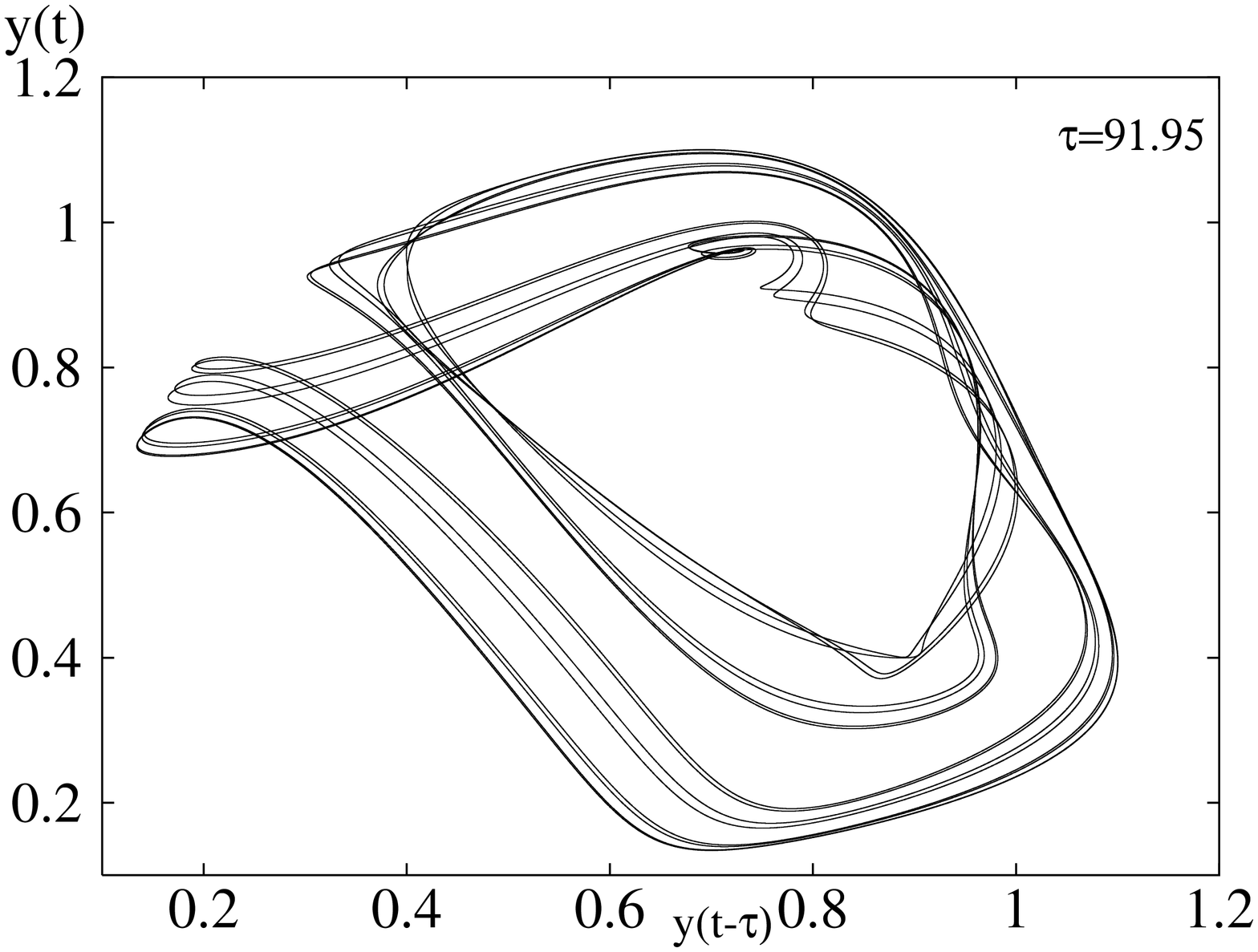}
\includegraphics[angle=270, width=5.cm]{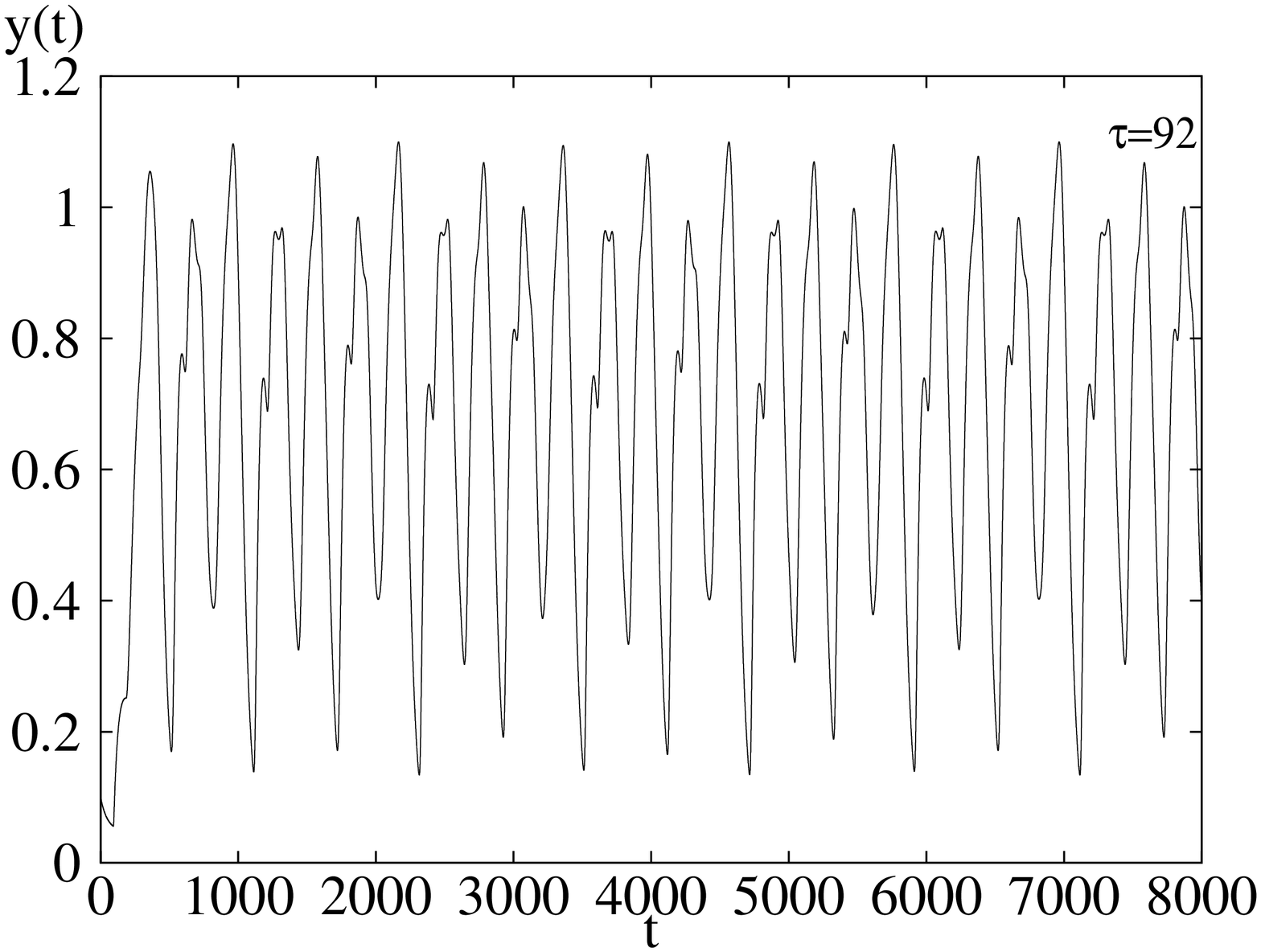}
\includegraphics[angle=270, width=5.cm]{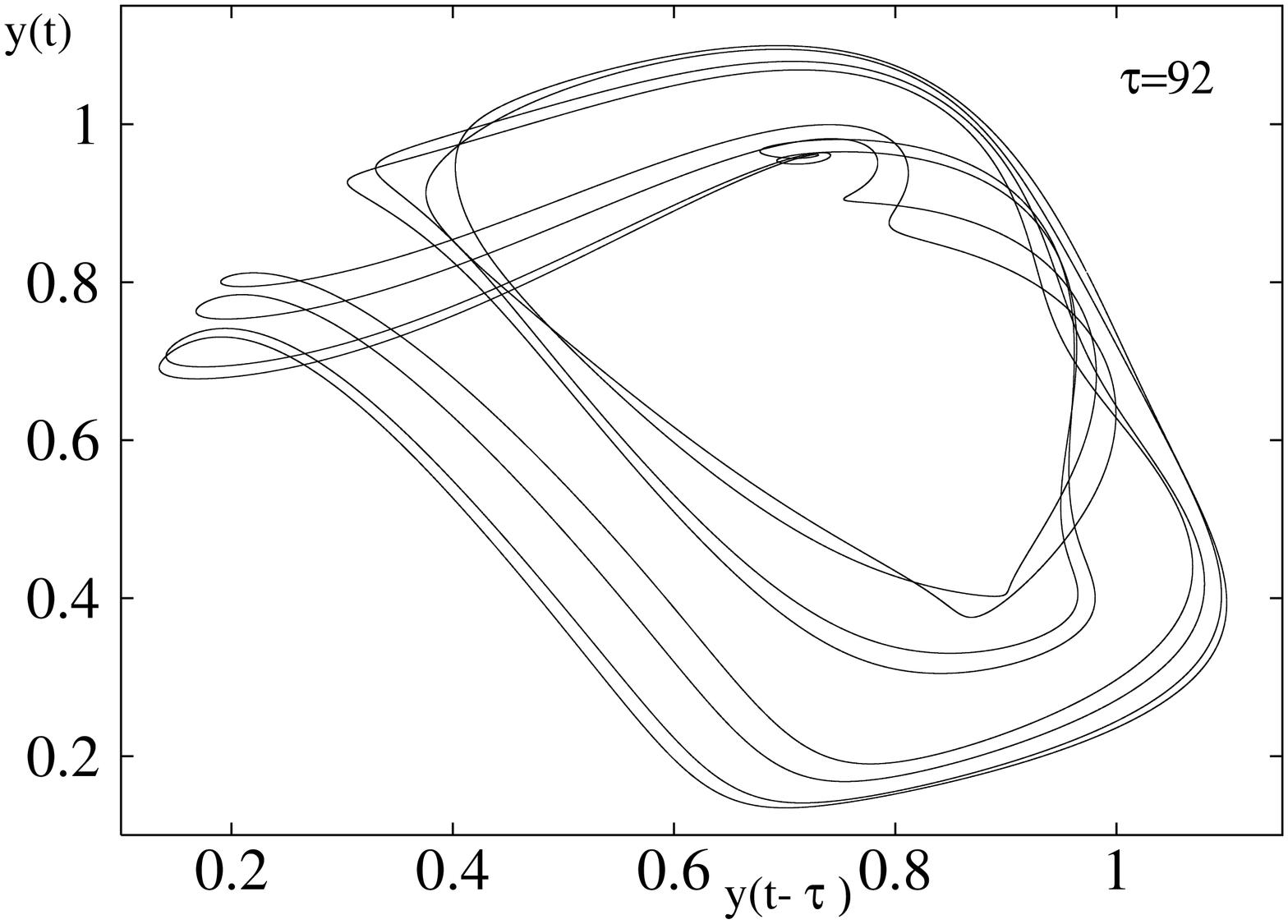}
\includegraphics[angle=270, width=5.cm]{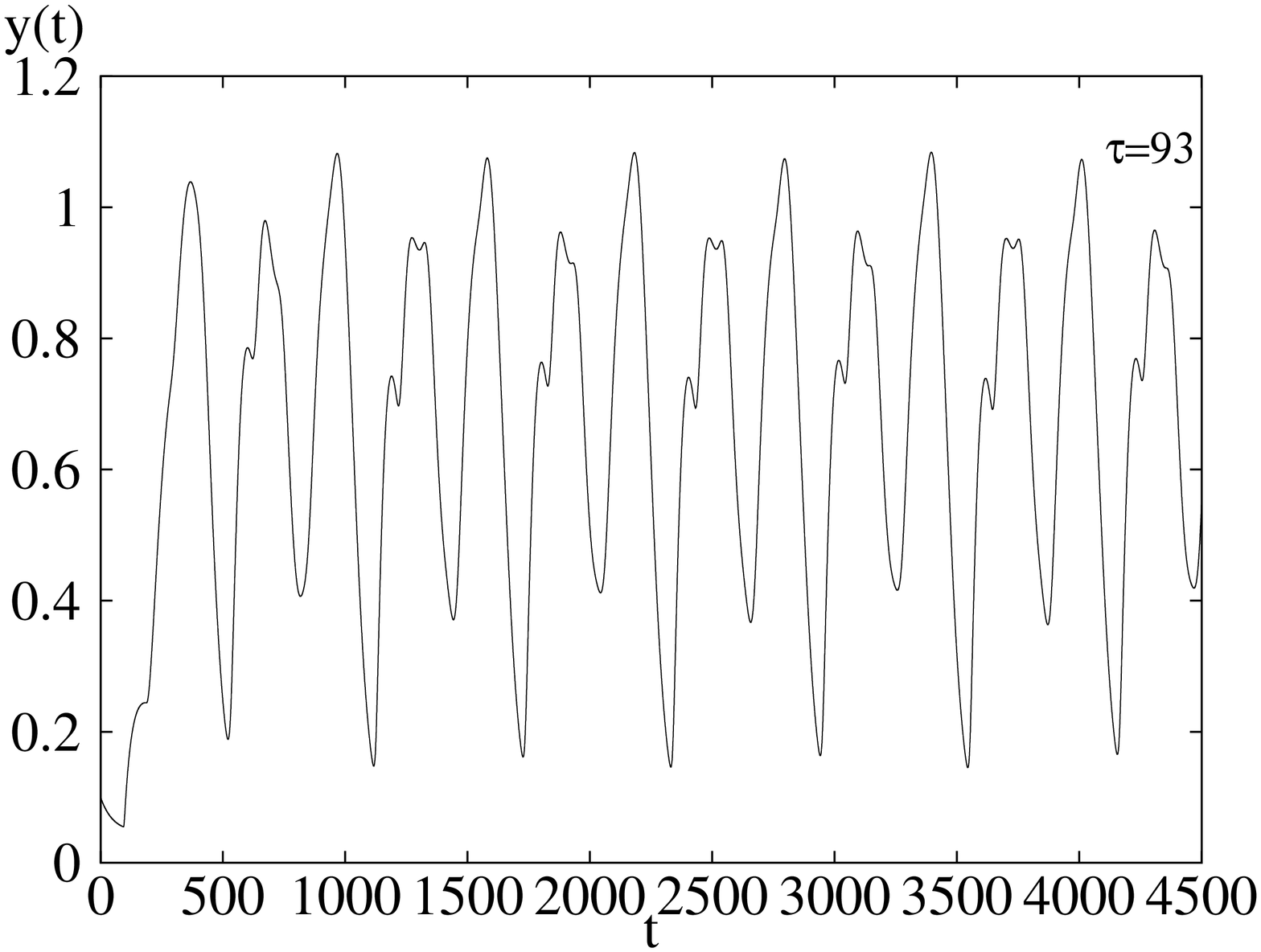}
\includegraphics[angle=270, width=5.cm]{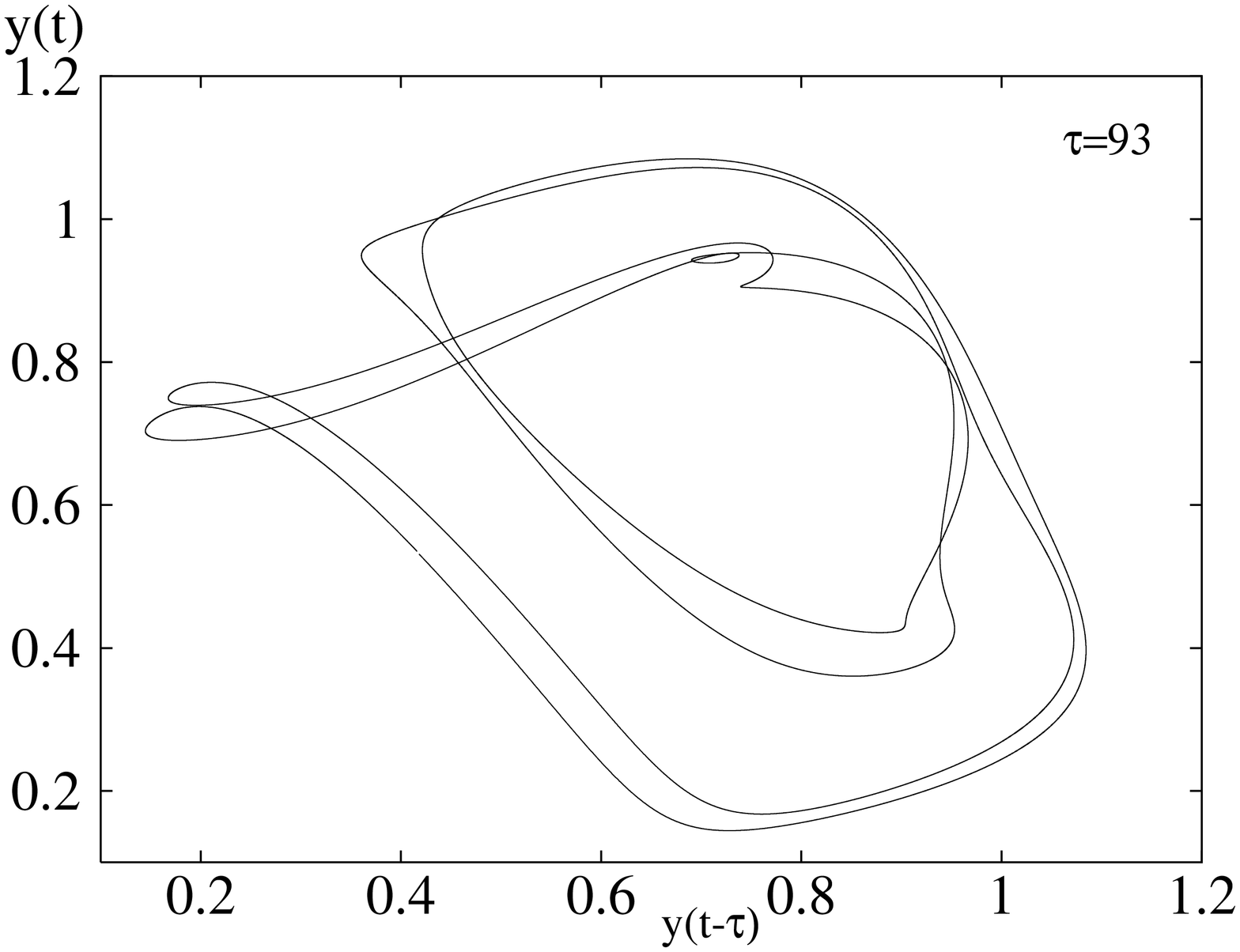}
\includegraphics[angle=270, width=5.cm]{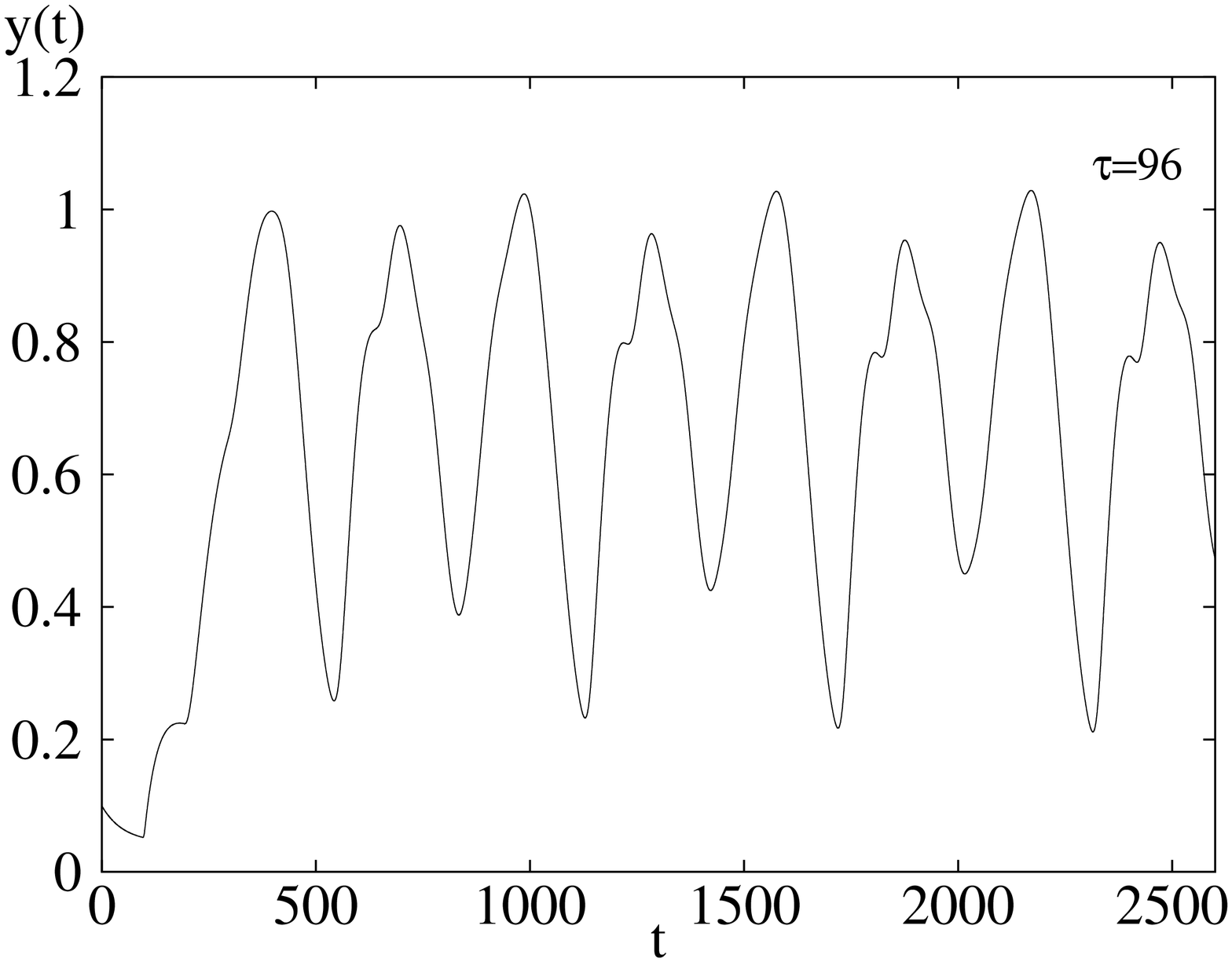}
\includegraphics[angle=270, width=5.cm]{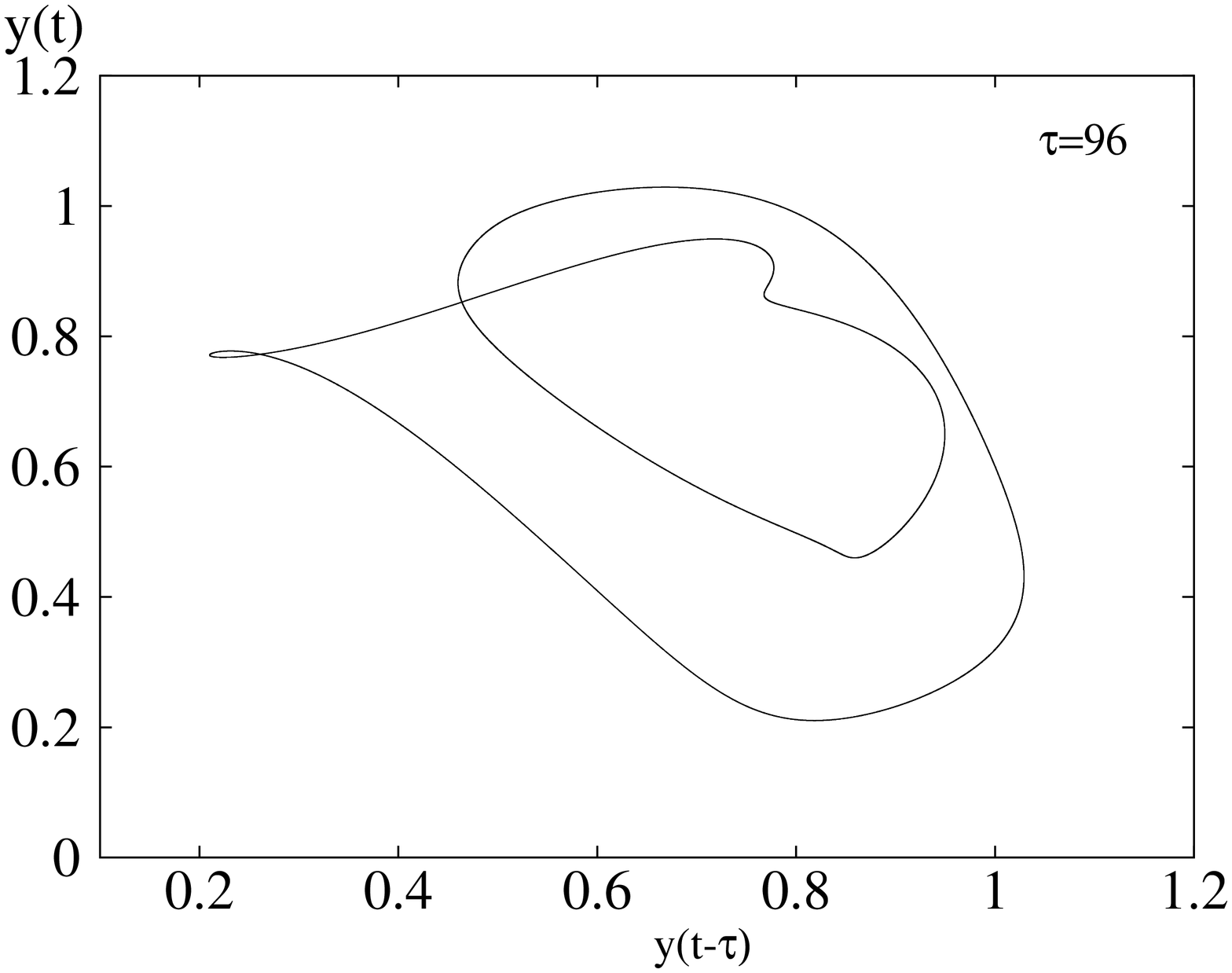}
\includegraphics[angle=270, width=5.cm]{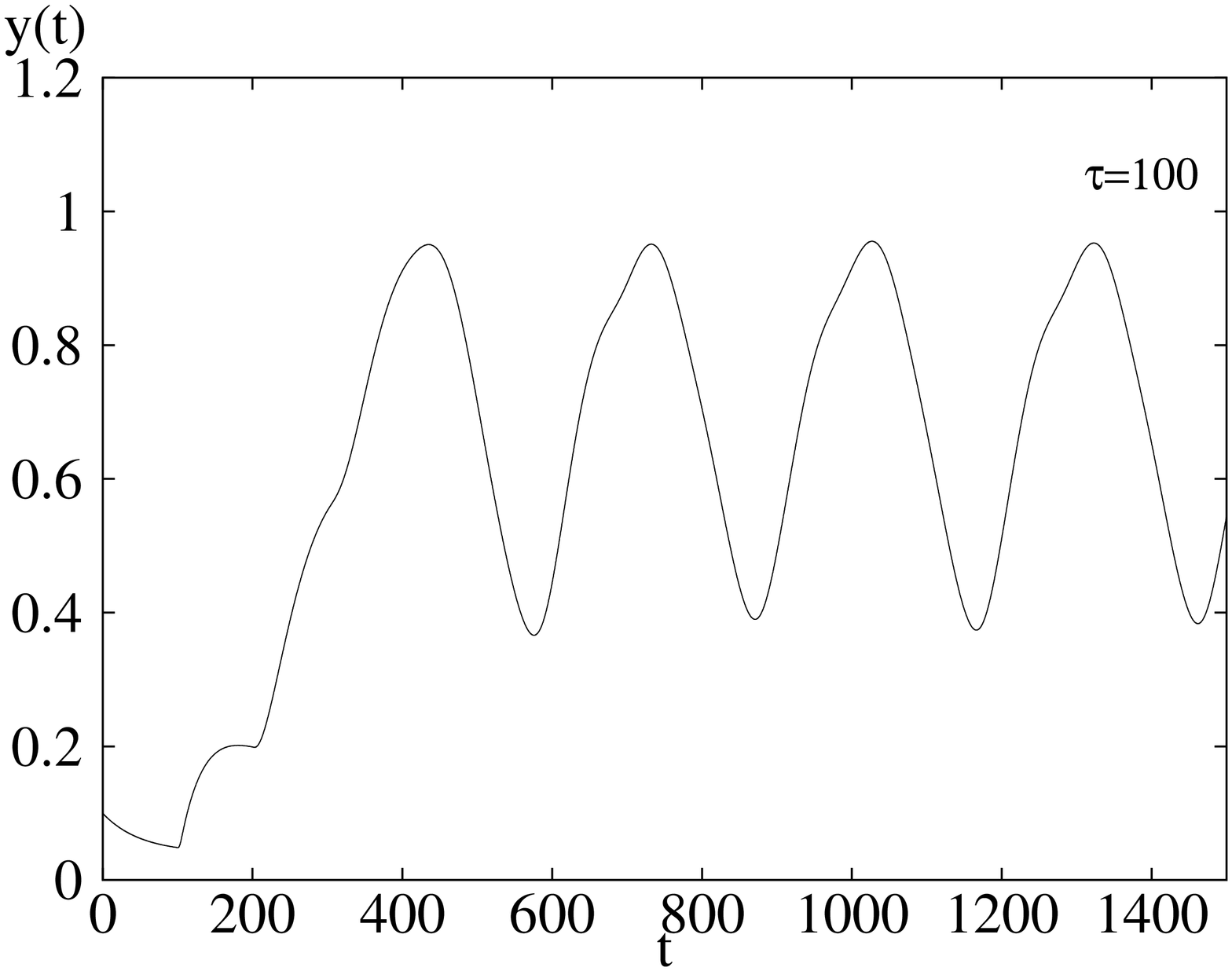}
\includegraphics[angle=270, width=5.cm]{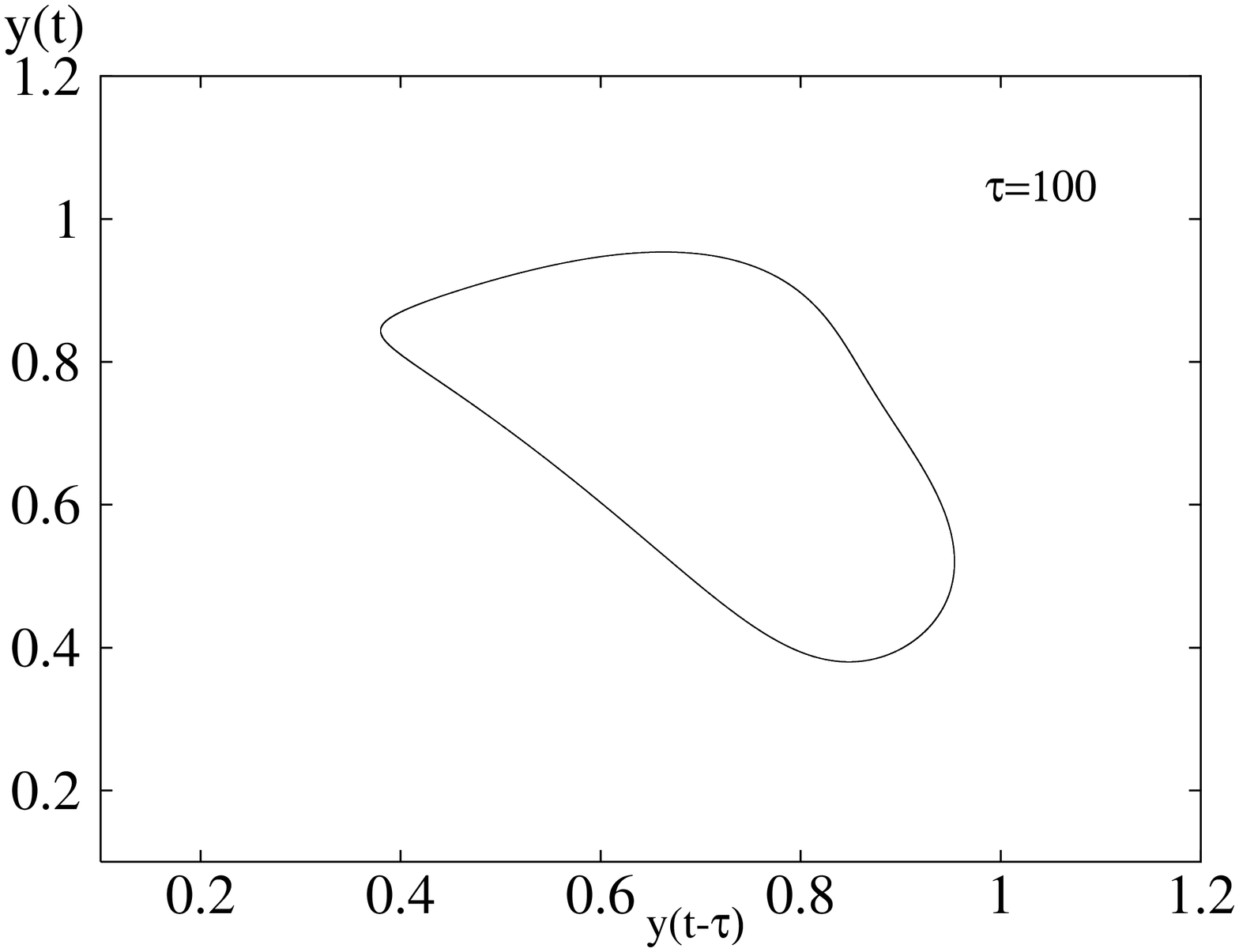}
   \end{center}
   \caption[ ] { {\scriptsize\bf (LEFT)} Time series starting  from the initial data $x(t)=y(t)=0.1, t\in[-\tau,0]$
   indicating how quickly the orbit gets close to the periodic attractor and
   {\scriptsize\bf (RIGHT)} time delay embeddings of
   the periodic attractors, demonstrating the sequence of  period
	 halfing bifurcations
   initiating  from the left for values of $\tau$ between  $91.5$ and $92$,
	 and at
   $\tau\approx 92.2, \  93.2, \ 98.3$. Graphs shown are for values of $\tau$  between these bifurcations:
    $\tau=91.95,\ 92,\
   93,\ 96  $, \ and \  $100$, \ with periods approximately equal to:
   $4794.3, \  2398.3, \ 1211.2,\  557.6, $ and  $295.3$, respectively. } \label{PD_right}
\end{figure}

Figure~\ref{zoom-in_80_100}  suggests that
 there are sequences of period doubling bifurcations,
 one initiating from the left at  $\tau\approx 83$, $86$,  $86.6$,
 $\dots,$
 and one initiating from the right  at  $\tau\approx  98.3$,  $93.2$, $92.2$,
 and $\tau$ between $92$ and $91.85$.
 To demonstrate these sequences,  time series ($y(t)$ versus
 $t$) and  time delay embeddings
 ($y(t)$ versus $y(t-\tau)$) at values of $\tau$ between these
 bifurcations are  shown in Figures~\ref{PD_left} and \ref{PD_right}.

\begin{figure}
    \begin{center}
	    \includegraphics[width=35ex]{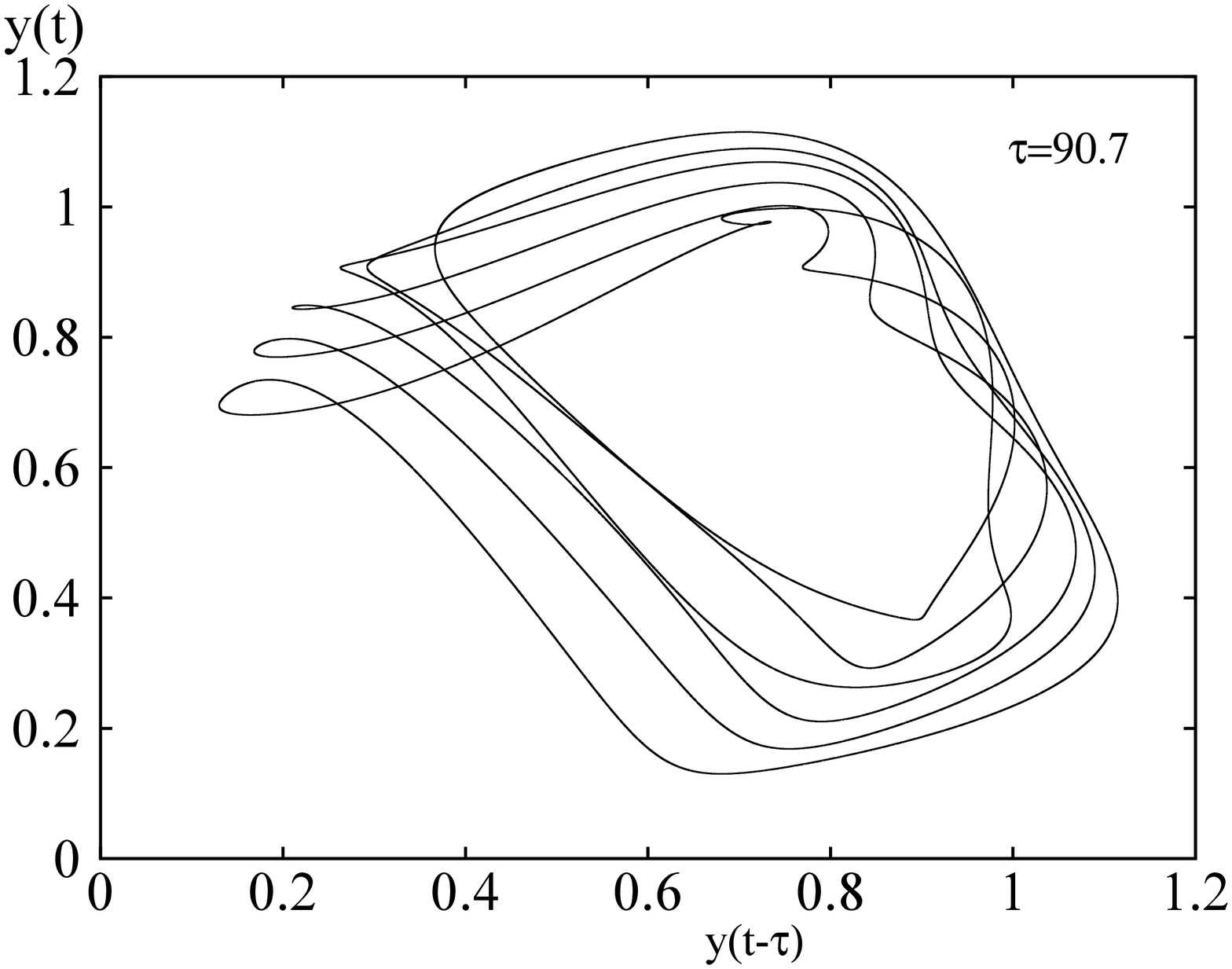}
  \end{center}
     \caption [two times three] { Time delay embedding starting at initial
     data $x(t)=y(t)=0.1, t\in[\tau,0]$ for $\tau=90.7$ showing a
   periodic attractor with period approximately
      $1800$, having $6(\neq2^n)$ loops for some integer $n$.}
     \label{period_two_times_three}
\end{figure}

    Figure~\ref{zoom-in_80_100} also suggests that between these sequences of period doubling
    bifurcations there is a window of values of $\tau$ at which there
     are  periodic attractors that do not have a
     period that  results from a bifurcation with period approximately
     equal to $2^n$  for some integer $n$, as well as chaotic dynamics.
     An example of the former is illustrated in
     Figure~\ref{period_two_times_three}.  For $\tau=90.7$, the
     time-series embedding of a periodic orbit with period approximately
    equal to  1800 time steps involving six loops (2 times 3)  is
    shown.

\begin{figure}[htbp!]
   \begin{center}
\includegraphics[angle=270, width=6cm]{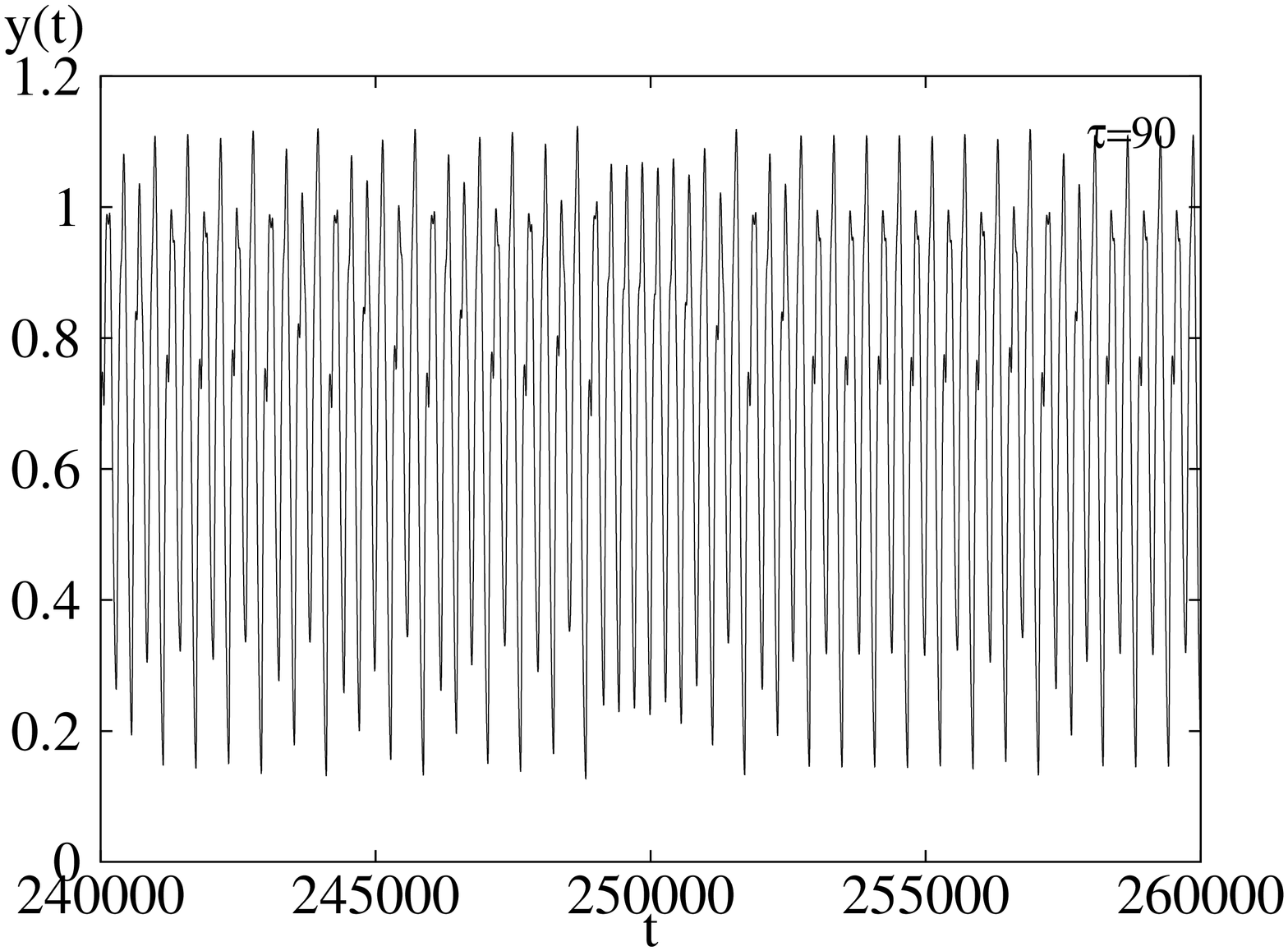}
\includegraphics[angle=270, width=6cm]{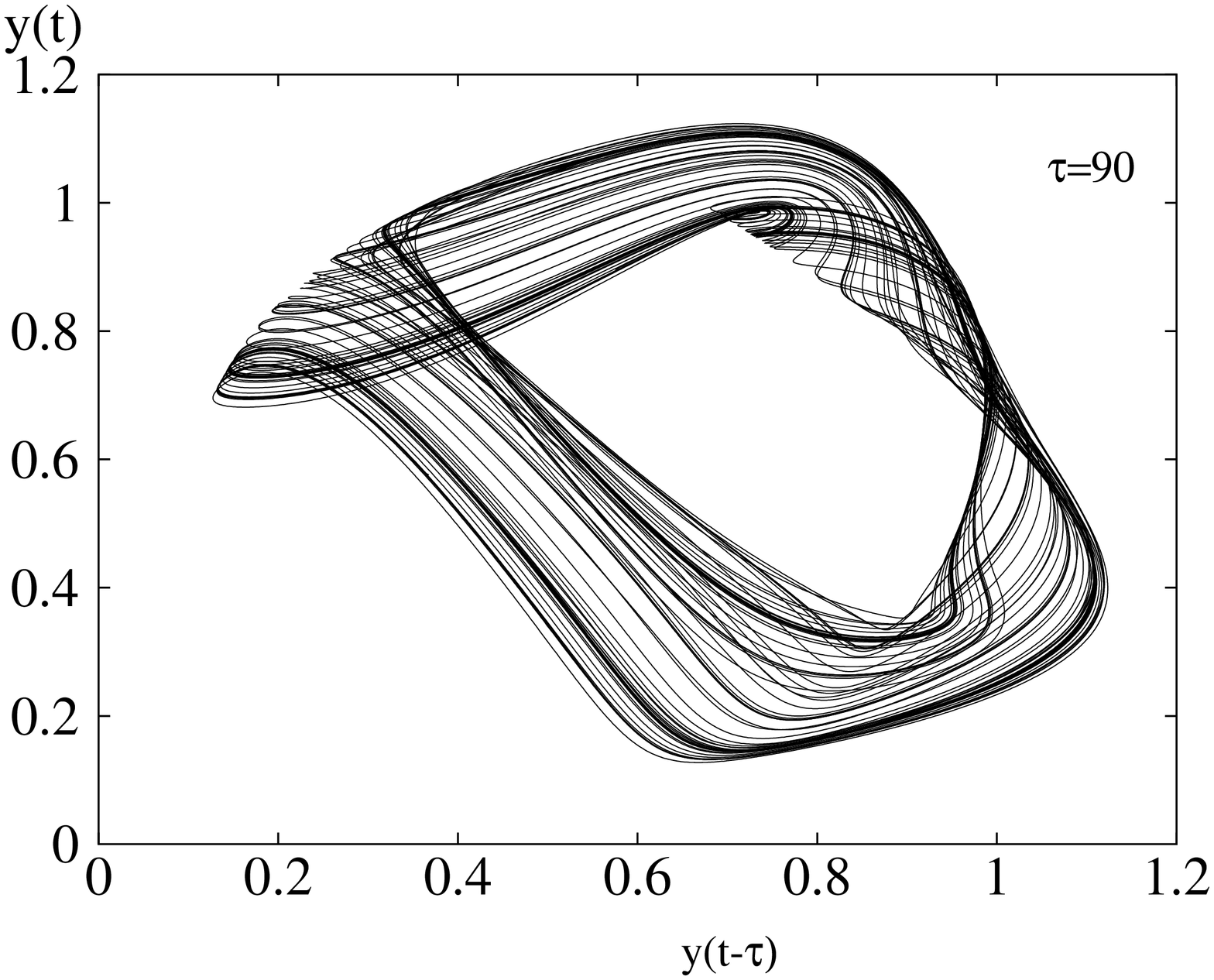}
   \end{center}
   \caption[ ]{ {\scriptsize\bf (LEFT)} Time series for $\tau=90$  starting  from the initial data
   $x(t)=y(t)=0.1, t\in[-\tau,0]$.
   {\scriptsize\bf (RIGHT)} Time delay embedding of
   the strange attractor for $\tau=90$.  Only the portion of the orbit from $t=240,000-260,000$
   is shown.} \label{chaotic_attractor}
 \end{figure}

\begin{figure}[htbp!]
	\begin{center}
	\includegraphics[angle=270, width=6cm]{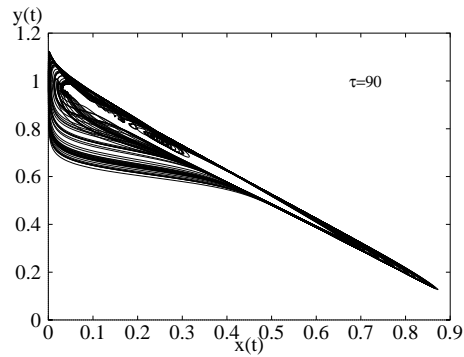}
	\caption[]{
	The strange attractor, for $\tau=90$,  shown in Figure~\ref{chaotic_attractor} in
	$(x,y)$-space.   Only the
   portion of the orbit from $t=240,000$ \,  to \, $260,000$ is shown.
   }\label{chaos_xy}
\end{center}
\end{figure}

\begin{figure}[htbp!]
	\begin{center}
%	\includegraphics[angle=270,
%	width=6cm]{Figures/pp_ts_tau_90_x_y.eps}
	\includegraphics[angle=0, width=7cm,height=5cm]{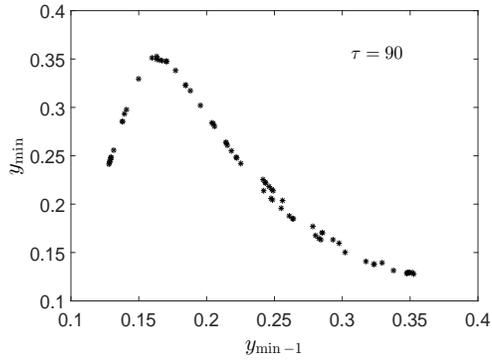} %\includegraphics[trim=10mm 30mm 10mm 30mm, clip, scale=0.75]{classical_rightsideup.pdf}
	\caption[]{% {\scriptsize\bf (LEFT)}
   %The strange attractor in $(x(t),y(t))$-space for $\tau=90$.  Only the
   %portion of the orbit from $t=240,000-260,000$ is shown.
   %{\scriptsize\bf (RIGHT)}
   The return map for $\tau=90$, computing the minimum value of $y(t)$
   %for $y(t)<0.7$
   as a function of
the preceding minimum value  of $y(t)$ for $y(t)<0.7$ in both cases, using the data in Figure~\ref{chaotic_attractor}. }\label{return_map}
\end{center}
\end{figure}

\begin{figure}[htbp!]
	\begin{center}
	\includegraphics[angle=0,width=6.75cm]{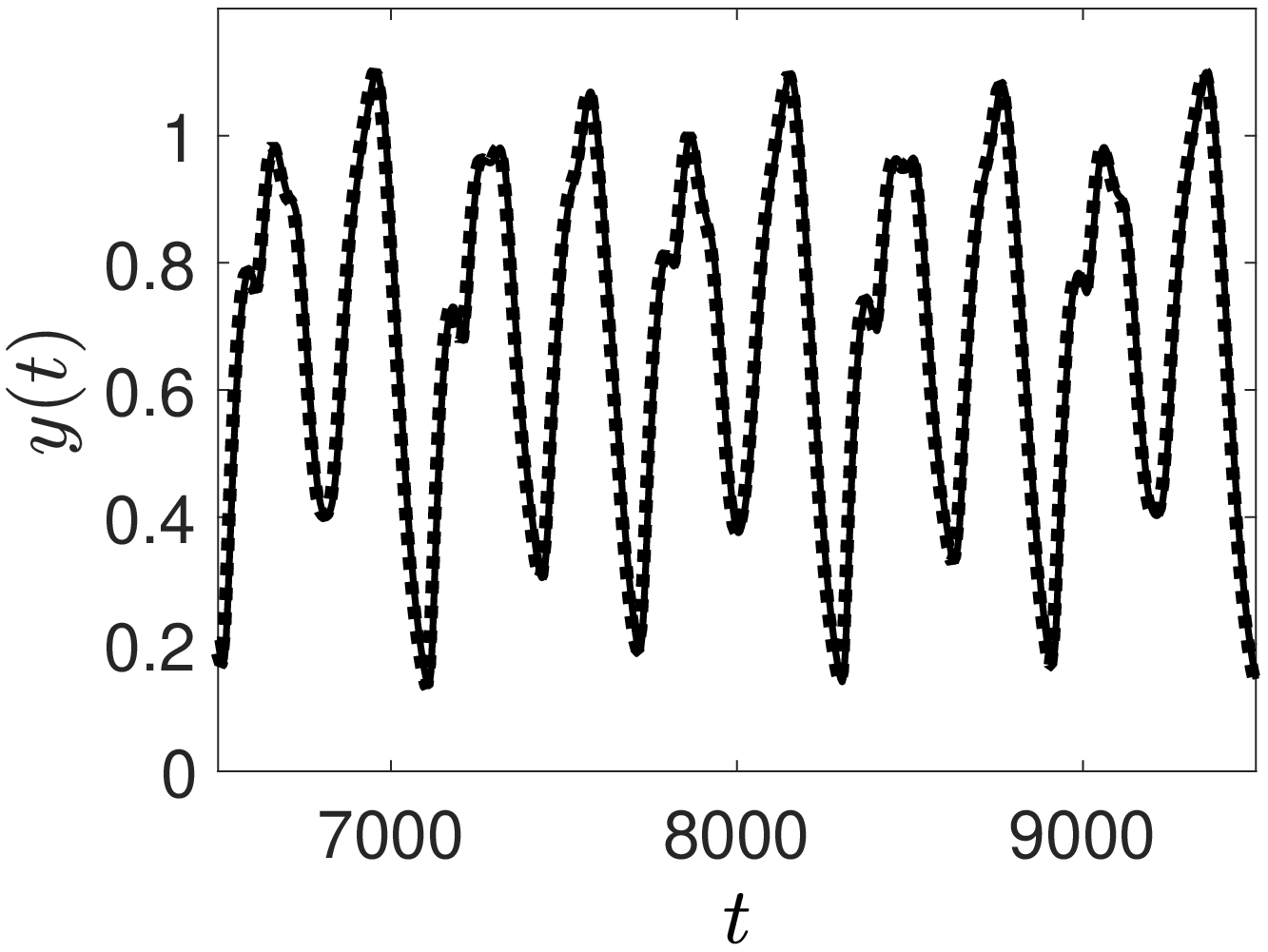}
		%tau_92_zoom_6500_9500_option2.eps}
	\includegraphics[angle=0,width=6.75cm]{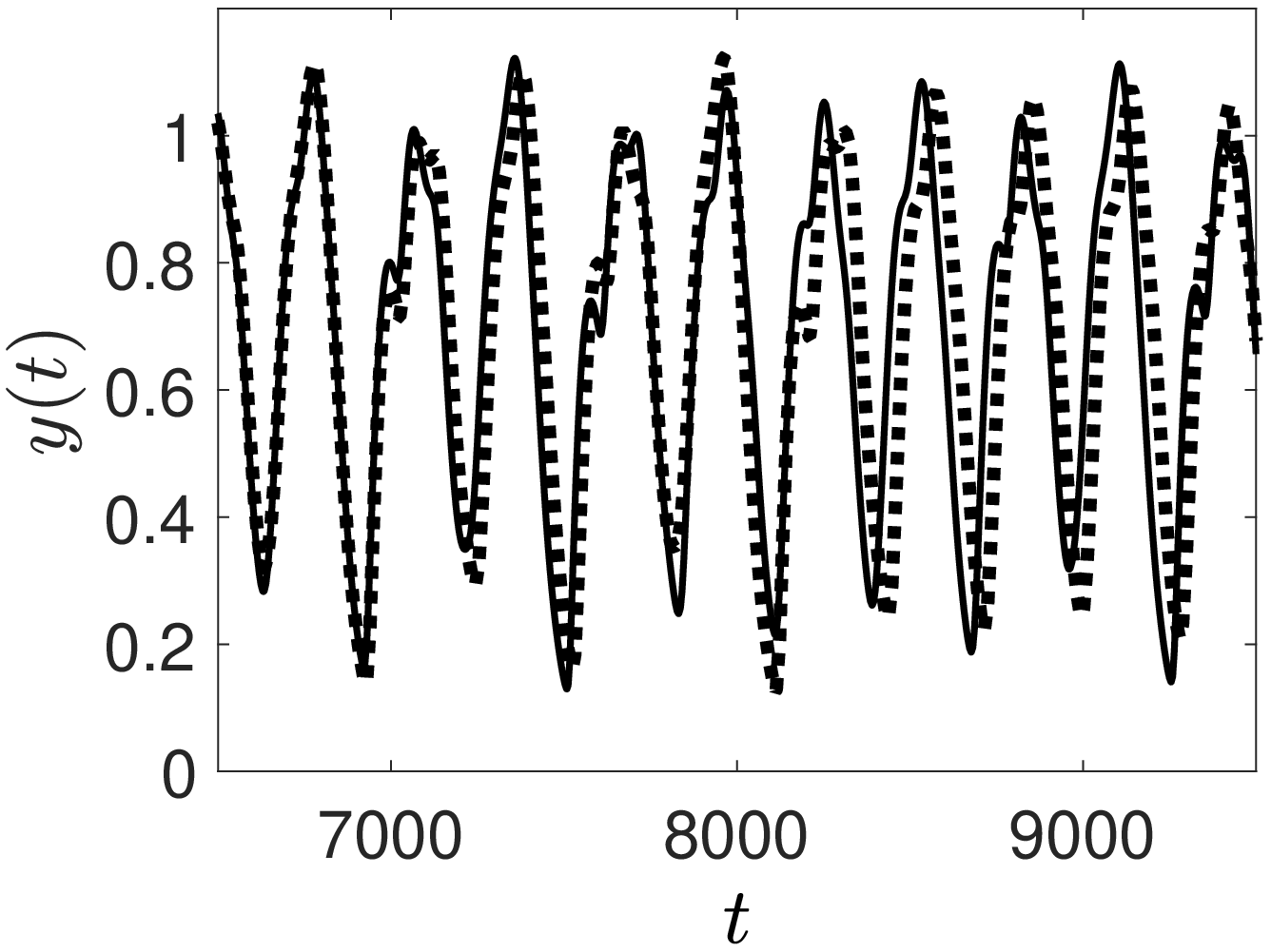}
		%tau_90_zoom_6500_9500_option2.eps}
	%\includegraphics[trim=10mm 30mm 10mm 30mm, clip, scale=0.75]{classical_rightsideup.pdf}
\caption[]{ Time series {\scriptsize\bf (LEFT)} for the solution that converges to a periodic attractor when $\tau=92$, and {\scriptsize\bf (RIGHT)}
for the solution that converges to a strange attractor when $\tau=90$, demonstrating that there is no sensitivity to initial data in the
former case, but that there is sensitivity in the latter case. Initial data used for the solid
curves:
$x(t)=y(t)=0.1$ for $t\in[-\tau,0]$, and for the dotted curves:  $x(t)=0.11$ and
$y(t)=0.1$ for $t\in[\tau,0]$.} \label{sensitivity}
\end{center}
\end{figure}

The time series and the time delay embedding of a chaotic attractor for
$\tau=90$ is shown in Figure~\ref{chaotic_attractor}. The projection of
this attractor into $(x(t),y(t))$-space is
also shown in Figure~\ref{chaos_xy}. This strange attractor resembles the
chaotic attractor of the well-known Mackey-Glass equation (\cite{Mackey:2009}, Figure 2).
The return map
shown in Figure~\ref{return_map}, for $\tau=90$, also resembles the return map for the
Mackey-Glass equation (\cite{Mackey:2009}, Figure 14) in the case of chaotic dynamics.
Sensitivity to initial data is a hallmark of chaotic dynamics.
Figure~\ref{sensitivity}{\scriptsize\bf (RIGHT)}
demonstrates that there is sensitivity to initial data in the case of the solution for
$\tau=90$ that
converges to the strange attractor,  shown in Figure~\ref{chaotic_attractor}.
To show that this is not just a numerical artifact, in
Figure~\ref{sensitivity}{\scriptsize\bf (LEFT)} we show that, as
expected,  there is no sensitivity for the solution for $\tau=92$ that
converges to the  periodic solution shown in Figure~\ref{PD_right}.

This example demonstrates that including delay in a simple
predator-prey model that always has a globally asymptotically stable
equilibrium point in the absence of delay,
cannot only destabilize a
globally asymptotically stable equilibrium point, but can  even result in
the birth of a strange attractor.

%\clearpage
\section{Discussion and Conclusions} \label{discussion}

We investigated the effect of the time  required for predators to process
their prey on the possible  dynamics predicted by a mathematical
model of predator-prey interaction.  We incorporated a discrete
delay    to model this process in one of the simplest classical
predator-prey models,    one that only allows convergence to an
equilibrium  when this delay is ignored.
We showed that including the delay results in a model with much
richer dynamics.  By choosing one of the simplest models when delay is
ignored, one that predicts that no   oscillatory behaviour is possible, the
effect of the delay on the dynamics is emphasized.

 This model can also be interpreted as a
model of a stage-structured population  with the
delay modelling the  maturation time of the juveniles
(see  Gourley and Kuang \cite{Kuang2004}).

In the model we considered,   the prey are assumed to grow logistically in
the absence of the predator.   The interaction of the predator and prey
is described using a linear response function, often referred to as mass
action or   Holling type I.
It is well-known that when delay is ignored this is one of the simplest
predator-prey models for which all solutions converge
to a globally asymptotically stable equilibrium point for all choices of the
parameters.   Therefore, any resulting non-equilibrium dynamics would
then be solely attributable to the introduction of the delay in the
growth term of the predator.  We not only found
non-trivial periodic solutions, but also bistability, and chaotic dynamics.
It is then likely that there is  similar   rich dynamics
in most predator-prey models with any reasonable response function, when  such a delay is incorporated,  for
some selection of the parameters, including the form most used by
ecologists, the  Holling type II form.  This form
given mathematically by $f(x)=mx/(1+bx))$, can be considered a   generalization of the Holling type I form, obtained by simply adding an extra parameter, $b$.
However, we feel that demonstrating  that this wide range of dynamics is
possible for even for one of the simplest models gives more compelling
evidence that delay should not be ignored when making policy decisions.

Understanding how changes in average
temperature   might affect   survivability of endangered populations or
result in invasions by undesirable populations is   important.
Since temperature can affect how quickly predators process the prey
that they capture, based on our results there might be important
implications   for  populations in the wild.   In most predator
populations,   the processing time, $\tau$, is faster when it is warmer and
slower when it is colder.  Our results may help us  understand how a change in
average temperature might
influence the dynamics of   particular predator-prey systems of interest.

Our study
suggests that we need to be careful when measuring population sizes
in the wild and predicting the general health of the population based
upon whether the population seems to be increasing or decreasing.
Short term indications  that a population size is changing
 in a system with oscillatory dynamics  may be misleading
and predicting future population size may be impossible without more
information. It would be necessary to have some idea of the period of the
intrinsic oscillations  if the population is suspected to
vary periodically.  If the dynamics are suspected to be chaotic, this  may be
even  more complicated, due to sensitivity to initial data.

If a predator-prey system has the potential to have  chaotic dynamics,
based on our results, is there anything that is predictable?    Can such
analyses suggest how to prevent   extinctions or
invasions due to a change in average temperature that could result in a
change in the processing time of the prey by the predator? We give some observations based on the predictions of our model and the example we considered in Section~\ref{sec:example}.  However,  more work would need to be done to determine whether  these predictions are consistent for  more realistic
models,  and if so,   long term observations would have to be made by
ecologists  to determine if they are relevant for   populations
in the wild.

  First,     at the one extreme,
if the processing time is too long, the predator population would not
be expected to survive,  since it
is obvious that if the processing time is longer than the life span of
the predator,  the predator population has no chance to avoid
extinction.   If there was
no such  threshold for extinction predicted by the model, the model should be abandoned.  It is
therefore very important to use the term
$e^{-s\tau}y(t-\tau)$ and not just $y(t-\tau)$ in order to account for the predators that do not survive
long enough to affect growth of the population,  in the equation
describing the growth of the predator in
model (\ref{ppreydelay}). Due to this term,
 in our model there is such a threshold, $\tau_c$.

 From the orbit diagram (see Figure~\ref{bif_diag_full}),
for the example considered in Section~\ref{sec:example}, we summarize some
observations.  The dynamics for both populations are   oscillatory for a wide
range of processing times,  and non-oscillatory for only relatively
(very) short or relatively long processing times  (i.e.  before the first Hopf bifurcation at
$\tau\approx 1.9 $ and after the last one at $\tau\approx 108 $).
For relatively long processing times, slightly larger than
the value of $\tau$ at the final Hopf
bifurcation of the coexistence equilibrium,  the population is no longer
oscillatory.  The size of the predator
population decreases relatively slowly as the processing time increases
further.
Although  the size of the predator population gets smaller as the
processing time increases, it does not get  much smaller,
until the processing time gets close to the threshold  for extinction
$\tau_c$.
If we had shown the diagram extended to the threshold $\tau_c=170$, one would
see that as the processing time gets close to  $\tau_c$, the size of the
predator population suddenly decreases
relatively quickly to zero.  So our model predicts that for a predator
population with fairly long processing times to
begin with,    cooling of the environment could be expected to be
detrimental with respect to the survivability of the predator
population. Thus, this effect on the size of the predator population might be minor if the delay is close to its Hopf
bifurcation value, but  could be  drastic if it  is close to the extinction value.

Our example also suggests that   the
predator population may not be oscillatory when it becomes
endangered and hence close to extinction, i.e., for excessively long
processing times. It would  be interesting to
investigate  if this also holds for the  model,
with Holling type I response functions replaced by Holling type II.
In the   model with Holling type II response functions, if the  carrying capacity of the environment for
the prey is
not affected by the cooling and it is relatively high, then the model
predicts that the predator
population could still be   oscillatory, even if  the processing
time for the predator is ignored.  However,   cooling might also be expected to
reduce  the carrying capacity for the prey, moving the parameters to a
range where that   model would also predict non-oscillatory dynamics (the paradox of enrichment \cite{Rosenzweig1971}).  So once again, perhaps alarm
bells should be sounded   when  there is cooling and the
predator population has been non-oscillatory and  appears to decline rapidly
as the average yearly temperature declines.

On the other hand, at the other extreme,  our model predicts a Hopf
bifurcation at a relatively small value of the delay, resulting in
the birth of a family of periodic orbits  with
amplitude  increasing very quickly as the delay
increases with
both the prey
and predator populations spending time very close to zero.
This remains the case for  small, but
intermediate values of the delay,  (before the two  saddle node of limit
cycles bifurcations near $\tau=80$).
For this range of $\tau$, these populations are therefore   very
susceptible to stochastic extinctions.  If the processing time was
originally very small (below the critical value for the first Hopf
bifurcation)  and cooling made it longer, again a stochastic extinction might be
likely.  Similarly, if $\tau$ was close to the first of the
two  saddle node of limit
cycles bifurcations near $\tau=80$), warming of the average temperature
could result in a stochastic extinction  of one  or both of the populations.
Since the predator cannot survive without the prey in our model, even if
it was the  prey population  that experienced the stochastic
extinction, the predator population would eventually die out as well.
Between the smallest value of $\tau$ at which there is a period doubling
bifurcation, and the value of $\tau$ at the largest Hopf bifurcation,
cooling would probably be advantageous to the predator population size.

In summary, it seems that whether  cooling is beneficial or detrimental to
  the size of the predator population depends on where the delay is on
  the bifurcation diagram, and it is very likely that
  this is  very difficult to determine. As well,
  if   values of $\tau$ were to lie in the chaotic region, then changes
  in the environment that changed the value of $\tau$ to obtain regular
  oscillatory dynamics may or may not be preferable.  As well,
 as can be seen by the various
attractors  shown in Figures~\ref{PD_left}-\ref{chaos_xy}, and the orbit
diagrams in Figure~\ref{zoom-in_80_100},  certain
properties of the system in the chaotic region such as
the maximum and minimum values of the predator
size  were fairly insensitive to the change in the delay.
 This is only a toy model. However, it suggests that ignoring delay
in a model can result in incorrect predictions. Here, delay could change
the dynamics from convergence to 
a globally asymptotically stable equilibrium to wild
oscillations.  It is most likely that in nature, it would not be
possible to distinguish from data, a priori, if the dynamics were chaotic or
periodic, but  more importantly it would not necessarily be predictable
what the effect of an increase or a decrease in the delay would be. Hence,
our analysis indicates that one should be extra cautious if trying to
manipulate the delay to achieve a certain result based on model
predictions.  

Finally it is worth pointing out that
   the resulting strange attractor in the predator-prey model studied
here bares such a close resemblance to the Mackey-Glass attractor, a model
involving  a single delay differential equation
to model a simple  feedback
system for  respiratory control or hematopoietic diseases. Understanding
whether there is a deeper significiance to this relationship may give us
a better understanding of the class of possible strange attractors
and warrants further investigation.

%%

%\clearpage

%\noindent  {\bf References}
%\bibliographystyle{plain}
%\bibliography{mybib}
\bibliographystyle{siam}

%\bibliography{mybib}

\section*{Funding}
The Research of Gail S. K. Wolkowicz was partially supported by
 Natural Sciences and Engineering (NSERC) Discovery Grant \# 9358 and
 Accelerator supplement.

\appendix
\section{Proofs}
\subsection{Proof of Proposition~\ref{bounded}} \label{app:bounded}
\noindent
{\bf Proof.}
Parts 1 \& 2. Assume that $(\phi(t),\psi(t))\in X$. Since the
right hand side of  (\ref{dimenless}) is Lipschitz, there exist $h>0$
such that solutions exist and are unique for  all $t\in[0,h]$.
We will show that solutions are bounded and hence do not blow up in
finite time. Therefore, by
the standard results for  existence and uniqueness of solutions for delay
differential equations in Driver~\cite{Driver62}, it will follow that
solutions exist  and are unique for all $t\geq 0$.

 Since
$\phi(0)\geqslant 0$, and
the face where $x(t)=0$ is invariant,   that $x(t)\geqslant 0$ for all
 $t\geqslant 0$ follows by uniqueness of solutions in forward time for initial
value problems.  That $y(t)\geqslant 0$ for all $t>0$, follows since
$y'(t)\geqslant -s y(t)$ for all $t\in[0,\tau]$ and so $y(t)\geqslant 0$ for all
$t\in [-\tau,\tau]$.  Arguing inductively on  $[n\tau, (n+1)\tau], \,
n=1,2,\dots$, it follows that $y(t)\geqslant 0$ for all $t>0$, where
the  solution  exist.

 To show that  solutions are bounded above, consider the first equation of (\ref{dimenless})
\begin{equation*}
\dot x(t)=x(t)(1-x(t))-y(t)x(t) \leqslant x(t)(1-x(t)).
\end{equation*}
It is well-known that for  the logistic equation
$\dot z(t)=z(t)\left(1-z(t)\right)$,
given any $\epsilon_0>0$, there exists $T>0$, such that $|z(t)|<1+\epsilon_0$
for all $t\geqslant T$.
Using a comparison principle (e.g. \cite{LakI} Theorem~1.4.1, page 15)
 $x(t)\leqslant z(t)$,
and so $0\leqslant x(t)<1+\epsilon_0$ for all $t\geqslant T$.

To prove  that $y(t)$ is bounded above,  define
\begin{equation}
w(t)=Y e^{-s\tau} x(t-\tau) + y(t).\label{doubleu}
\end{equation}
Then
\begin{equation*}
\begin{aligned}
  \dot{w}(t)
&= Y e^{-s\tau}\frac{dx(t-\tau)}{dt}+\frac{dy(t)}{dt},\\
&= -s y(t) + Y e^{-s\tau}x(t-\tau)\left(1-x(t-\tau)\right),\\
&= -s w(t) + Y e^{-s\tau} x(t-\tau)\left(s+1-x(t-\tau)\right), \\
&\leqslant -s w(t)+\frac{1}{4}Ye^{-s\tau}(s+1)^2,
\end{aligned}
\end{equation*}
since $\left(x(t-\tau)-\frac{s+1}{2}\right)^2\geqslant 0$
%\begin{equation*}
implies  that $x(t-\tau)\left(s+1-x(t-\tau)\right)\leqslant
\frac{(s+1)^2}{4}$.
%\end{equation*}
Therefore, since
$z(t)=z(0)e^{-st}+\frac{1}{4s}Ye^{-s\tau}(s+1)^2(1-e^{-s t})$ is the
solution of the initial value problem
\begin{equation*}
\dot z(t)= -sz(t)+\frac{1}{4}Ye^{-s\tau}(s+1)^2, \qquad
z(0)=w(0)\geqslant 0,
\end{equation*}
using a comparison principle it follows that $w(t)\leqslant z(t)$ for
all $t>0$.
Consequently, by (\ref{doubleu}), $y(t)\leqslant w(t)\leqslant
w(0)e^{-st}+\frac{1}{4s}Ye^{-s\tau}(s+1)^2(1-e^{-st})$.
%\begin{comment}
%By (\ref{doubleu}),
%\begin{equation*}
%\begin{aligned}
%Y e^{-s\tau} x(t-\tau) + y(t) =w(t)
%&\leqslant w(0)+\frac{1}{4s}Ye^{-s\tau}(s+1)^2 \\
%&\leqslant \left(Y e^{-s\tau} x(-\tau) + y(0)\right)+\frac{1}{4s}Y
%e^{-s\tau}\left(s+1\right)^2.
%\end{aligned}
%\end{equation*}
%\end{comment}
Therefore, $y(t)$ is bounded.
\medskip

Part 3.  Assume that $(\phi(t),\psi(t))\in X^0$. Since $\phi(0)>0$, and
the face where $x(t)=0$ is invariant, it follows that $x(t)>0$ for all
 $t\geqslant 0$.  Since there exists $\theta\in[-\tau,0]$ such that
 $\phi(\theta)\psi(\theta)>0$, take $T=\tau+\theta$, and note that
 $0\leqslant T\leqslant\tau$.  Since $y(T)\geqslant 0$,
 either $y(T)>0$ or $y(T)=0$ and by (\ref{dimenless}),
 $y'(T)=-sy(T)+Ye^{-s\tau}y(\theta)x(\theta)= 0
 +Ye^{-s\tau}\psi(\theta)\phi(\theta) >0$. In both cases,
  there exists $\epsilon>0$ such that $y(t)>0$ for all
 $t\in(T,T+\epsilon]$.
%From (\ref{dimenless}),
%$y'(t)\geqslant  -sy(t)$ and for all $t< T$, since
%$t-\tau\in[-\tau,0]$ if $t\leqslant T$ and so $x(t-\tau)y(t-\tau)
%\geqslant 0$.
%Therefore, $y(t)\geqslant \psi(0) e^{-s t}\geqslant 0$ for all $t\in[0,T]$
%and  $y'(T)=-sy(T)+Ye^{-s\tau}y(\theta)x(\theta)>0.$ Hence, there
%exists $\epsilon>0$ such that $y(t)> 0 $  for all $t\in (T,T+\epsilon]$.
But this implies that $y'(t)\geqslant -sy(t)$ for all
$t\in[-\tau,T+\epsilon+\tau]$.    Therefore, $y(t)>0$ for all
$t\in[T+\epsilon,T+\epsilon+\tau]$.  By repeating this argument, it follows
that $y(t)>0$ for all $t>T+\epsilon$ for any $\epsilon>0$.
$\hfill\Box$

\subsection{ Proof of Theorem~\ref{th:gasEi}} \label{app:gasEi}

\noindent
{\bf Proof.} Part 1. Evaluating  (\ref{chareqn})
at  $E_0$ gives
%\begin{equation*}
$P(\lambda)|_{E_0}=(\lambda+s)(\lambda-1)=0,$
%\end{equation*}
which has two real roots $\lambda=-s$ and $\lambda=1$.
Therefore, $E_0$ is a saddle.

%$\hfill\Box$
%\begin{theorem}%theorem 3
%Consider (\ref{dimenless}). Equilibrium $E_1$ is unstable if
%$0\leqslant \tau <\tau_{c}$ and globally asymptotically stable if
%$\tau> \tau_{c}$.
%\end{theorem}

%\noindent
%{\bf Proof.}

Part 2.(a)  Evaluating (\ref{chareqn}) at $E_1$ gives
% \begin{equation*}
$ P(\lambda)|_{E_1}=(\lambda + 1)(\lambda+s -
Ye^{-(s+\lambda)\tau})=0.$
%\end{equation*}
One of the roots is $\lambda=-1$. The other roots satisfy
%\begin{equation}
$g(\lambda,\tau) := (\lambda + s)e^{(\lambda + s)\tau} = Y.$
%\end{equation}
 For any fixed $0 \leqslant \tau < \tau_c$,  there is a
real root $\lambda(\tau)>0$, such that $g(\lambda(\tau),\tau)=Y$, since $g(0,\tau)<Y$ and $g(\lambda,\tau)\rightarrow
\infty$ as $\lambda \rightarrow \infty$.
%The left hand side of (\ref{xepx}) is a monotone
%increasing function in $\lambda$ for any fixed $\tau$. It takes the
%value $se^{s\tau}$ at $\lambda=0$ and tends to positive infinity as
%$\lambda\rightarrow +\infty$. Since $0\leqslant\tau<\tau_c$,
%$se^{s\tau}<Y$. By the Intermediate Value Theorem, there exists a
%unique $\lambda(\tau)>0$ such that equation (\ref{xepx}) holds and
%so $P(\lambda)|_{E_1}=0$ has at least one positive root
%$\lambda(\tau)$.
Hence, $E_1$ is unstable.
%for $0\leqslant\tau<\tau_c$.

Part 2.(b) Next assume that $\tau>\tau_c$. We prove that $E_1$ is globally
asymptotically stable.  Since   $\frac{se^{s\tau}}{Y}>1$,
$\epsilon_0:=\frac{1}{2}\left(\frac{s}{Y}e^{s\tau}-1\right)>0$, and so
 by Proposition~\ref{bounded}, there exists a $T>0$ such that $x(t)<
1+\epsilon_0$ for all $t>T$. Therefore, for all sufficiently large $t$,
$
Y e^{-s\tau} x(t-\tau) <Y e^{-s\tau}(1+\epsilon_0)  = Y e^{-s\tau}\left(1+\frac{1}{2}\left(\frac{s}{Y}e^{s\tau}-1\right)\right)
=\frac{1}{2}Y e^{-s\tau} +\frac{s}{2}<s,
$
and so the second equation of (\ref{dimenless}) can be written
$\dot y(t)=-s y(t) + b(t)y(t-\tau),$
where  $b(t):=Y e^{-s\tau} x(t-\tau)<s$. Choosing $\alpha=s/2$
in Kuang \cite{Kuang}, (Example 5.1, Chapter 2, page 32), since
$b^2(t)<s^2=4(s-\alpha)\alpha=4(s-s/2)(s/2)$,
$y(t)\rightarrow 0$ as $t\rightarrow \infty$. Hence, for any
$\epsilon> 0$, there exists $T_1$ such that $0 < y(t) < \epsilon$ for $t >
T_1$. From the first equation of (\ref{dimenless}), for any
$0<\epsilon<1$,
$x(t)\left(1-x(t)-\epsilon\right) < \dot x(t) <
x(t)\left(1-x(t)\right).$
Note that all solutions of $\dot z(t) = z(t)(1-z(t))$ converge to
$z(t)=1$ and all solutions of $\dot z(t) = z(t)(1-z(t)-\epsilon)$
converge to  $z(t)=1-\epsilon$  as $t\rightarrow \infty$.
By a standard comparison principle, for any solution $x(t)$ of
(\ref{dimenless}), $x(t)\rightarrow 1$ as $t\rightarrow \infty$.
Therefore, $E_1$ is globally asymptotically stable.

Part 3. Existence of $E_+$ follows from (\ref{notation1}).

Part 3.(a) When $\tau=0$, model (\ref{dimenless}) reduces to a well
studied model involving only ordinary differential equations. That $E_+$ is globally
asymptotically stable in this case is well known and can easily be
proved using phase plane analysis and  the Bendixson-negative criterion
to rule out periodic solutions.

Part 3.(b)  The proof is similar to the approach used in Chapter 5.7 of
Smith and Thieme \cite{persistencebook}. The predator
reproduction number at a constant level of $x$ is given by
$\cal{R}(x)=\frac{Ye^{-s\tau}}{s}x.$   Then,  since $E_+$ exists,
\begin{equation} \label{rep_number}
  \cal{R}(1)>1 \quad \mbox{and} \quad  \cal{R}(x_+)=1.
\end{equation}

Let,
$x^{\infty}\equiv\limsup_{t\rightarrow \infty} \, x(t)$ and
$y^{\infty}\equiv\limsup_{t\rightarrow \infty} \, y(t)$.
Then, we claim
\begin{equation} \label{eq:min}
x^\infty\geq \frac{s e^{s\tau}}{Y}.
\end{equation}
Suppose not, i.e., that
$x^{\infty}$ is smaller than this minimum.  Applying the fluctuation
lemma (see Hirsch et. al. \cite{HHG1985} or Smith and Thieme
\cite{persistencebook})
  to the $y$ equation in (\ref{dimenless}),  there exists a monotone
  increasing  sequence of
  times $\{t_n\} \rightarrow \infty$ as $n\rightarrow \infty$ such that
  $ \dot{y}(t_n)= 0$ and $y(t_n)\rightarrow  y^{\infty}$ as
  $n\rightarrow \infty$.
  Therefore, $0= \dot{y}(t_n)=
  -sy(t_n)+Ye^{-s\tau}x(t_n-\tau)y(t_n-\tau)$. Letting $n\rightarrow
  \infty$, we obtain
$0\leq  (-s+Ye^{-s\tau} x^{\infty})y^{\infty}$, since
$\limsup_{n\rightarrow \infty} x(t_n-\tau)y(t_n-\tau) \leq
x^{\infty}y^{\infty}$.  The term in the
brackets in this inequality is negative, since we are assuming
(\ref{eq:min}) does not hold. It follows that
$y^{\infty}=0$.  But then, applying the fluctuation lemma to the $x$ equation
in (\ref{dimenless}), it follows that   $x^{\infty}=1>\frac{s e^{s\tau}}{Y}$,
contradicting our assumption that (\ref{eq:min}) is not satisfied.
Hence, (\ref{eq:min}) holds.

Let $\Phi:\mathbb{R}_+\times X\rightarrow X$ denote the semiflow
generated by  model (\ref{dimenless}).   Consider  the function
$\rho: X\rightarrow \mathbb{R}_+$,  defined for $(\phi(t),\psi(t))\in
X^0$ by  $\rho((\phi,\psi))=\phi(0)$ so that
$\rho(\Phi(t,(\phi,\psi)))=x(t)$.  By part 3 of
Proposition~\ref{bounded},
$\rho((\phi,\psi)))>0$ implies
that $\rho(\Phi(t,(\phi,\psi)))>0$ for all $t\geq 0$.  Therefore, the
semi-flow is uniformly weakly $\rho$-persistent. By a similar argument
to that given in Theorem 5.29 of \cite{persistencebook} together with
Proposition~\ref{bounded},  $\Phi$ is a continuous semiflow that has a compact attractor of bounded
sets. Therefore,  $\phi$ is uniformly $\rho$-persistent by Theorem 5.2 of
\cite{persistencebook}, and hence there exists $\epsilon_1>0$ such that
$x_\infty\equiv \liminf_{t\rightarrow \infty} \, x(t) \geq \epsilon_1$, for all
solutions with $x(0)>0$.

Before, we show that if $(\phi,\psi)\in X^0$, there exists $\epsilon_2>0$
such that $\liminf_{t\rightarrow \infty} y(t)>\epsilon_2$, we use Laplace
transforms to show that
if $(\phi,\psi)\in X^0$, then
%if $y(\bar{t}>0$ for some $\bar{t}\geq 0$, then
$\cal{R}(x_\infty) \leq 1$.
We denote the  Laplace transform of a function  $f(t)$ as
$\widehat{f}(\lambda)=\int_0^{\infty} \, e^{-\lambda t} f(t) \, dt$,
and note that the Laplace transform of $y(t)$
exists for all $\lambda>0$,   since $y(t)$ is bounded by
Proposition~\ref{bounded}. Taking the Laplace transform on both sides of
the $y$ equation of (\ref{dimenless}) and simplifying we obtain:
$$(\lambda+s)\widehat{y}(\lambda)=y(0)+Y e^{-(s+\lambda)\tau}
\widehat{xy}(\lambda)+Y e^{-(s+\lambda)\tau}\int_{-\tau}^0 \,
e^{-\lambda \tau} x(t) y(t) \, dt.$$
There exists $\delta>0$ such that $x(t)\geq (x_{\infty}-\delta)$ for all
$t\geq0$, and so
$\widehat{xy}(\lambda)\geq(x_{\infty}-\delta)\widehat{y}(\lambda).$
$$(\lambda+s)\widehat{y}(\lambda)\geq   Y
e^{-(s+\lambda)\tau}(x_{\infty} -\delta)\widehat{y}(\lambda).$$
By part 3 of Proposition~\ref{bounded}, $y(t)>0$ for all sufficiently large
t, and hence $\widehat{y}(\lambda)>0$, we can divide by
$\widehat{y}(\lambda)$ to obtain
$$(\lambda+s)\geq   Y
e^{-(s+\lambda)\tau}(x_{\infty} -\delta).$$  After a shift of time, if
necessary, we can take the limit as
$\delta, \lambda \rightarrow 0^+$ to obtain
$s \geq Ye^{-s\tau}x_{\infty}.$  This is equivalent to
\begin{equation} \label{Rxinfty}
  \cal{R}(x_{\infty})\leq 1.
\end{equation}

 Define
$\rho:X\rightarrow\mathbb{R}_+$ by
$\rho((\phi,\psi))=\min\left\{\phi(0),\psi(0)\right\}$.\\
Then, $\rho(\Phi(t,(\phi,\psi)))=\min\{x(t),y(t)\}.$
Suppose $(\phi,\psi)\in X^0$ and  $y(t)$ is not uniformly persistent.
By part 3 of Proposition~\ref{bounded}, shifting time if necessary, there
is no loss of generality if we assume that $\psi(0)=y(0)>0$.
Therefore, $\rho(\Phi(t,(\phi,\psi)))>0$ for all $t\geq 0$.
Since by Proposition~\ref{bounded}, $\Phi$ has a compact attractor of
bounded sets, by Theorem~5.2 in  \cite{persistencebook}
 we need only show that $\Phi$ is uniformly weakly $\rho$-persistent.
 Suppose not.  Recall that we have
already shown that $x_{\infty}>\epsilon_1>0.$ Take any
$\epsilon\in(0,\epsilon_1)$.  Then, there exists a solution with $x(0)>0$ and $y(0)>0$
 such that $y^{\infty}  < \epsilon.$
 Apply the fluctuation lemma to the $x$ equation of (\ref{dimenless}).
 Therefore, $0>x_{\infty}(1-x_{\infty}-\epsilon).$  Taking $\epsilon>0$
 sufficiently small,  since
   $\cal{R}(1)>1$ by  (\ref{rep_number}) and $\cal{R}(x)$ is
 increasing,  it follows that $\cal{R}(x_\infty)>1$,
 contradicting (\ref{Rxinfty}).
$\hfill\Box$

\subsection{Proof of Lemma~\ref{lem:imaginarycross}: Roots of the characteristic equation cannot bifurcate in
from infinity} \label{app:imaginarycross}
\noindent
{\bf Proof.} In Kuang \cite{Kuang} (Theorem~1.4, Chapter 3, page 66),
taking $n=2$ and
$g(\lambda,\tau)=p(\tau)\lambda+(q\lambda+c(\tau))e^{-\lambda\tau}+\alpha(\tau),$
since
\begin{equation*}
\limsup_{\mathrm{Re}\lambda>0, |\lambda|\rightarrow
\infty}|\lambda^{-2} g(\lambda,\tau)|= 0<1,
\end{equation*}
no root of (\ref{CharE+}) with positive real part can enter from
infinity as $\tau$ increases from $0$. Since, when $\tau=0$  all roots
have negative real parts, the result follows.
%roots with positive
%real part can only appear by crossing the imaginary axis
%as $\tau$ increases.
$\hfill\Box$

\subsection{Proof of
Theorem~\ref{omega+posicond}}\label{app:omega+posicond}
\noindent
{\bf Proof.} Using the quadratic formula to solve  for $\omega^2$ in
(\ref{omegapower4}),  the roots must satisfy
$$%\begin{equation}
\omega_{\pm}^2=\frac{1}{2}\bigg(q^2-p^2(\tau)+2\alpha(\tau) \pm
\sqrt{(q^2-p^2(\tau)+2\alpha(\tau))^2
-4\left(\alpha^2(\tau)-c^2(\tau)\right)}\bigg).\label{omega+-2}
$$%\end{equation}
Since
\begin{equation*}
q^2-p^2(\tau)+2\alpha(\tau) =
s^2-s^2\left(1+\frac{e^{s\tau}}{Y}\right)^2+2\frac{s^2 e^{s\tau}}{Y}
= -\left(\frac{se^{s\tau}}{Y}\right)^2< 0,
\end{equation*}
$\omega_-^2$ is either complex or negative for any $\tau\geqslant 0$, and
$\omega_+^2$ is positive, if, and only if,
\begin{equation*}
\begin{aligned}
(\alpha^2(\tau)-c^2(\tau))=-3 s^2\left(\frac{se^{s\tau}}{Y}-\frac{1}{3}\right)
\left(\frac{se^{s\tau}}{Y}-1\right)<0.
\end{aligned}
\end{equation*}
This is the case, if, and only if,
%\begin{equation*}
$\frac{se^{s\tau}}{Y}<\frac{1}{3}$ or %\qquad\mbox{or}\qquad
$\frac{se^{s\tau}}{Y}> 1.$
%\end{equation*}
However, $E_+$ only exists when $x_+(\tau)=\frac{se^{s\tau}}{Y}<1$.
Therefore, a real positive root exists if, and only if,
$x_+(\tau)=\frac{se^{s\tau}}{Y}<\frac{1}{3}$. This
implies that $\tau<\tau^*$.
%or equivalently, $x_+(\tau)<\frac{1}{3}$.
Hence, for $\tau\in[0,\tau^*)$, a real positive root $\omega_+(\tau)$
exists, and is defined explicitly by (\ref{chpt1omegapm}).

If $\tau=\tau^*$, then $\omega_+(\tau)=0$, and if
$\tau>\tau^*$, then either $\omega_+(\tau)$ is not real or $E_+$ does not
exist.
$\hfill\Box$

\subsection{Properties of the functions $h_1(\omega,\tau)$ and
$h_2(\omega,\tau)$ defined in (\ref{h1_h2_defined})}
 \label{app:h1_h2}
\noindent

\begin{lemma}\label{lem:h1_pos}
  Assume that $E_+$ exists and that $\frac{s}{Y}<\frac{1}{3}$.
  \begin{enumerate}
    \item   If $\tau\in[0,\tau^*]$,  then $h_1^2(\omega_+(\tau),\tau)+h_2^2(\omega_+(\tau),\tau)=1$.
    \item
 $h_1(\omega_+(\tau^*),\tau^*)=0.$  If $ \ 0\leqslant\tau<\tau^*$, then
 $h_1(\omega_+(\tau),\tau)>0$.
 % (where $h_1$ is defined in
 % (\ref{h1defined}) and $\omega_+(\tau)$ is defined in
  %(\ref{chpt1omegapm})).
 \item  $h_2(\omega_+(\tau^*),\tau^*)=-1.$ If $ \ 0\leqslant\tau<
   \tau^*$, then $-1 <
  h_2(\omega_+(\tau),\tau)< 1$.
 % (where $h_2$ is defined in
  %(\ref{h2defined}) and $\omega_+(\tau)$ is defined in
 % (\ref{chpt1omegapm})).
\end{enumerate}
\end{lemma}

\noindent
{\bf Proof.} Part 1.  Since $h_1(\omega_+(\tau),\tau)$ is equal to the
right-hand side of (\ref{chpt1sincos:sin}), and
$h_2(\omega_+(\tau),\tau)$ is equal to the right-hand side of
(\ref{chpt1sincos:cos}), where $\omega=\omega_+(\tau)$ satisfies (\ref{omegapower4}), it follows that
 $h_1^2(\omega_+(\tau),\tau)+h_2^2(\omega_+(\tau),\tau)=1$.
\medskip

Part 2. By Theorem~\ref{omega+posicond}, $\omega_+(\tau^*)=0$, and so
$h_1(\omega_+(\tau^*),\tau^*)=0.$ Since for $t\in[0,\tau^*)$,  $x_+(\tau)<\frac{1}{3}$  and  by
Theorem~\ref{omega+posicond}, $\omega_+(\tau)>0$,
 it follows that the denominator of
$h_1(\tau,\omega_+(\tau))$ is always positive.  Hence, $h_1(\omega_+(\tau),\tau)>0$, if, and only if,
  \begin{eqnarray*}
  0&<& s(1-x_+(\tau))+x_+(\tau)(1-2x_+(\tau))-\omega_+^2(\tau)\\
%&=&s+x_+(\tau)-s x_+(\tau)-2x_+^2(\tau)-\omega^2>0.
&=&
s(1-x_+(\tau))+x_+(\tau)(1-2x_+(\tau))\\
&&-\frac{1}{2}\left(-x_+^2(\tau)+\sqrt{
x_+^4(\tau)+4s^2(3x_+(\tau)-1)(x_+(\tau)-1)}\right).
\end{eqnarray*}

This is equivalent to,
$$\frac{1}{2}\sqrt{x_+^4(\tau)+4s^2(3x_+(\tau)-1)(x_+(\tau)-1)}
<s(1-x_+(\tau))+x_+(\tau)(1-\frac{3}{2}x_+(\tau)).$$
Since $x_+(\tau)<\frac{1}{3}$, both sides of the above inequality are
positive. Squaring, both sides yields,
\begin{eqnarray*}
&&\frac{1}{4}\left (x_+^4(\tau)+4s^2(1-3x_+(\tau))(1-x_+(\tau))\right)
\\
&&\quad < s^2(1-x_+(\tau))^2+x_+^2(\tau)(1-\frac{3}{2}x_+(\tau))^2\\
&& \quad \quad +
\quad 2sx_+(\tau)(1-x_+(\tau)(1-\frac{3}{2}x_+(\tau)).
\end{eqnarray*}
But this is equivalent to,
\begin{eqnarray*}
  0&<& 2s^2 x_+(\tau)(1-x_+(\tau))+x_+^2(\tau) (1-2x_+(\tau))(1-x_+(\tau))\\
  && + \quad 2sx_+(\tau)(1-x_+(\tau)(1-\frac{3}{2}x_+(\tau)).
\end{eqnarray*}
This last inequality is satisfied, since by part 2 of
Theorem~\ref{omega+posicond},
$x_+(\tau)<\frac{1}{3}$, and so
 $h_1(\omega_+(\tau),\tau)>0$ for $\tau\in[0,\tau^*)$.
\medskip

Part 3.   That $h_2((\omega_+(\tau^*),\tau^*)=-1$, follows from a straightforward calculation, after substituting $\omega_+(\tau^*)=0$ and
 $x_+(\tau^*)=\frac{1}{3}$ in (\ref{h2defined}). For
 $\tau\in[0,\tau^*)$, the result follows immediately, by parts 1 and 2.
 $\hfill\Box$
\medskip

\subsection{Proof of Theorem~\ref{chpt3intersections}}
 \label{app:chpt3intersections}

\noindent
{\bf Proof.} Assume that $\bar{\tau}\in(0,\tau^*)$.
By Remark~\ref{rem:char_pure_imag},
(\ref{CharE+}) has a  pair of pure
imaginary eigenvalues, if and only if, $\bar{\tau}=\tau_n^j\in(0,\tau^*)$, for
some $n\geqslant 0, \ 0\leqslant j \leqslant j_n$.  Since a necessary condition for a root of (\ref{CharE+})
to exist is that
 (\ref{omegapower4}) holds, and hence
$\omega(\bar{\tau})= \omega_+(\bar{\tau})$ given by
(\ref{chpt1omegapm}),  there can be at most one pair of pure imaginary
roots for each such $\tau=\bar{\tau}$, and hence  if any such roots exist, they are simple, and no other root of
(\ref{CharE+}) is an integer multiple.
%By Theorem~\ref{omega+posicond}, $\omega_+(\tau)> 0$.
% By Lemma \ref{omegataul},
% if there exists an integer
%$n\geqslant 0$ such that $\ \theta(\tau)+2 n\pi$ intersects
%$\tau\omega_+(\tau)$ at some $\tau_{n}\in (0,\tau^*)$, then
%$(\tau_n,\omega_+(\tau_n))$ is a solution of (\ref{chpt1sincos}),
%and therefore the characteristic equation (\ref{CharE+}) at $E_+$ has a
%pair of pure imaginary eigenvalues,
%$\lambda=\pm i\omega_+(\tau_n)$. For any such $\tau$, there is only one
%pair of pure imaginary roots,
%no other root of (\ref{CharE+})
%is an integer multiple of $\pm i \omega_+(\tau_n)$.
%
%Next, consider the transversality requirement for Hopf
%bifurcation, i.e., whether
%$\frac{\mathrm{d}\mathrm{Re}(\lambda(\tau))}{\mathrm{d}\tau}\big|_{\tau=\tau_n^j}\neq 0$.

In Beretta and  Kuang,
\cite{Kuang2002}  (Theorem~4.1, equation (4.1) p.1157),
it is shown that
 $$ \mathrm{sign}\left(\frac{\mathrm{d}}{\mathrm{d} \tau}
 \mathrm{Re}(\lambda(\tau)) \right)\Big|_{\tau=\tau_{n^j}}
=
\mathrm{sign}\left(\frac{\mathrm{d}}{\mathrm{d}\tau}\left[\frac{\tau
\omega_+(\tau)-(\theta(\tau)+2n\pi)}{\omega_+(\tau)}
\right]\right)\Big|_{\tau=\tau_n^j}.$$
After differentiating, and noting that $\omega_+(\tau_n^j)>0$ and
$\tau_n^j\omega_+(\tau_n^j)-(\theta(\tau_n^j)+2n\pi)=0$, it is easy to see that the term on the right has the same sign as
$\left(\frac{\mathrm{d}}{\mathrm{d}\tau}\left[\tau
\omega_+(\tau)-(\theta(\tau)+2n\pi)
\right]\right)\big|_{\tau=\tau_n^j}$.
It follows that transversality holds whenever the graphs
of $\tau \omega_+(\tau)$ and $\theta(\tau)$ have different slopes
at   $\tau_n^j$.
$\hfill\Box$
\medskip

\subsection{Proof of Corollary~\ref{cor:2k+1intersection}}
 \label{app:2k+1intersection}

\noindent
{\bf Proof.}  Parts 1, 2 and 3 follow immediately from
Theorem~\ref{chpt3intersections}, since for $\tau\in[0,\tau^*]$, both curves are continuous,
 by part 3 of Lemma~\ref{lem:h1_pos} and
Remark~\ref{rem:characterization},  $\theta(\tau)\in(0,\pi)$,
and $\tau\omega_+(\tau)\geqslant 0$ with equality  for $\tau=0$ and
$\tau=\tau^*$.  (See Figure~\ref{fig:thetatauomega2} for typical examples.)

Next consider local stability of $E_+$ when it exists, i.e.
$\tau\in[0,\tau_c)$. By part 3(a) of Theorem~\ref{th:gasEi},
 $E_+$ is globally asymptotically stable when $\tau=0$.  By
Lemma~\ref{lem:imaginarycross},
$E_+$ can only change stability for $0\leqslant\tau<\tau_c$, if a pair
of roots cross the imaginary axis. As $\tau$ increases from $0$,  by
part 3, the first possible such crossing occurs for $\tau=\tau_0^1.$
Hence, $E_+$ is locally asymptotically stable for $\tau\in[0,\tau^1_0)$.
When $\tau=\tau_c$,  $E_+$ and $E_1$ coalesce, and by part 2(b) of
Theorem~\ref{th:gasEi}, $E_1$ is globally asymptotically stable for
$\tau>\tau_c$. Since again by
Lemma~\ref{lem:imaginarycross},
$E_+$ can only change stability for $0\leqslant\tau<\tau_c$, if a pair
of roots cross the imaginary axis, and $\tau_0^{j_0}$ is the last
possible such crossing, $E_+$ must be locally asymptotically stable for
$\tau\in(\tau_0^{j_0},\tau_c)$.
$\hfill\Box$
\medskip

\end{document}